\documentclass{amsart}
\usepackage[mathscr]{eucal}
\usepackage[dvipdfm]{graphicx}
\usepackage{amssymb, amscd, color, a4wide}
\theoremstyle{plain}
\newtheorem{theorem}{Theorem}[section]
\newtheorem{corollary}{Corollary}[section]
\newtheorem{lemma}{Lemma}[section]
\newtheorem{proposition}{Proposition}[section]
\newtheorem{remark}{Remark}[section]

\newtheorem{definition}{Definition}[section]

\numberwithin{equation}{section}
\def\Vec#1{\mbox{\boldmath $#1$}}
\newcommand{\C}{\mathbb C}
\newcommand{\diag}{\operatorname{diag}}

\begin{document}
\title[CMC surfaces in $\mathbb{H}^3$]
{Constant mean curvature surfaces in hyperbolic 3-space via loop groups}
\author[J. F.~Dorfmeister]{Josef F. Dorfmeister}
\dedicatory{Dedicated to the memory of Hongyou Wu}
\address{Fakult\"at f\"ur Mathematik, 
TU-M\"unchen, 
Boltzmann str. 3,
D-85747, 
Garching, 
Germany}
\email{dorfm@ma.tum.de}
\author[J.~Inoguchi]{Jun-ichi Inoguchi}
\address{Department of Mathematical Sciences, 
Faculty of Science,
Yamagata University, 
Yamagata, 990--8560, Japan}
\email{inoguchi@sci.kj.yamagata-u.ac.jp}
\thanks{The second named author is partially supported by Kakenhi
18540068, 21546067}
\author[S.-P.~Kobayashi]{Shimpei Kobayashi}
\address{
Graduate School of Science and Technology,
Hirosaki University,
Hirosaki, 036-8561, Japan}
\email{shimpei@cc.hirosaki-u.ac.jp}
\thanks{The third named author is partially supported by Kakenhi
20740045}
\date{\today}
\begin{abstract}
 In hyperbolic $3$-space $\mathbb{H}^3$ surfaces of constant 
 mean curvature $H$ come in three types, corresponding to the cases 
 $0 \leq H < 1$, $H = 1$, $H > 1$.
 Via the Lawson correspondence the latter two cases
 correspond to constant mean curvature surfaces in Euclidean
 $3$-space $\mathbb{E}^3$ with $H=0$ and $H \neq 0$, respectively.
 These surface classes have been investigated intensively
 in the literature. For the case $0 \leq H < 1$ there is 
 no Lawson correspondence in Euclidean space and there
 are relatively few publications. Examples have been difficult 
 to construct. In this paper we present a generalized Weierstrass type 
 representation for surfaces of constant mean curvature in
 $\mathbb{H}^3$ with particular emphasis on
 the case of mean curvature $0\leq  H < 1$.
 In particular, the generalized Weierstrass type 
 representation presented in this paper enables us to 
 construct simultaneously minimal surfaces 
 ($H=0$) and non-minimal constant mean curvature surfaces ($0<H<1$).
\end{abstract}
\maketitle
\section*{Introduction}
 Harmonic maps into (semi-simple) Riemannian symmetric spaces
 can be studied by infinite dimensional Lie group theory (\textit{loop
 groups}).
 The dressing action of the loop group enables one to construct
 ``non-trivial solutions'' from a trivial solution and, meromorphically, 
 new solutions from old/known ones. This is based on a
 Weierstrass type representation for harmonic maps which
  was established by the first named author of the present paper, 
 F.~Pedit and H.~Wu \cite{DPW} in terms of loop group decompositions.

 As an application of the loop group method, constant mean curvature
 surfaces in Euclidean  $3$-space $\mathbb{E}^3$, pseudo-spherical
 surfaces in $\mathbb{E}^3$ and their indefinite analogues (spacelike or
 timelike surfaces in Minkowski $3$-space) are studied extensively
 \cite{DIT}, \cite{I}. The loop group method to construct such
 \textit{integrable surfaces} is frequently called \textit{DPW method} or 
 \textit{generalized Weierstrass type representation} for surface
 geometry.

 In Euclidean geometry, the starting point of the DPW method is that
 integrable surfaces are naturally associated with harmonic maps into the
 $2$-sphere $\mathbb{S}^2$  with respect to appropriate conformal
 structures. More precisely, a surface in $\mathbb{E}^3$ is of constant
 mean curvature if and only if its Gauss map is harmonic. This
 characterization of constancy of mean curvature is referred to as the
 \textit{Ruh-Vilms property}. On the other hand, let $f:M\to
 \mathbb{E}^3$ be a surface in $\mathbb{E}^3$ with negative or positive
 Gaussian curvature. Then the Gaussian curvature is constant if
 and only if its Gauss map is harmonic with respect to the conformal
 structure determined by the second fundamental form. 
 This is another instance of the Ruh-Vilms property.
 Analogously, in Minkowski 3-space the corresponding integrable
 surfaces are associated with harmonic maps into hyperbolic $2$-space
 $\mathbb{H}^2$ or de Sitter $2$-space $\mathbb{S}^{1, 1}$.

 In M{\"o}bius geometry, a conformally immersed surface in
 M{\"o}bius $3$-space $\mathcal{M}^3$ is M{\"o}bius minimal (Willmore)
 if and only if its conformal Gauss map (central sphere congruence)
 is a harmonic map into the
 moduli space of oriented $2$-spheres. Note that via the projective
 light-cone model of the M{\"obius} $3$-space, the
 moduli space of all oriented $2$-spheres (including point spheres)
 is identified with the de Sitter $4$-space $\mathbb{S}^{1, 3}$, a
 semi-simple Lorentzian symmetric space \cite{Bryant:Willmore}.

 F.E.~Burstall and U.~Hertrich-Jeromin \cite{BHJ} proved the Ruh-Vilms
 property for certain surfaces in Lie sphere geometry as well as
 projective differential geometry. More precisely,  let $f:M \to
 \mathcal{M}^3$ be a surface in M{\"o}bius space. Then $M$ is said to be
 a \textit{Lie minimal surface} if $M$ is a critical point of the
 Lie-area functional. Note that the Lie-area functional is invariant 
 under the M{\"o}bius group  as well as the Lie sphere transformation group.
 On the other hand, for every surface in real projective $3$-space
 $\mathbb{R}P^3$, a projectively invariant area functional (projective
 area) has been introduced. A critical point of the projective area is
 called a \textit{projective minimal surface}.
 Both the Lie-minimality and projective minimality are characterized by
 the harmonicity of  appropriate Gauss maps taking values in a certain 
 non-compact Grassmann manifold equipped with an invariant indefinite
 semi-Riemannian metric.

 In addition, E.~Musso and L.~Nicolodi \cite{MN} showed that in Laguerre
 geometry, Laguerre-minimality of surfaces is characterized by the
 harmonicity of the Laguerre Gauss map (\textit{middle sphere
 congruence}). 
 Thus all these  surfaces classes are defined by the Ruh-Vilms property.

 Here we should remark that all the target manifolds of \textit{harmonic
 Gauss maps} in Lie-sphere, M\"obius, projective and Laguerre
 geometries are semi-simple semi-Riemannian 
 \textit{symmetric spaces} equipped with invariant indefinite metrics.

 It would be interesting to characterize all surface geometries to
 which the DPW method is applicable. So far only partial results are
 known: See \textit{e.g.}, the survey \cite{D} by the first named
 author and the recent classification, carried out by the third named
 author, of surfaces in 3-dimensional space forms
 which are \textit{real forms} of complex constant mean curvature
 surfaces  \cite{Kobayashi}.
 
 In this article, we give a \textit{new} example of such a surface
 geometry. More precisely we shall show that the DPW method is
 also applicable to surfaces in hyperbolic $3$-space
 $\mathbb{H}^3$ of constant mean curvature $0 \leq H <1$.

 The geometry of a CMC surface in $\mathbb{H}^3$ depends on
 the range of the mean curvature.
 In case $H>1$, there exist compact CMC surfaces  in $\mathbb{H}^3$.
 It is well known that there exist natural, locally bijective correspondences 
 (so-called \textit{Lawson correspondences}) between CMC surfaces in 
 3-dimensional space forms.
 In particular, CMC surfaces in $\mathbb{H}^3$ with $H\geq 1$ have 
 corresponding CMC or minimal surfaces (\textit{Lawson correspondents}) 
 in Euclidean 3-space or the $3$-sphere. Based on the Lawson correspondences,
 M.~Kilian, W.~Rossman, N.~Schmitt and the third named author \cite{KKRS}
 studied CMC surfaces in $\mathbb{H}^3$ with $H>1$ via the DPW method. 
 In case $H=1$, the Lawson correspondent in $\mathbb E^3$ 
 is a minimal surface, and  R.~Bryant \cite{Br} 
 gave the Weierstrass type representation  
 in terms of a  holomorphic differential equation.

 On the other hand, by the maximum principle, there are no compact 
 CMC surfaces in $\mathbb{H}^3$ (without boundary) such that 
 $0\leq H \leq1$. Moreover, 
 CMC surfaces in $\mathbb{H}^3$ with $0\leq H<1$ have 
 neither Lawson correspondents in Euclidean 3-space $\mathbb{E}^3$ nor in
 the $3$-sphere  $\mathbb{S}^3$. 
 This means that CMC surface geometry with $0\leq H<1$
 in hyperbolic space has special and unusual features and hence 
 it is of interest. 

 There is another motivation for the study of CMC surfaces with $0<H<1$. 
 D.~Chopp and J.~Velling \cite{CV} found numerical evidence that 
 $\mathbb{H}^3$ can be foliated by CMC discs with $0<H<1$ that share a common 
 Jordan curve boundary in the ideal boundary of $\mathbb{H}^3$. 

 In contrast to the case $H\geq 1$, there are relatively few papers on 
 CMC surfaces in $\mathbb{H}^3$ with  $0\leq H<1$.
 Initial studies on the existence and regularity of CMC surfaces in 
 $\mathbb{H}^3$ with $0\leq H<1$ can be found in  
 M.~Anderson \cite{Anderson}, 
 K.~Uhlenbeck \cite{Uhlenbeck}, 
 K.~Polthier \cite{Polthier}, 
 Y.~Tonegawa \cite{Tonegawa}
 and J.~Velling \cite{Velling}.

 In \cite{BaBo}, M.~Babich and A.~Bobenko considered minimal surfaces 
 of finite type in $\mathbb H^3$ and described them explicitly in terms 
 of theta functions.

 One of the remarkable points of the present study is that, --as in the
 already known cases--
 we associate a certain kind of \textit{Gauss map} to constant mean
 curvature surfaces in $\mathbb H^3$. And the constancy of the mean
 curvature is equivalent to the harmonicity of the Gauss map. 
 Thus the Ruh-Vilms property holds.
 However, unlike all surface geometries discussed by the DPW method
 so far, our target space is \textit{not} a symmetric space but a 
 real $5$-dimensional $4$-symmetric space.

 In general, the loop group approach can not be applied to harmonic maps
 into non-symmetric homogeneous spaces. 
 Thus it is not obvious how  to apply the loop group method to our case.  
 More precisely, the loop group method has two main 
 ingredients. One is a zero curvature representation and the other one 
 is a loop group decomposition.
 The zero curvature representation is equivalent with the existence of a 
 flat connections and this representation enables us to use 
 loop groups. A loop group decomposition is at the heart 
 of the DPW method and recovers surfaces from holomorphic potentials.

 In this paper, the Gauss map of CMC surfaces in $\mathbb H^3$
 has as target space a $4$-symmetric 
 space. Thus the DPW method is, a priori, not applicable.
 However it turns out that these Gauss maps are Legendre maps. 
 In the situation under consideration this suffices to obtain a 
 zero curvature representation for all Gauss maps of CMC
 surfaces in $\mathbb H^3$ (and in particular for mean curvature 
 $H$ satisfying $0 \leq H <1$).

 Applying the DPW method in this setting, we obtain a potential
 in a $4$-graded Lie algebra. However, not every potential in this 
 $4$-graded Lie algebra will admit the application of the DPW method, 
 since the admissibility condition
 will not be satisfied in general \cite{BP}.

 It is a fortunate accident that the Maurer-Cartan form of 
 the moving frame of every CMC immersion into $\mathbb H^3$
 with $0 \leq H <1$ is also contained in same real form of 
 $\Lambda \mathrm{SL}_2 \C_{\sigma} \times \Lambda \mathrm{SL}_2 \C_{\sigma}$
 and this is part of some ``symmetric space setting'' 
 (which actually is not induced from any finite dimensional symmetric space).
 This enables us to apply the DPW method via standard way and therefore 
 two main featured of the loop group method can be established and 
 the loop group method is applicable.

 The reason for the applicability of the DPW method in our case
 also can be explained by using the Lawson correspondence which 
 transforms our original harmonic Gauss map into a {\it primitive} 
 map, to which the DPW method is applicable.
 By the Lawson correspondence for a CMC surface 
 with mean curvature $0 \leq H <1$ in $\mathbb H^3$, 
 the CMC surface can be considered 
 as a minimal surface in $\mathbb H^3 (c)$ for some $c <0$, where  
 $\mathbb H^3 (c)$ denotes the hyperbolic $3$-space of sectional 
 curvature $c$, \cite{AA:CMCH3}.
 At this point it is important to note that 
 the Gauss map of a CMC surface in $\mathbb H^3(c)$ is primitive
 if and only if the surface is minimal. 
 Finally we recall that primitive maps 
 into $k$-symmetric ($k>2$) spaces
 form a special class of harmonic maps and it is known that the DPW method 
 is applicable to primitive maps \cite{DMPW}. Moreover, we interpret
 our original Gauss map of a CMC surface in $\mathbb H^3$,
 which is only a Legendre harmonic map, as a primitive 
 Gauss map of some CMC surface in $\mathbb H^3(c)$, 
 and thus the loop group method becomes applicable.

 We would like to point out that the DPW method for CMC surfaces 
 in $\mathbb{H}^3$ presented in \cite{KKRS} is based on the 
 Lawson correspondence, so the authors of 
 \cite{KKRS} can treat CMC surfaces in $\mathbb{H}^3$ with mean curvature $H>1$.
 In the present study, we give a unified approach to both
 cases: $H>1$ and $0 \leq H  <1$. The case $H =1$ has been studied
 extensively already, \textit{e.g.}, \cite{Br}, \cite{UY},
 and will not be considered in  this paper.

 This paper is organized as follows.
 After establishing the requisite knowledge on harmonic maps into normal 
 semi-Riemannian homogeneous spaces and homogeneous geometry of
 hyperbolic space in sections \ref{sc:Harmonic}--\ref{sc:Hyperbolic},
 we shall devote  sections \ref{sc:SurfacesH3}--\ref{sc:2by2} to 
 surface geometry in $\mathbb{H}^3$ in terms of 
 $\mathrm{SL}_{2}\mathbb{C}$-valued functions. 
 In section \ref{sc:CMCH>1}, 
 we shall  give a loop group formulation of CMC surfaces with $H>1$. 
 This formulation is different from the one used in \cite{KKRS}. 
 We shall clarify the Lawson correspondences between CMC surfaces 
 in $\mathbb{H}^3$ with $H>1$ and  CMC surfaces in $\mathbb{E}^3$ 
 in terms of loop groups.

 In the next section \ref{sc:CMCH<1}, we shall give a loop group formulation
 for CMC surfaces with  $0\leq H<1$. The two cases discussed in section 
 \ref{sc:CMCH>1} and section  \ref{sc:CMCH<1} are distinguished 
 by the automorphism of the Kac-Moody Lie algebra which characterizes
 the Lax equations of corresponding types of CMC surfaces.  
 
 The key tool of the present study is a ``contact geometric characterization'' 
 of CMC surfaces. 
 For every (oriented) surface $f:M\to \mathbb{H}^3$, there exists 
 a smooth map $F$ into the unit tangent sphere bundle $\mathrm{U}\mathbb{H}^3$ 
 of $\mathbb{H}^3$. 
 The map $F$ is referred to as  the \textit{Gauss map} of $M$.
 One can see that the Gauss map $F$ of a surface $M$ is a Legendre map
 in the sense of V.I.~Arnold \cite{Arnold}, that is, 
 it is tangent to the canonical contact structure. 
 Moreover, the constancy of the mean curvature is equivalent 
 to the harmonicity of that Gauss map. 
 Thus the Ruh-Vilms property holds.
%
  
 In section \ref{sc:4symUnittangent}, we shall characterize harmonic Gauss maps
 in a way which is different from the Ruh-Vilms property.
 The Legendre property of the Gauss map will be characterized in terms of 
 $4$-symmetric structure of $\mathrm{U}\mathbb{H}^3$. 
 These characterizations in terms of contact geometry and 
 $4$-symmetric structure yield a zero-curvature representation for
 Legendre harmonic maps.  Based on these results, in sections 
 \ref{sc:Potential}--\ref{sc:DPW}, 
 we shall give a DPW method for Legendre harmonic maps 
 (and hence for CMC surfaces with $0\leq H<1$). 
 In the final section, we shall exhibit some examples of CMC surfaces 
 with $0\leq H<1$ via the DPW method established in this paper.

 In the surface geometry of $\mathbb{H}^3$, several notions of 
 \textit{Gauss map}  have been introduced.  
 For the convenience of the reader, we collect in the appendix 
 several notions of Gauss map that have been used by 
 different authors. We shall explain how these other 
 Gauss maps can be derived  from our Gauss map.  
 As a side result we obtain that only the Gauss map considered in 
 this paper is suitable for a DPW method of CMC surfaces in $\mathbb H^3$.

\textbf{Acknowledgements:} We would like to thank Idrisse Khemar 
 for helpful comments on a preliminary 
 version of this paper. Part of this
 work was carried out during the workshop ``Surface Theory: Research in
 Pairs'' at Kloster Sch\"{o}ntal, March 2008, funded by DFG Grant
 DO776. This work was started when the second named author visited the
 University of Kansas in 2000. He would like to express his sincere thanks
 to the Department of Mathematics. 
 Some of the results of this article were 
 reported at  the workshop ``Progress in Surface Theory'' held at
 Mathematisches Forschungsinstitut Oberwolfach, May,  2010.
\section{Harmonic maps into normal semi-Riemannian homogeneous spaces}
\label{sc:Harmonic}
\subsection{}
 Let $G/H$ be a reductive homogeneous space with semi-simple Lie group
 $G$. We equip $G/H$ with a $G$-invariant semi-Riemannian metric which is
 derived from (a constant multiple of) the Killing form of $G$.
 Assume that the Lie algebra $\mathfrak{h}$ of $H$ is non-degenerate with 
 respect to the induced scalar product.  
 Then the orthogonal complement $\mathfrak{p}$ of $\mathfrak{h}$ is 
 non-degenerate and can be identified with the tangent space of $G/H$ 
 at the origin $o=H$.  
 The Lie algebra $\mathfrak{g}$ is decomposed into the direct sum:
$$
\mathfrak{g}=\mathfrak{h}\oplus \mathfrak{p}
$$
 of linear subspaces.  
 The resulting homogeneous semi-Riemannian manifold is a 
 \textit{normal semi-Riemannian homogeneous space} \cite{CE}.

\subsection{}
A smooth map $\psi:M \to N$ of a Riemann surface $M$ into a semi-Riemannian 
manifold $N$ is said to be a \textit{harmonic map} if its tension field  
$\mathrm{tr}\>(\nabla \mathrm{d}\psi)$ vanishes \cite{Ur}. 
When the target space $N$ is a normal semi-Riemannian homogeneous space $G/H$,
the harmonic map equation for $\psi$ has a particularly simple form.

Now let $\psi:\mathbb{D} \to G/H$ be a smooth map from a simply connected
 domain $\mathbb D \subset \C $ into a normal semi-Riemannian homogeneous space. 
Take a frame $\Psi:\mathbb{D} \to G$ of
 $\psi$ and put $\alpha:=\Psi^{-1}\mathrm{d}\Psi$. Then we have the
 identity (\textit{Maurer-Cartan equation}):
$$
\mathrm{d}\alpha+\frac{1}{2}[\alpha \wedge \alpha]=0.
$$
 Decompose $\alpha$ along the Lie algebra decomposition
 $\mathfrak{g}=\mathfrak{h}\oplus\mathfrak{p}$ as
$$
\alpha=\alpha_{\mathfrak{h}}+\alpha_{\mathfrak{p}},
\ \
\alpha_{\mathfrak{h}}\in
\mathfrak{h},
\ \
\alpha_{\mathfrak{p}}
\in
\mathfrak{p}.
$$
We decompose $\alpha_{\mathfrak{p}}$ 
with respect to the conformal
structure of $\mathbb{D}$ as
$$
\alpha_{\mathfrak p}=
\alpha_{\mathfrak p}^{\prime}+
\alpha_{\mathfrak p}^{\prime \prime}.
$$
Here $\alpha_{\mathfrak{p}}^{\prime}$ and 
$\alpha_{\mathfrak p}^{\prime \prime}$
are the $(1,0)$ and $(0,1)$ part of 
$\alpha_{\mathfrak p}$, respectively.

The harmonicity of $\psi$ is equivalent to
\begin{equation}\label{harmonicity}
\mathrm{d}(*\alpha_{\mathfrak p})+
[\alpha \wedge *\alpha_{\mathfrak p}]
=0.
\end{equation}
Here $*$ denotes the Hodge star operator of $\mathbb{D}$.
The Maurer-Cartan equation is split into its $\mathfrak{h}$-component and 
$\mathfrak{p}$-component:
\begin{equation}\label{MC-h}
\mathrm{d}\alpha_{\mathfrak{h}}+\frac{1}{2}
[\alpha_{\mathfrak{h}}\wedge \alpha_{\mathfrak{h}}]
+[\alpha_{\mathfrak{p}}^{\prime}\wedge 
\alpha_{\mathfrak{p}}^{\prime\prime}]_{\mathfrak{h}}=0,
\end{equation}
\begin{equation}\label{MC-p}
\mathrm{d}\alpha
_{\mathfrak{p}}^{\prime}
+
[
\alpha_{\mathfrak{h}}
\wedge 
\alpha_{\mathfrak{p}}^{\prime}
]
+
\mathrm{d}\alpha_{\mathfrak{p}}^{\prime\prime}
+
[
\alpha_{\mathfrak{h}}
\wedge 
\alpha_{\mathfrak{p}}^{\prime\prime}
]
+
[
\alpha_{\mathfrak{p}}^{\prime}
\wedge 
\alpha_{\mathfrak{p}}^{\prime\prime}
]_{\mathfrak{p}}
=0.
\end{equation}
Hence for a harmonic map $\psi:\mathbb{D}\to G/H$ with a framing $\Psi$, 
the pull-back $1$-form $\alpha=\Psi^{-1}\mathrm{d}\Psi$ satisfies
(\ref{harmonicity}), (\ref{MC-h}) and (\ref{MC-p}). 
Combining (\ref{harmonicity}) and (\ref{MC-p}), we have
\begin{equation}\label{harm+MC}
\mathrm{d}\alpha^{\prime}_{\mathfrak p}+
[\alpha_{\mathfrak h} \wedge 
\alpha_{\mathfrak p}^{\prime}]=
-\frac{1}{2}[\alpha_{\mathfrak p}^{\prime}
\wedge 
\alpha_{\mathfrak p}^{\prime \prime}
]_{\mathfrak p}.
\end{equation}
One can easily check that 
the harmonic map equation for $\psi$ combined with
the Maurer-Cartan equation is equivalent to the system (\ref{MC-h}) and (\ref{harm+MC}).

Assume that
\begin{equation}\label{admissible}
[\alpha_{\mathfrak p}^{\prime}
\wedge 
\alpha_{\mathfrak p}^{\prime \prime}
]_{\mathfrak p}=0.
\end{equation}
Then the harmonic map equation together with the Maurer-Cartan equation
is reduced to the system of equations:
$$
\mathrm{d}\alpha_{\mathfrak p}^{\prime}+
[\alpha_{\mathfrak h} \wedge 
\alpha_{\mathfrak p}^{\prime}]=0,
$$
$$
\mathrm{d}\alpha_{\mathfrak h}+
\frac{1}{2}
[\alpha_{\mathfrak h} \wedge
\alpha_{\mathfrak h}]
+[\alpha_{\mathfrak p}^{\prime}
\wedge
\alpha_{\mathfrak p}^{\prime \prime}
]=0.
$$
This system of equations is equivalent to the
following \textit{zero-curvature representation}:
$$
\mathrm{d}\alpha_{\lambda}+
\frac{1}{2}[\alpha_{\lambda}
\wedge \alpha_{\lambda}]
=0,
$$
where $\alpha_{\lambda}:=
\alpha_{\mathfrak h}
+\lambda^{-1}\alpha_{\mathfrak p}^{\prime}+
\lambda\> \alpha_{\mathfrak p}^{\prime \prime}$
with $\lambda \in S^1$. 

\begin{proposition}
 Let $\mathbb D$ be a region in  $\C$ and 
 $\psi:\mathbb{D} \to G/H$ a harmonic map
 which satisfies the admissibile condition {\rm (\ref{admissible})}.
 Then the loop of connections $\mathrm{d}+\alpha_{\lambda}$ is flat for
 all $\lambda$. Namely:
\begin{equation}\label{eq:flatconnections}
\mathrm{d}\alpha_{\lambda}+\frac{1}{2}
[\alpha_{\lambda}\wedge
\alpha_{\lambda}]=0
\end{equation}
 for all $\lambda$.
 Conversely assume that $\mathbb{D}$ is simply connected. Let
 $\alpha_{\lambda}=\alpha_{\mathfrak h}+\lambda^{-1} \alpha_{\mathfrak
 p}^{\prime}+\lambda\alpha_{\mathfrak  p}^{\prime\prime}$ be an
 $S^1$-family of $\mathfrak{g}$-valued $1$-forms 
 satisfying \eqref{eq:flatconnections} for all $\lambda \in S^1$.
%
Then there exists a $1$-parameter family of maps
$\Psi_{\lambda} :\mathbb{D} \to G$ such that
$$
\Psi_{\lambda}^{-1}{\mathrm d}\Psi_{\lambda}=\alpha_{\lambda}
\;\;\;\mbox{and}\;\;\;
\psi_{\lambda}=\Psi_{\lambda}\;{\mbox{\rm mod}}\; H:\mathbb{D} \to G/H
$$
is harmonic for all $\lambda$.
\end{proposition}

When the target space $G/H$ is a semi-Riemannian symmetric space, then the
admissible condition 
is fulfilled automatically for any $\psi$, since
$[\mathfrak{p},\mathfrak{p}]\subset \mathfrak{h}$.

\subsection{Primitive maps}
Let $G$ be a semi-simple Lie group 
with automorphism $\tau$ of order $k> 2$.
A normal semi-Riemannian 
homogeneous space $G/H$ is said to be a 
(regular)  \textit{semi-Riemannian $k$-symmetric space}
if $G_{\tau}^{\circ}\subset H \subset G_{\tau}$.
Here $G_{\tau}$ is the Lie subgroup of all fixed points of 
$\tau$ and 
$G_{\tau}^\circ$ 
the identity component of it.

We denote the induced Lie algebra automorphism of 
$\mathfrak{g}$ by the same letter $\tau$. 
Now we have the eigenspace decomposition of 
the complexified Lie algebra 
$\mathfrak{g}^{\mathbb{C}}$;
$$
\mathfrak{g}^{\mathbb{C}}=\sum_{j\in \mathbb{Z}_k}
\mathfrak{g}_{j}^{\mathbb{C}}.
$$ 
Here $\mathfrak{g}_j^{\mathbb{C}}$ is the eigenspace of $\tau$ with
eigenvalue $\zeta^j$. Here $\zeta$ is the primitive $k$-th root of unity.
In particular, $\mathfrak{g}_0^{\mathbb{C}}=\mathfrak{h}$ and $\overline{\mathfrak{g}_{-1}^{\mathbb{C}}}=\mathfrak{g}_1^{\mathbb{C}}$.
\begin{definition}[\cite{BP}]
{\rm
Let $\psi:M\to G/H$ be a smooth map of a Riemann surface into a regular 
semi-Riemannian $k$-symmetric space. Then $\psi$ is said to be a
\textit{primitive map} if any frame $\Psi$ has
$\alpha_{\mathfrak{p}}^{\prime}$  taking value in $\mathfrak{g}_{-1}^{\mathbb{C}}$.
}
\end{definition} 
One can see that every primitive map satisfies the admissible condition 
(\ref{admissible}). In fact,
$$
[\alpha_{\mathfrak{p}}^{\prime}
\wedge
\alpha_{\mathfrak{p}}^{\prime\prime}]
=[\alpha_{-1}\wedge\alpha_{1}]
\in \mathfrak{g}_0^{\mathbb{C}}=
\mathfrak{h}.
$$ 
For more information on primitive maps, we refer to \cite{BP} and
\cite{DMPW}.

\section{Hyperbolic space}\label{sc:Hyperbolic}
\subsection{}
 Let $\mathbb{E}^{1,n}$ be the $n+1$-dimensional Minkowski space
 with scalar product.
$$
\langle \Vec{x},\Vec{y}\rangle=-x_{0}y_{0}+x_{1}y_{1}+\cdots+x_{n}y_{n}.
$$
 We denote by $\Vec{e}_0=(1, 0, \dots, 0), \Vec{e}_1=(0, 1, \dots, 0), \dots,
 \Vec{e}_n=(0, 0, \dots, 1)$ the basis of $\mathbb{E}^{1,n}$.
 We will identify Euclidean $n$-space $\mathbb E^{n} = \mathbb
 E^{0,n} = \{ \Vec{x} \in \mathbb{E}^{1,n} \;| \;x_0 =0\}$. 
 Let $\epsilon$ denote the \textit{signature matrix} defined by
 $\epsilon=\diag (-1,1,\cdots,1)$. Then we define the
 \textit{Lorentz group} $\mathrm{O}_{1,n}$ by
$$
 \mathrm{O}_{1,n}=\{A\in\mathrm{GL}_{n+1}\mathbb{R}
 \
 \vert
 \
 A^{t}\epsilon A=\epsilon\}.
$$
 Note that $\mathrm{O}_{1,n}$ acts isometrically on $\mathbb{E}^{1,n}$,
 and has four connected components. We denote by $\mathrm{SO}^{+}_{1,n}$
 the identity component of $\mathrm{O}_{1,n}$. The Lie algebra
 $\mathfrak{so}_{1,n}$ of $\mathrm{SO}^+_{1,n}$ is given by
$$
 \mathfrak{so}_{1,n} =
 \left
 \{
 X=
 \left(
 \begin{array}{cc}
 0 & \Vec{x}^{t}\\
 \Vec{x} & b
 \end{array}
 \right)
 \
 \biggr 
 \vert
 \
 b\in
 \mathfrak{o}_{n},
 \
 \Vec{x}\in \mathbb{R}^{n}
 \
 \right\}.
$$
 We equip the Lie algebra $\mathfrak{so}_{1,n}$ with an invariant
 scalar product $\langle \cdot,\cdot\rangle$. 
\begin{equation}\label{scalarproduct}
 \langle X,Y\rangle=\frac{1}{2}\mathrm{tr}\>(XY),\ \
 X,Y \in \mathfrak{so}_{1,n}.
\end{equation}
 Since $\mathrm{O}_{1,n}$ is non-compact, this scalar product
 is indefinite. More precisely it has signature $(n(n-1)/2,n)$.
\subsection{}
 The Lie group $\mathrm{SO}^{+}_{1,n}$ acts transitively and 
 isometrically on the hyperbolic $n$-space
$$
 \mathbb{H}^{n}=\{\Vec{x}\in \mathbb{E}^{1,n}
 \
 \vert
 \ \langle \Vec{x},\Vec{x}\rangle=-1,\
 x_{0}>0\
 \}.
$$
 The isotropy subgroup of 
 $\mathrm{SO}^{+}_{1,n}$ at $\Vec{e}_{0}$ is
$$
\left\{
\left(
\begin{array}{cc}
1 & 0
\\
0 & a 
\end{array}
\right)
\
\biggr
\vert
\
a \in \mathrm{SO}_{n}
\right\}.
$$
Hence $\mathbb{H}^{n} \cong \mathrm{SO}^{+}_{1,n}/ \mathrm{SO}_{n}$.
The tangent space $T_{\Vec{e}_0}\mathbb{H}^n$ of $\mathbb{H}^n$ at $\Vec{e}_0$ 
is identified with the following linear subspace of $\mathfrak{so}_{1,n}$;
$$
\left\{
\left(
\begin{array}{cc}
0 & \Vec{x}^{t}\\
\Vec{x} & 0
\end{array}
\right)
\
\biggr 
\vert
\
\Vec{x}\in \mathbb{R}^{n}
\
\right\}.
$$
 On $T_{\Vec{e}_0}\mathbb{H}^{n}$, the scalar product
 (\ref{scalarproduct}) is positive
 definite. Moreover the Riemannian metric on
 $\mathrm{SO}^{+}_{1,n}/\mathrm{SO}_{n}$ induced from
 (\ref{scalarproduct}) is of constant curvature $-1$.

 Next we define an involution $\sigma_{H}$ of $\mathrm{SO}^{+}_{1,n}$ by
 $\sigma_{H}=\mathrm{Ad}(\epsilon)$. Then
 $(\mathrm{SO}^{+}_{1,n},\sigma_{H})$ is a symmetric pair.  
\begin{remark}
{\rm
The Killing form $\varphi$ of $\mathfrak{so}_{1,n}$
is 
$$
\varphi(X,Y)=-(n-1)\mathrm{tr}\>(XY),
\ \
X,Y 
\in \mathfrak{so}_{1,n}.
$$
 We equip the tangent space $T_{\Vec{e}_0}\mathbb{H}^{n} \subset
 \mathfrak{so}_{1,n}$ with the inner product $\langle\cdot,\cdot\rangle$
 given  by
$$
\langle
 X,Y\rangle=-\frac{1}{2(n-1)}\varphi(X,Y)=\frac{1}{2}\mathrm{tr}\>(XY).
$$
}
\end{remark}

\subsection{}\label{sc:unit tangent sphere bundle}
Let us denote by $\mathrm{U}\mathbb{H}^n$ the 
\textit{unit tangent sphere bundle} of $\mathbb{H}^n$.
Namely, $\mathrm{U}\mathbb{H}^n$ is the manifold of all unit tangent 
vectors of $\mathbb{H}^n$.
Then $\mathrm{U}\mathbb{H}^n$ is identified with the submanifold
$$
\left\{
(\Vec{x},\Vec{v})\
\vert
\
\langle \Vec{x},\Vec{x}\rangle=-1,\
\langle \Vec{v},\Vec{v}\rangle=1,\
\langle \Vec{x},\Vec{v}\rangle=0,\
x_{0}>0
\right\}
$$
 of $\mathbb{E}^{1,n}\times \mathbb{E}^{1,n}$. The tangent space
 $T_{(\Vec{x},\Vec{v})}\mathrm{U}\mathbb{H}^{n}$ at a point
 $(\Vec{x},\Vec{v})$ is expressed as
$$
T_{(\Vec{x},\Vec{v})}\mathrm{U}\mathbb{H}^{n}=
\{(X,V)\in \mathbb{E}^{1,n}\times 
\mathbb{E}^{1,n}
\
\vert
\
\langle \Vec{x},X\rangle=0,\
\langle \Vec{v},V\rangle=0,\ 
\langle\Vec{x}, V \rangle + \langle \Vec{v},
X\rangle =0 
\}.
$$
 Define a $1$-form $\omega$ on $\mathrm{U}\mathbb{H}^{n}$ by
$$
\omega_{(\Vec{x},\Vec{v})}(X,V)=\langle X,\Vec{v}\rangle=-\langle
\Vec{x},V\rangle.
$$
 Then one can see that $\omega$ is a \textit{contact form} on
 $\mathrm{U}\mathbb{H}^{n}$, \textit{i.e.},
 $(\mathrm{d}\omega)^{n-1}\wedge \omega\not=0$. The distribution 
$$
\mathcal{D}_{(\Vec{x},\Vec{v})}:=
\{(X,V)\in T_{(\Vec{x},\Vec{v})}\mathrm{U}
\mathbb{H}^n
\
\vert
\
\omega_{(\Vec{x},\Vec{v})}(X,V)=0
\}
$$
 is called the \textit{canonical contact structure} of
 $\mathrm{U}\mathbb{H}^{n}$, see section \ref{sc:contactmfd}.

 The Lorentz group $\mathrm{SO}^+_{1,n}$ acts on $\mathrm{U}\mathbb H^n$ via
 $A \cdot (\Vec{x}, \Vec{v}) = (A \Vec{x},A\Vec{v})$. It is easy to see that
 under this action the unit tangent sphere bundle
 $\mathrm{U}\mathbb{H}^{n}$ is a homogeneous space 
 of $\mathrm{SO}^{+}_{1,n}$. The isotropy subgroup at
 $(\Vec{e}_{0},\Vec{e}_1)$ is
$$
\left\{
\left(
\begin{array}{cc}
\mathbf{1} & 0\\
0 & b 
\end{array}
\right)
\
\biggr
\vert
\
b \in \mathrm{SO}_{n-1}
\right\}.
$$
 Hence $\mathrm{U}\mathbb{H}^n \cong
 \mathrm{SO}^{+}_{1,n}/\mathrm{SO}_{n-1}$. 
 The Lie algebra of this isotropy subgroup is
$$
\left\{
\left(
\begin{array}{cc}
0 & 0 \\
0 & t
\end{array}
\right)
\
\biggr
\vert
\
t\in \mathfrak{so}_{n-1}
\right\}.
$$
 The tangent space $T_{(\Vec{e}_0,\Vec{e}_1)}(\mathrm{U}\mathbb{H}^{n})$
 is identified with
$$
\left
\{
X=
\left(
\begin{array}{lll}
0 & x_1 &\Vec{x}^t \\ 
x_1 &0 & \Vec{y}^t \\ 
\Vec{x} &- \Vec{y} &\Vec{0}
\end{array}
\right)
\
\biggr 
\vert
\
x_1\in \mathbb R, \Vec{x}, \Vec{y} \in \mathbb R^{n-1}
\
\right\}\subset \mathfrak{so}_{1,n}.
$$
 The semi-Riemannian metric induced on the homogeneous space
 $\mathrm{U}\mathbb{H}^n \cong \mathrm{SO}^{+}_{1,n}/\mathrm{SO}_{n-1}$
 via the scalar product (\ref{scalarproduct}) on $\mathfrak{so}_{1,n}$
 has signature $(n-1,n)$. One can see that $\mathrm{U}\mathbb{H}^{n}$ is
 a normal semi-Riemannian homogeneous space (hence it is naturally
 reductive) but not a semi-Riemannian symmetric space.

\subsection{}
 Let us denote by $\mathrm{Gr}_{1,1}(\mathbb{E}^{1,n})$ the
 Grassmann manifold
 of all oriented \textit{timelike} planes in
 $\mathbb{E}^{1,n}$. There exists a natural projection
 $\pi_1:\mathrm{U}\mathbb{H}^{n}\to
 \mathrm{Gr}_{1,1}(\mathbb{E}^{1,n})$:
$$
\pi_1 (\Vec{x},\Vec{v})= \Vec{x} \wedge \Vec{v}. 
$$
 The Grassmann manifold
 is a homogeneous space of
 $\mathrm{SO}^{+}_{1,n}$. In fact, $\mathrm{SO}^+_{1,n}$ acts
 isometrically and transitively on $\mathrm{Gr}_{1, 1} (\mathbb
 E^{1,n})$ by 
 $$
 A \cdot (\Vec{x} \wedge \Vec{v}) = (A \Vec{x}) \wedge ( A \Vec{v}).
 $$
 The isotropy subgroup at $\Vec{e}_{0}\wedge
 \Vec{e}_{1}$ is $\mathrm{SO}_{1,1}\times \mathrm{SO}_{n-1}$.  The 
 semi-Riemannian metric induced from the Killing form of $\mathrm{SO}^+_{1,n}$ 
 has signature $(n-1,n-1)$.

 The contact form $\omega$ on $\mathrm{U}\mathbb{H}^n$ induces 
 a symplectic form $\Omega$ on the Grassmann manifold 
 $\mathrm{Gr}_{1,1}(\mathbb{E}^{1,n})$ so that 
 $\pi_{1}^{*}\Omega=\mathrm{d}\omega$.

 The Grassmann manifold
 $\mathrm{Gr}_{1,1}(\mathbb{E}^{1,n})$ admits an invariant product
 structure $\mathrm{P}$, \textit{i.e.}, an endomorphism field 
 $\mathrm{P}$ satisfying $\mathrm{P}^2=\mathrm{Id}$, 
 compatible with the metric. Moreover, $\mathrm{P}$ is parallel
 with respect to the Levi-Civita connection. The resulting homogeneous
 space $\mathrm{SO}^{+}_{1,n}/\mathrm{SO}_{1,1}\times \mathrm{SO}_{n-1}$ is
 an indefinite para-K\"ahler symmetric space \cite{KK}. 
 The symplectic form $\Omega$ is related to the para-K\"ahler 
 structure by 
 $$
 \Omega(X, Y) = 2 \langle X, \mathrm{P}Y \rangle
 $$
 for all vector fields $X$ and $Y$, see \cite{Honda}.

\subsection{}
 From now on we will concentrate on the case $n=3$. Then 
$$
\mathfrak{so}_{1,3}=
\left\{X=
\left(
\begin{array}{cccc}
0 & x_1 & x_2 & x_3\\
x_1 & 0 & x_{12} & x_{13}\\
x_1 & -x_{12} & 0 & x_{23}\\
x_{3} & -x_{13} & -x_{23} & 0
\end{array}
\right)
=
\left(
\begin{array}{cc}
0 & \Vec{x}^{t}\\
\Vec{x} & b^{x}
\end{array}
\right)
\right\}.
$$
 Here we put 
$$
\Vec{x}=(x_1,x_2,x_3)^{t} \in \mathbb R^3 \subset \mathbb
E^{1,3}\ \ \mbox{and} \ \
b^{x}=\left(
\begin{array}{ccc}
0 & x_{12} & x_{13}\\
 -x_{12} & 0 & x_{23}\\
 -x_{13} & -x_{23} & 0
\end{array}
\right).
$$
 Then one can check that
$$
XY=\left(
\begin{array}{cc}
\Vec{x}^{t}\Vec{y} & \Vec{x}^{t}b^{y}
\\
b^{x}\Vec{y} & b^{x}b^{y}+\Vec{x}\Vec{y}^{t}
\end{array}
\right)
$$
 and
\begin{equation}\label{eq:innerproduct}
\frac{1}{2} \mathrm{tr}\>(XY)=\frac{1}{2}\mathrm{tr}\>(b^xb^y) +\Vec{x}^{t}\Vec{y}
= - (x_{12}y_{12}+x_{13}y_{13}+x_{23}y_{23}) + \langle \Vec{x},\Vec{y}\rangle
\end{equation}
 for 
$$
X=\left(
\begin{array}{cc}
0 & \Vec{x}^{t}\\
\Vec{x} & b^{x}
\end{array}
\right),
\ \
Y=
\left(
\begin{array}{cc}
0 & \Vec{y}^{t}\\
\Vec{y} & b^{y}
\end{array}
\right).
$$
 Now we identify $\mathfrak{so}_{1,3}$ with $\mathbb{R}^{6}$ by
\begin{equation}\label{eq:identifyso}
\left(
\begin{array}{cccc}
0 & x_1 & x_2 & x_3\\
x_1 & 0 & x_{12} & x_{13}\\
x_{2} & -x_{12} & 0 & x_{23}\\
x_{3} & -x_{13} & -x_{23} & 0
\end{array}
\right)
\longleftrightarrow 
(x_{12},x_{13},x_{23}, x_1,x_2,x_3).
\end{equation}
 Then $\mathfrak{so}_{1,3}$ is identified with the semi-Euclidean
 $6$-space $\mathbb{E}^{3,3}=\left(\mathbb{R}^{6},\langle
 \cdot,\cdot\rangle\right)$ with scalar product 
$$
 \langle \cdot,\cdot\rangle = 
- \mathrm{d}x_{12}^2
- \mathrm{d}x_{13}^2
- \mathrm{d}x_{23}^2
+ \mathrm{d}x_{1}^2
+ \mathrm{d}x_{2}^2
+ \mathrm{d}x_{3}^2
\;.
$$ 
 The isotropy subgroup of $\mathrm{SO}^{+}_{1,3}$ acting on
 $\mathrm{U}\mathbb{H}^{3}$ at $(\Vec{e}_0,\Vec{e}_1)$ is 
$$
\left\{
\left(
\begin{array}{cccc}
1 & 0 & 0 & 0 \\
0 & 1 & 0 & 0\\
0 & 0 & \cos \theta  & -\sin \theta \\
0 & 0 & \sin \theta & \cos \theta
\end{array}
\right)
\right\}
\cong \mathrm{SO}_2
$$
with Lie algebra
$$
\left\{
\left(
\begin{array}{cccc}
0 & 0 & 0 & 0 \\
0 & 0 & 0 & 0\\
0 & 0 & 0 & -t \\
0 & 0 & t & 0
\end{array}
\right)
\
\biggr 
\vert
\
t \in \mathbb{R}
\right\}\cong \mathfrak{so}_2.
$$
 The tangent space $T_{(\Vec{e}_0,\Vec{e}_1)} \mathrm{U}\mathbb{H}^3$ is
 naturally  identified with the complement
$$
\left\{
\left(
\begin{array}{cccc}
0 & x_1 & x_2 & x_3\\
x_1& 0 & x_{12} & x_{13}\\
x_2 & -x_{12} & 0 & 0 \\
x_{3} & -x_{13} & 0 & 0
\end{array}
\right)
\right\}
\cong \mathbb{R}^{5}
$$
 of the isotropy algebra in $\mathfrak{so}_{1,3}$. The scalar product 
 $\langle X,Y\rangle$ of
 $X$, $Y\in T_{(\Vec{e}_0,\Vec{e}_1)}\mathrm{U}\mathbb{H}^3$ is computed
 as
\begin{equation}\label{eq:scalarUH3}
\langle X,Y\rangle=-(x_{12}y_{12}+x_{13}y_{13})+x_{1}y_{1}
+x_{2}y_{2}+x_{3}y_{3}.
\end{equation}
 Let $\pi_{2}:\mathrm{U}\mathbb{H}^{3}\to \mathbb{H}^3$ denote the
 natural projection. Then the vertical subspace $\mathcal{V}_2$  of
 $T_{(\Vec{e}_0,\Vec{e}_1)}\mathrm{U}\mathbb{H}^3$ with
 respect to $\pi_2$  at $(\Vec{e}_0,\Vec{e}_1)$ is
$$
\mathcal{V}_2 (\Vec{e}_0,\Vec{e}_1) =
\left\{
\left(
\begin{array}{cccc}
0 & 0 & 0 & 0\\
0 & 0 & x_{12} & x_{13}\\
0 & -x_{12} & 0 & 0\\
0 & -x_{13} & 0 & 0
\end{array}
\right)
\right\}
\cong \mathbb{E}^{2,0}.
$$
 Here $\mathbb E^{2,0}$ is $\mathbb R^2$ with scalar
 product $-\mathrm{d}x_{12}^2-\mathrm{d}x_{13}^2$.
 Since the restriction of the scalar product \eqref{eq:scalarUH3} to
 $\mathcal{V}_2$ is non-degenerate, the horizontal
 subspace $\mathcal{H}_2$ can be defined by
 $\mathcal{H}_2 =\mathcal{V}_2^{\perp}$. 
 The submersion $\pi_2$ satisfies
 $$
 \langle \mathrm{d}\pi_2(X), \mathrm{d}\pi_2(Y)\rangle = 4 \langle X, Y\rangle
 $$
 for any vector fields $X$ and $Y$ on $\mathrm{U}\mathbb H^3$.

 Next we consider the Grassmann manifold $\mathrm{Gr}_{1,1} (\mathbb
 E^{1,3})$.
 The isotropy subgroup of $\mathrm{SO}^{+}_{1,3}$ acting on
 $\mathrm{Gr}_{1,1}(\mathbb{E}^{1,3})$ at $\Vec{e}_{0}\wedge
 \Vec{e}_1$ is $\mathrm{SO}_{1,1}\times \mathrm{SO}_2$ with Lie algebra
$$
\left\{
\left(
\begin{array}{cccc}
0 & s & 0 & 0\\
s & 0 & 0 & 0\\
0 & 0 & 0 & t\\
0 & 0 & -t & 0
\end{array}
\right)
\right\}.
$$
 Therefore, the tangent space $T_{\Vec{e}_{0}\wedge
 \Vec{e}_1}\mathrm{Gr}_{1,1}(\mathbb{E}^{1,3})$
  can be identified with a subspace of
 $\mathfrak{so}_{1,3}$ 
$$
T_{\Vec{e}_{0}\wedge
 \Vec{e}_1}\mathrm{Gr}_{1,1}(\mathbb{E}^{1,3})
\cong
\left\{
\left(
\begin{array}{cccc}
0 & 0 & x_2 & x_3\\
0 & 0 & x_{12} & x_{13}\\
x_{2} & -x_{12} & 0 & 0\\
x_{3} & -x_{13} & 0 & 0
\end{array}
\right)\right\}
\cong
\mathbb{E}^{2,2},
$$
 which is complementary to the isotropy algebra.
 Here $\mathbb E^{2,2}$ is a semi-Euclidean $4$-space with 
 scalar product $-\mathrm{d}x_{12}^2-\mathrm{d}x_{13}^2
 +\mathrm{d}x_{2}^2 +\mathrm{d}x_{3}^2$.
 Let $\pi_{1}:\mathrm{U}\mathbb{H}^3 \to
 \mathrm{Gr}_{1,1}(\mathbb{E}^{1,3})$ denote the natural
 projection. Then the vertical subspace $\mathcal V_1$ 
 of $T_{(\Vec{e}_0,\Vec{e}_1)} \mathrm{U}\mathbb{H}^3$ with respect to
 $\pi_1$ at $(\Vec{e}_0,\Vec{e}_1)$ is
$$
\mathcal{V}_1 (\Vec{e}_0,\Vec{e}_1) =
\left\{
\left(
\begin{array}{cccc}
0 & x_1 & 0 & 0\\
x_1 & 0 & 0 & 0\\
0 & 0 & 0 & 0\\
0 & 0 & 0 & 0
\end{array}
\right)
\right\}\cong \mathbb{E}^{1}.
$$
 Since the restriction of the scalar product \eqref{eq:scalarUH3} to
 $\mathcal{V}_1$ is non-degenerate, the horizontal
 subspace of $\mathcal{H}_1$ can be defined by
 $\mathcal{H}_1 =\mathcal{V}_1^{\perp}$. 
 Moreover, it is easy to see 
 that $\mathrm{d} \pi_1$ restricted to $\mathcal {H}_1$ preserves the scalar
 product. Thus $\pi_1$ is a semi-Riemannian submersion, see appendix
 \ref{sc:submersion}.
\subsection{}\label{sc:Geo}
 Next we consider $\mathrm{Geo}(\mathbb{H}^3)$ the space of all oriented
 geodesics in $\mathbb{H}^3$. Take a point $\gamma\in
 \mathrm{Geo}(\mathbb{H}^3)$, then  $\gamma$ is given by the
 intersection of $\mathbb{H}^3$ with a timelike plane $W\in
 \mathrm{Gr}_{1,1}(\mathbb{E}^{1,3})$. By identifying $\gamma$ with $W$,
 the space $\mathrm{Geo}(\mathbb{H}^3)$ is identified with the 
 Grassmann manifold
 $\mathrm{Gr}_{1,1}(\mathbb{E}^{1,3})$.
 In contrast to the general case of $\mathrm{Gr}_{1,1}(\mathbb{E}^{1,n})$,
 in the case $n=3$, the Grassmann
 $\mathrm{Gr}_{1,1}(\mathbb{E}^{1,n})$ admits an invariant 
 complex structure $J$ compatible with the metric such that 
 the resulting homogeneous
 space $\mathrm{SO}^{+}_{1,3}/\mathrm{SO}_{1,1}\times \mathrm{SO}_{2}$ is
 an indefinite K\"ahler symmetric space.  
 The K\"ahler structure is related to the symplectic form $\Omega$ by 
 $$
 \Omega(X, Y) = -2 \langle X, JY\rangle
 $$
 for any vector fields $X$ and $Y$, see \cite{Honda}. 
%
\subsection{}\label{sc:fibering}
 In addition to the fibrations of $\mathrm{U}\mathbb H^3$ mentioned
 above there also exists a fibration onto the de Sitter $3$-space
 $\mathbb{S}^{1, 2}$: $\pi_3: \mathrm{U}\mathbb{H}^3
 \to \mathbb{S}^{1, 2}, \pi_3 (\Vec{x}, \Vec{v}) = \Vec{v}\in \mathbb
 {S}^{1, 2} \subset
 \mathbb{E}^{1,3}$, where we consider $\mathrm{U}\mathbb H^3$ again as 
 a subspace of $\mathbb E^{1,3} \times \mathbb E^{1,3}$ (see also appendix
 \ref{App:Obata}). Altogether the unit tangent sphere bundle
 $\mathrm{U}\mathbb{H}^3$ has the fibrations:
\begin{align*}
&\pi_{1}:\mathrm{U}\mathbb{H}^3\to \mathrm{Gr}_{1,1}(\mathbb{E}^{1,3});
\ \ \pi_{1}(\Vec{x},\Vec{v})=\Vec{x} \wedge \Vec{v}, \\
&\pi_{2}:\mathrm{U}\mathbb{H}^3\to \mathbb{H}^{3};
\ \
\pi_{2}(\Vec{x},\Vec{v})=\Vec{x}\in \mathbb{H}^3\subset
\mathbb{E}^{1,3}, \\
&\pi_{3}:\mathrm{U}\mathbb{H}^3\to \mathbb{S}^{1, 2};
\ \
\pi_{3}(\Vec{x},\Vec{v})=\Vec{v}\in \mathbb{S}^{1, 2}\subset
\mathbb{E}^{1,3}.
\end{align*}
These fibrations are realized as homogeneous projections:
\begin{align*}
&\pi_{1}:\mathrm{SO}^{+}_{1,3}/\mathrm{SO}_2\to
 \mathrm{SO}^{+}_{1,3}/\mathrm{SO}_{1,1}\times \mathrm{SO}_2, \\
&\pi_{2}:\mathrm{SO}^{+}_{1,3}/\mathrm{SO}_2\to \mathrm{SO}^{+}_{1,3}/
\mathrm{SO}_3,\\
&\pi_{3}:\mathrm{SO}^{+}_{1,3}/\mathrm{SO}_2\to \mathrm{SO}^{+}_{1,3}/
\mathrm{SO}_{1,2}.
\end{align*}

\section{Surfaces in $\mathbb{H}^3$}\label{sc:SurfacesH3}

\subsection{}
 Let $f:M\to \mathbb{H}^3\subset \mathbb{E}^{1,3}$ be a conformal
 immersion of a Riemann surface with unit normal vector field $n$.
  Clearly, by replacing, if necessary, the unit normal $n$ of an
 immersion by $-n$, we can assume that the mean curvature satisfies
 $H\geq 0$.
 Let $\mathbb D$ denote the universal cover of $M$. Since $\mathbb D
 \cong S^2$ can only occur for totally umbilic CMC immersions with $H>1$
 and cannot occur for CMC immersions with $0\leq H \leq 1$, \cite{Spivak}, 
 we can assume that $\mathbb D \subset \C$ being an open. 
 (Usually we will assume $\mathbb D = \C$ or $\mathbb D=$ open unit disk.) 
 Hence, without loss of generality,
 $\mathbb D$ is open in $\C$, and the first
 fundamental form $\mathrm{I}$ is written as
\begin{equation}
\mathrm{I}=e^{u}\mathrm{d}z\mathrm{d}\bar{z}.
\end{equation}
 The \textit{Hopf differential} of $(M,f)$ is a quadratic differential on 
 $M$ defined by 
$$
Q\>\mathrm{d}z^{2},
\ \
Q=\langle f_{zz},n\rangle.
$$
 The \textit{Gauss-Codazzi equations} of $(M,f)$ are given by
\begin{equation}
u_{z\bar{z}}+\frac{1}{2}(H^2-1)e^{u}-2|Q|^{2}e^{-u}=0, \;\;
Q_{\bar z}=\frac{1}{2}H_{z}e^{u}.
\end{equation}
 Here $H$ denotes the \textit{mean curvature} of $(M,f)$, 
 which is explicitly given by 
 $H = 2 e^{-u} \langle f_{z \bar z}, n\rangle$.
 The constancy of the mean curvature is characterized as follows.

\begin{proposition}
 Let $M$ be a Riemann surface. 
 A conformal immersion $f:M\to \mathbb{H}^3$ is of constant mean curvature
 if and only if its Hopf differential is holomorphic.
\end{proposition}
 Assume that $H$ is constant, then the Gauss-Codazzi equations
 are invariant under the deformation
$$
Q\longmapsto \lambda^{-1}Q,\ \ 
\lambda\in S^1.
$$ 
 Hence, on the region $\mathbb{D}$, there exists a $1$-parameter 
 deformation family of conformal constant mean curvature immersions 
 $\{f_{\lambda}\}$ through $f_1=f$. All these immersions have the 
 same induced metric and mean curvature. The family $\{f_\lambda\}$ 
 is referred to as the 
 \textit{associated family} of the original immersion $f$. 

\subsection{}\label{sc:tangentGauss}
 Let again $f:M\to \mathbb{H}^3$ be a conformal immersion with unit
 normal $n$. For the purposes of this paper it will be important to
 consider the map
\begin{equation}\label{eq:mapsinUH3}
F:=(f,n):M\to \mathrm{U}\mathbb{H}^3.
\end{equation}
 In our paper, the smooth map $F$ will be  called the \textit{Gauss
 map} of $f$. The Gauss map satisfies  $\langle
 \mathrm{d}f,n\rangle=0$. In section \ref{sc:unit tangent sphere bundle}
 we have introduced the canonical contact form
 $\omega$ of $\mathrm{U}\mathbb{H}^3$. By the definition of $\omega$, we
 have
$$
F^{*}\omega=\langle \mathrm{d}f,n\rangle.
$$
 Hence the Gauss map $F$ satisfies the \textit{Legendre condition}:
$$
 F^{*}\omega=0.
$$
 Note that the Gauss map $F$ is also called the \textit{Legendre
 lift} of $f$ as in \cite{Arnold}.
 The Legendre property will be discussed 
 in section \ref{sc:contactmfd}. The
 following result which is very important for this paper is due to
 T.~Ishihara \cite{Ishihara}. 
\begin{proposition}{\rm (T.~Ishihara)}\label{prop:Ishihara}
 Let $M$ be a Riemann surface. 
 A conformal immersion $f:M \to \mathbb{H}^3$ has constant mean curvature  
 if and only if its Gauss map is harmonic with respect to the metric 
 induced from \eqref{scalarproduct}. 
\end{proposition} 
 In section \ref{sc:Harmonicity}, we will give a proof of this result in
 terms of frames and the Sym formula.

\section{The $2 \times 2$-matrix model for immersions into $\mathbb H^3$}\label{sc:2by2}
\subsection{}\label{sc:Minkowski}
 The Minkowski 4-space $\mathbb{E}^{1,3}$ is identified with the space
 $\mathrm{Her}_2 \C$ of all complex Hermitian $2\times
 2$-matrices: 
\begin{equation}\label{eq:basis}
\mathbb E^{1,3} \cong \mathrm{Her}_2 \C=\left \{
\xi=\left (
\begin{array}{cc}
\xi_0+\xi_1 &
\xi_3-i\xi_2 \\
\xi_3+i\xi_2 &
\xi_0-\xi_1 
\end{array}
\right)\ 
\Biggr \vert
\
\xi_0,
\xi_1,\xi_2,\xi_3
\in \mathbb{R}\
\right \}.
\end{equation}
 The space $\mathrm{Her}_2 \C$ is spanned by the orthonormal basis
$$
\Vec{e}_{0}=\left(
\begin{array}{cc}
1 & 0\\
0 & 1
\end{array}
\right),
\ \
\Vec{e}_{1}=\left(
\begin{array}{cc}
1 & 0\\
0 & -1
\end{array}
\right),
\
\Vec{e}_{2}=\left(
\begin{array}{cc}
0 & -i\\
i & 0
\end{array}
\right),
\
\Vec{e}_{3}=\left(
\begin{array}{cc}
0 & 1\\
1 & 0
\end{array}
\right).
$$
 It is easy to see that $-\det \xi=-\xi_0^2+\xi_1^2+\xi_2^2+ \xi_3^2$ 
 for $\xi \in \mathrm{Her}_2 \C$.  Thus the Lorentzian metric
 of $\mathrm{Her}_2 \C$ is described as 
$$
\langle \xi, \eta\rangle=-\frac{1}{2} {\mathrm{tr}}( \xi \Vec{e}_2 \eta^t \Vec{e}_2), \ \ \xi, \eta \in \mathrm{Her}_2 \C.
$$
 In particular, we have
$$
 \langle \xi, \xi\rangle = -\det \xi, \ \ \xi \in \mathrm{Her}_2 \C.
$$
 Thus  we have the identification: 
$$
\mathbb{H}^3=\{ \xi \in \mathrm{Her}_2 \C \ \vert \ \det \xi=1,\
 \mathrm{tr}\: \xi>0\}.
$$
 The special linear group $G=\mathrm{SL}_{2}\C$ acts isometrically and
 transitively on the hyperbolic $3$-space via the action:
\begin{equation}\label{eq:actionA}
\mathrm{A}:\mathrm{SL}_{2}\C\times \mathbb{H}^3\to 
\mathbb{H}^3,\
(g,\xi) \mapsto g\> \xi \>g^{*},
\end{equation}
 where $g^*$  denotes $\bar g^t$.
 The isotropy subgroup of this action at $\Vec{e}_{0}$ 
 is the special unitary group $\mathrm{SU}_2$. Hence
 $\mathbb{H}^3$ is represented by
 $\mathbb{H}^3=G/K=\mathrm{SL}_{2}\C/\mathrm{SU}_2$ as a 
 \textit{Riemannian symmetric space}. The natural projection $\pi:G\to
 \mathbb{H}^3$ is given explicitly by $\pi(g)=gg^{*},\ g\in G$. In other
 words, $\mathbb{H}^3$ is represented as
$$
\mathbb{H}^{3}=\{gg^{*}\
\vert
\
g \in G\}.
$$
\begin{remark}
{\rm 
 It is important to note that in this context the simple Lie group
 $\mathrm{SL}_{2}\C$ is regarded as a simple real Lie group 
 and as a double covering of the special Lorentz group
 $\mathrm{SO}^{+}_{1,3}$. The real Lie algebra $\mathfrak{sl}_{2}\C$
 is spanned by the basis 
\begin{equation}\label{orthonormalbasis}
\{i\Vec{e}_{1},i\Vec{e}_{2},i\Vec{e}_{3}, \Vec{e}_{1},\Vec{e}_{2},
 \Vec{e}_{3}\}.
\end{equation}
 The bi-invariant semi-Riemannian metric $\langle\cdot,\cdot \rangle$ on
 $G$ corresponding to the scalar product
 (\ref{eq:innerproduct}) via the 
 isomorphism \eqref{eq:basis} has the signature $(-,-,-,+,+,+)$. The
 tangent space $\mathfrak{m}=T_{\Vec{e}_0}\mathbb{H}^3$ is given by
$$
\mathfrak{m}=\mathfrak{sl}_{2}\C
\cap \mathrm{Her}_2 \C=
\left\{
\left(
\begin{array}{cc}
a & \bar{b}
\\
b & -a
\end{array}
\right)
\
\biggr
\vert
\ a \in \mathbb{R},\
b\in \C
\right\}=
\mathbb{R}\Vec{e}_{1}\oplus \mathbb{R}\Vec{e}_{2}
\oplus \mathbb{R}\Vec{e}_3.
$$
}
\end{remark}

\subsection{}\label{sc:tangentspherebundle}
 The unit tangent sphere bundle $\mathrm{U}\mathbb{H}^3$ is represented
 as
$$
\mathrm{U}\mathbb{H}^{3}=
\{ (\Vec{x},\Vec{v})\in \mathrm{Her}_2 \C\times 
\mathrm{Her}_2 \C
\
\vert
\
\det \Vec{x}=1, \mathrm{tr}\>\Vec{x}>0,\>
\det \Vec{v}=-1, \> \langle \Vec{x},\Vec{v}\rangle =0\}. 
$$
 The special liner group $G=\mathrm{SL}_{2}\C$ acts
 isometrically and transitively on $\mathrm{U}\mathbb{H}^{3}$ via the
 action:
$$
g\cdot (\Vec{x},\Vec{v})=(g \Vec{x} g^{*}, g \Vec{v} g^{*}).
$$
 The isotropy subgroup of $G$ at $(\Vec{e}_{0},\Vec{e}_{1})$ is 
$$
H=\left\{
\left(
\begin{array}{cc}
e^{i\theta} & 0\\
0 & e^{-i\theta}
\end{array}
\right)
\right\}=\mathrm{U}_1.
$$
 Thus the unit tangent sphere bundle $\mathrm{U}\mathbb H^3$ is as 
 homogeneous space $G/H = \mathrm{SL}_2 \C/\mathrm{U}_1$.
 The Lie algebra $\mathfrak{h}$ of $H$ is
\begin{equation}\label{eq:h}
\mathfrak{h}=\left\{
\left(
\begin{array}{cc}
ia_{2} & 0\\
0 & -ia_{2}
\end{array}
\right)
\
\biggr
\vert
\
a_{2}\in \mathbb{R}\>
\right\}=\mathfrak{u}_1=\mathbb{R}(i\Vec{e}_{1}).
\end{equation}
 The tangent space $\mathfrak{p}:=T_{(\Vec{e}_0,\Vec{e}_1)}
 \mathrm{U}\mathbb{H}^3$ is given by
\begin{equation}\label{eq:p}
\mathfrak{p}=
\left\{
\left(
\begin{array}{cc}
a_{1} & b_{1}+ib_{2}\\
c_{1}+ic_{2} & -a_{1}
\end{array}
\right)
\
\biggr
\vert
\
a_{1},b_{1},b_{2},
c_{1},c_{2}
\in \mathbb{R}
\>
\right\}.
\end{equation}
 Note that
\begin{eqnarray*}
\mathfrak{p}&=&
\{\mathbb{R}\Vec{e}_1\oplus
\mathbb{R}\Vec{e}_2
\oplus
\mathbb{R}\Vec{e}_3\}
\oplus
\{
\mathbb{R}(i\Vec{e}_2)\oplus
\mathbb{R}(i\Vec{e}_3)
\}\\
&=& \mathfrak{m}\oplus 
\{
\mathbb{R}(i\Vec{e}_2)\oplus
\mathbb{R}(i\Vec{e}_3)
\}.
\end{eqnarray*}
 This equation shows that the horizontal distribution $\mathcal{H}_{2}$
 and the vertical distribution $\mathcal{V}_{2}$ with respect to the fibering
 $\pi_2 : \mathrm{U}\mathbb{H}^{3}\to \mathbb{H}^3$ are generated via
 the identifications \eqref{eq:identifyso} and \eqref{eq:basis} by
$$
\mathcal{H}_{2}=\mathfrak{m},\ \
\mathcal{V}_{2}=\mathbb{R}(i\Vec{e}_2)\oplus
\mathbb{R}(i\Vec{e}_3).
$$

\subsection{}
 Next we consider the fibering $\pi_{1}:\mathrm{U}
 \mathbb{H}^3\to \mathrm{Gr}_{1,1}(\mathbb{E}^{1,3})$.
 From section \ref{sc:fibering}, the Grassmann manifold
 $\mathrm{Gr}_{1,1}(\mathbb{E}^{1,3})$ is a homogeneous space of $G = {\rm
 SL}_2 \C$. The isotropy subgroup of $G$ at $\Vec{e}_{0}\wedge \Vec{e}_{1}$
 is 
$$
D=\mathrm{GL}_{1}\mathbb{C}=\left\{\left(
\begin{array}{cc}
w &0\\
0 & 1/w
\end{array}
\right)
\
\biggr 
\vert
\
w\in \C^{\times}
\right\}
\cong \C^{\times}.
$$ 
 The Lie algebra $\mathfrak{d}$ of $D$ is
\begin{equation}\label{eq:d}
\mathfrak{d}=
\left\{\left(
\begin{array}{cc}
v &0\\
0 & -v
\end{array}
\right)
\
\biggr 
\vert
\
v\in \C
\right\}=
\mathbb{R}\Vec{e}_{1}\oplus \mathbb{R}(i\Vec{e}_1).
\end{equation}
 The tangent space $\mathfrak{q}:=T_{\Vec{e}_{0}\wedge
 \Vec{e}_1}\mathrm{Gr}_{1,1} (\mathbb{E}^{1,3})$ of $\mathrm{Gr}_{1,1}
 (\mathbb{E}^{1,3})$ at $\Vec{e}_{0}\wedge \Vec{e}_1$ is
\begin{equation}\label{eq:q}
 \mathfrak{q}=\mathbb{R}\Vec{e}_{2}\oplus \mathbb{R}(i\Vec{e}_{2})
 \oplus \mathbb{R}\Vec{e}_{3}\oplus \mathbb{R}(i\Vec{e}_{3}).
\end{equation}
 The horizontal and vertical distributions with respect to $\pi_1$ are
 generated by
$$
\mathcal{H}_{1}=\mathfrak{q},
\ \
\mathcal{V}_{1}=\mathbb{R}\Vec{e}_{1}.
$$

\subsection{}
 Next we consider the fibering $\pi_{3}:\mathrm{U}
 \mathbb{H}^3\to \mathbb S^{1,2}$.
 From section \ref{sc:fibering}, the de Sitter $3$ space $\mathbb
 S^{1,2}$ is a homogeneous space of $G = {\rm
 SL}_2 \C$. The isotropy subgroup of $G$ at $\Vec{e}_{1}$ is 
 $\mathrm{SU}_{1,1}$. Hence
 $\mathbb{S}^{1, 2}$ is represented by
 $\mathbb{S}^{1, 2}=G/K=\mathrm{SL}_{2}\C/\mathrm{SU}_{1,1}$ as a 
 \textit{Lorentzian symmetric space}. The natural projection 
 $\pi:G\to \mathbb{S}^{1, 2}$ is given explicitly by 
 $\pi(g)=g\Vec{e}_1g^{*},\ g\in G$. In other
 words, $\mathbb{S}^{1, 2}$ is represented as
$$
\mathbb{S}^{1, 2}=\{g\Vec{e}_1g^{*}\
\vert
\
g \in G\}.
$$
 The horizontal and vertical distributions with respect to 
 $\pi_3$ are generated by
$$
\mathcal{H}_{3}=\mathbb{R}\Vec{e}_1 \oplus \mathbb{R}(i \Vec{e}_2)\oplus 
 \mathbb{R}(i\Vec{e}_3),
\ \
\mathcal{V}_{3}=\mathbb R\Vec{e}_2\oplus \mathbb R\Vec{e}_3.
$$
\begin{remark}
{\rm 
 Precisely speaking, to represent the de Sitter 3-space
 $\mathbb{S}^{1,2}$ as a Lorentzian symmetric space, we need to equip the 
 scalar product $-\langle \cdot,\cdot\rangle$ on the linear space
 $\mathbb{R}\Vec{e}_1 \oplus \mathbb{R}(i \Vec{e}_2)\oplus 
 \mathbb{R}(i\Vec{e}_3)$.
}
\end{remark}

\subsection{}
 Now let again  $f:M\to \mathbb{H}^3$ be a CMC surface with unit normal
 $n$. Take a simply connected complex coordinate region 
 $\mathbb D \subset \C$ as before. Denote by
 $(x,y)$ the associated isothermal coordinates, \textit{i.e.},
 $z=x+iy$, where the induced metric is expressed as
 $\mathrm{I}=e^{u} (\mathrm{d}x^2+\mathrm{d}y^2)$. 

 The \textit{coordinate frame} $\Psi$ of $f$ with respect to
 $(x,y)$ is a map from $\mathbb{D}$ into the Lorentz group
 $\mathrm{SO}^{+}_{1,3}$ defined by
$$
\Psi=(f,n,e^{-u/2}f_{y},e^{-u/2}f_{x}).
$$
 As mentioned above, the action \eqref{eq:actionA} induces a double covering
 $\tilde \pi:\mathrm{SL}_{2}\C\to \mathrm{SO}^{+}_{1,3}$. Since
 $\mathbb D$ is simply connected, the lift $\hat{\Phi}$ of the
 coordinate frame $\Psi$ to $\mathrm{SL}_{2}\C$ is
 determined uniquely by
\begin{equation}\label{eq:hatframe}
\Vec{e}_j \longmapsto \hat \Phi \Vec{e}_j \hat \Phi^*
\end{equation}
 up to sign. The lift $\hat{\Phi}$ satisfies the following
 Gauss-Weingarten formulas (see appendix
 \ref{sc:Gauss-Codazzi}): 
\begin{equation}\label{eq:hatMau1}
\hat{\Phi}^{-1}\hat{\Phi}_z=
\left(
\begin{array}{cc}
u_{z}/4 & \frac{1}{2}(H+1)e^{u/2} \\
-Q e^{-u/2} & -u_{z}/4
\end{array}
\right), \;\;
\hat{\Phi}^{-1}\hat{\Phi}_{\bar z}=
\left(
\begin{array}{cc}
-u_{\bar z}/4 &
{\bar Q}e^{-u/2} \\
-\frac{1}{2}(H-1)e^{u/2} &
u_{\bar z}/4
\end{array}
\right).
\end{equation}
 Let us consider the associated family $\{f_\lambda\}_{\lambda\in S^1}$
 of $f=f_1$. The immersion $f_{\lambda}$ has the Hopf differential
 $\lambda^{-1}Q\>\mathrm{d}z^2$. The corresponding
 $\mathrm{SL}_{2}\C$-valued frame is denoted by
 $\hat{\Phi}_{\lambda}$. For our purposes, it will be useful to perform
 the following gauge transformation
\begin{equation}\label{eq:tildephi}
\tilde{\Phi}:=\hat{\Phi}_{\lambda^2}
\left(
\begin{array}{cc}
\sqrt{\lambda} & 0\\
0 & 1/\sqrt{\lambda}
\end{array}
\right).
\end{equation}
 Then we obtain $\tilde U = \tilde{\Phi}^{-1}\tilde{\Phi}_{z}$
 and $\tilde V = \tilde{\Phi}^{-1}\tilde{\Phi}_{\bar z}$ with
\begin{equation}\label{Lax}
{\tilde U}=
\left(
\begin{array}{cc}
u_{z}/4 & \frac{\lambda^{-1}}
{2}(H+1)e^{u/2} \\
-\lambda^{-1}Q e^{-u/2} & -u_{z}/4
\end{array}
\right), \;\;
{\tilde V}=
\left(
\begin{array}{cc}
-u_{\bar z}/4 &
\lambda{\bar Q}e^{-u/2} \\
-\frac{\lambda}{2}(H-1)e^{u/2} &
u_{\bar z}/4
\end{array}
\right).
\end{equation}
%
%
 From \eqref{Lax}, $\tilde U$ and $\tilde V$
 are elements  of the \textit{loop algebra} of $\mathfrak{sl}_2 \C$:
$$
 \Lambda \mathfrak{sl}_2 \C_{\sigma} = \left\{
 g:S^1\to \mathfrak{sl}_{2}\C\ \vert \
 g(-\lambda)=\sigma g(\lambda)\>\right\},
$$
 where 
\begin{equation}\label{eq:sigma}
 \sigma = \mathrm{Ad}(\Vec{e}_1).
\end{equation}
 It will turn out to be useful to consider the analytic loops in
 $\Lambda \mathfrak{sl}_2 \C_{\sigma}$. Those loops will
 be denoted by $\tilde \Lambda \mathfrak{sl}_2 \C_{\sigma}$. 
 We will use a similar notation for loop groups.
 It is easy to see that 
 $\tilde U, \tilde V \in \tilde \Lambda \mathfrak{sl}_2 \C_{\sigma}$ 
 and $\tilde \Phi \in \tilde \Lambda \mathrm{SL}_2 \C_{\sigma}$.
 Then we obtain the {\it Sym formula} for CMC surfaces in $\mathbb H^3$
 which is easily seen to be equivalent to the corresponding formula 
 used in \cite{Bob:Russ},  \cite{BaBo}. 
\begin{proposition}\label{prop:Sym}
 Let $\tilde{\Phi}$ be a solution to \eqref{Lax}. Then
$$
{\tilde f}_{\lambda}={\tilde \Phi}{\tilde \Phi}^{*}
$$
 is a loop of immersions of constant mean curvature $H$ with unit normal
 vector field 
$$
 {\tilde n}_{\lambda}={\tilde \Phi} \Vec{e}_{1}\>
{\tilde \Phi}^{*}.
$$
 For $\lambda =1$ we obtain $\tilde f_{\lambda =1} = f$.
\end{proposition}
{\bf Proof}.
 A direct computation shows
\begin{eqnarray*}
\tilde{f}_{z}=
\frac{\lambda^{-1}}{2}e^{u/2}{\tilde\Phi}
(\Vec{e}_{3}+i \Vec{e}_{2})\tilde \Phi^*,\;\;
\tilde{f}_{\bar z}=
\frac{\lambda}{2}e^{u/2}{\tilde \Phi}
(\Vec{e}_{3}-i \Vec{e}_{2}){\tilde \Phi}^*.
\end{eqnarray*}
 From these equations we obtain
$$
\mathrm{I}_{\lambda}=e^{u}\mathrm{d}z\mathrm{d}{\bar z}.
$$ 
 Hence $\tilde{f}_{\lambda}$ is an immersion for all $\lambda \in
 S^1$. The unit normal for $\tilde f_{\lambda}$ is given by ${\tilde
 n}_{\lambda}=\tilde{\Phi}\Vec{e}_1\tilde \Phi^*$. Thus 
$$
{\tilde n}_{z}=
{\tilde \Phi}(U\Vec{e}_{1}+\Vec{e}_{1}V^{*}){\tilde \Phi}^*=
-\frac{\lambda^{-1}}{2}\left\{ He^{u/2}\tilde{\Phi}(\Vec{e}_{3}+\Vec{e}_{2})
\tilde{\Phi}^*+2 Q e^{-u/2}\tilde{\Phi}
(\Vec{e}_{3}-\Vec{e}_{2})\tilde{\Phi}^*\right\}.
$$
 From these equations one can check that each $\tilde{f}_{\lambda}$ has
 constant mean curvature $H$ and Hopf differential $\lambda^{-2}Qdz^2$.
 $\Box$
\begin{remark}
{\rm
 The construction above shows that all one needs for the construction
 of constant mean curvature surfaces are a real number $H$, a
 holomorphic function $Q$ and a real valued function $u$ such that 
 with $\tilde U$ and $\tilde V$ as in \eqref{Lax} the one 
 form $\tilde \alpha = \tilde U dz + \tilde V d\bar z$ is integrable.
}
\end{remark}
\subsection{}\label{sc:Harmonicity}
 In this section, we discuss the harmonicity of the Gauss map 
 associated with $\tilde f_{\lambda}$. Let $f:M\to \mathbb{H}^3$ be a
 conformal immersion as above.
 Let $\hat{\Phi}$ denote the $\mathrm{SL}_{2}\C$-valued frame
 which is a lift of the coordinate frame $\Psi$
 and put $\hat{\alpha}:=\hat{\Phi}^{-1}\mathrm{d}\hat{\Phi}$.
 Then we decompose $\hat{\alpha}$ as
$$
\hat{\alpha}=\hat{\alpha}_{\mathfrak{h}}+
\hat{\alpha}_{\mathfrak{p}}^{\prime}+
\hat{\alpha}_{\mathfrak{p}}^{\prime\prime}
$$
 according to the Lie algebra decomposition
 $\mathfrak{sl}_2 \C=\mathfrak{h}\oplus \mathfrak{p}$ as in
 \eqref{eq:h} and \eqref{eq:p} respectively. A direct computation shows
$$
[\hat{\alpha}^{\prime}_{\mathfrak{p}}\wedge
\hat{\alpha}^{\prime\prime}_{\mathfrak{p}}]
=-\frac{1}{4}
\left\{(H^2-1)e^{u}-4|Q|^{2}e^{-u}
\right\}\Vec{e}_{1}\>\mathrm{d}z\wedge \mathrm{d}\bar{z}.
$$
 This is contained in $\mathfrak h$. Therefore the $\mathfrak p$-part  of 
 $[\hat \alpha_{\mathfrak p}^{\prime} \wedge 
 \hat \alpha_{\mathfrak p}^{\prime \prime}]$ vanishes
 and $\hat{\alpha}$ satisfies the admissibility condition
 (\ref{admissible}). Moreover, it is easy to check that 
$$
\mathrm{d}(*\hat{\alpha}_{\mathfrak{p}})
+
[\hat{\alpha}^{\prime}\wedge\>*
\hat{\alpha}_{\mathfrak{p}}]
=-ie^{u/2}
\left(
\begin{array}{cc}
0 & H_{\bar z}\\
-H_{z} & 0
\end{array}
\right)\mathrm{d}z\wedge \mathrm{d}\bar{z}
$$
 holds. Since \eqref{harmonicity} describes the harmonicity of the 
 Gauss map, this formula implies Proposition \ref{prop:Ishihara}.

\begin{remark}
{\rm 
 On the unit tangent sphere bundle $\mathrm{U}\mathbb{H}^3$, we can
 define  Riemannian metrics so that the natural projection $\pi_2$ is a
 Riemannian submersion. One of such metrics is the \textit{Sasaki
 lift metric}. It is not difficult to see that the Gauss map $F$ of a
 non-minimal CMC
 surface is \textit{never} harmonic  with respect to the Sasaki lift
 metric. See \cite[p.~271, proof of Corollary]{JR}.

 In 3-dimensional homogeneous Riemannian spaces of
 non-constant curvature, the harmonicity of the Gauss map with respect to the 
 Sasaki lift metric is a very
 strong restriction for CMC surfaces. In fact, the
 only CMC surfaces with harmonic Gauss map in a 3-dimensional
 homogeneous Riemannian space with 4-dimensional isometry group are inverse
 images of geodesics under the Hopf-fibration or totally geodesic
 leaves. The latter case only occurs if the ambient space is a
 direct product space \cite{Sanini}, \cite{Tamura}. 

 CMC surfaces with harmonic Gauss maps in 
 3-dimensional homogeneous Riemannian spaces with 3-dimensional isometry
 group have been classified by J.~Van der Veken and the second named
 author \cite{IV}.
}
\end{remark}
\section{CMC surfaces with $H>1$}\label{sc:CMCH>1}
 As pointed out in the introduction, the case $H =1$ is
 special and has been investigated already intensively.
 Therefore this case will not be considered in this paper.
 In this section, we study CMC surfaces with mean curvature $H$ such
 that $H>1$. In this case, we may write $H=\mathrm{coth}\>q,\; q \in
 \mathbb R_{>0}$.

 We perform a gauge transformation:
\begin{equation}\label{eq:extendedH>1}
\Phi:=\tilde{\Phi}
\left(
\begin{array}{cc}
e^{q/4} & 0\\
0 & e^{-q/4}
\end{array}
\right).
\end{equation}
 We call $\Phi$ the \textit{extended frame} of a CMC immersion $f$ with
 $H>1$.  Moreover we put
$$
\mathscr{H}=e^{-q}(H+1) \in \mathbb R,
\ \
\nu=-e^{-q/2}\lambda.
$$
 Then $H = \coth q$ implies $\mathscr H =e^{-q}(H+1) = e^q (H-1)$. 
 Moreover, the Lax pair
$$
 U=\Phi^{-1}\Phi_{z},\  \
 V=\Phi^{-1}\Phi_{\bar z}
$$
 is given by 
\begin{equation}\label{LaxH>1}
U=\left(
\begin{array}{cc}
u_{z}/4 & -\frac{1}{2}\nu^{-1}\mathscr{H} e^{u/2}
\\
\nu^{-1}Q e^{-u/2} & -u_{z}/4
\end{array}
\right),
\ \
V=
\left(
\begin{array}{cc}
-u_{\bar z}/4 & 
-\nu \bar{Q}e^{-u/2}
\\
\frac{1}{2}\nu \mathscr{H}e^{u/2}
&
u_{\bar z}/4
\end{array}
\right).
\end{equation}
 Clearly, the matrices $U$ and $V$ are holomorphic in the parameter $\nu
 \in \C^{\times}$. 
 In particular, the gauged frame $\Phi$ defined in
 \eqref{eq:extendedH>1} above can be considered to be a holomorphic
 function in $\nu$, where $\nu$ is restricted to the circle of radius $ r
 = e^{-q/2}$.
 Noting that everything is holomorphic in $\nu \in \C^{\times}$, 
 it is straightforward to  check that the $1$-form 
 $\alpha = U \mathrm{d} z + V \mathrm{d} \bar z$ 
 is fixed by the following loop algebra automorphism:
\begin{equation}\label{eq:invc3}
 \tau_3:g(\nu)\longmapsto -g^{*}(1/\bar{\nu}).
\end{equation}
 This automorphism is said to be of type $C_3$ (almost compact
 automorphism of the third kind), \cite{Kobayashi}. The Maurer-Cartan form
 $\alpha=\Phi^{-1}\mathrm{d}\Phi$ has the decomposition
\begin{equation}\label{eq:decalpha}
\alpha=\nu^{-1}\alpha_{-1}+\alpha_{0}+\nu\alpha_{1}.
\end{equation}
 $\tau_3 (\alpha) = \alpha$ translates into 
$$
\overline{\alpha_{0}}=-\alpha_{0},
\ \
\alpha_{-1}= - \alpha_{1}^{*}.
$$
 The mapping $\Phi$ takes values in the twisted loop group
$$
\Lambda \mathrm{SL}_{2}\C_{\sigma,\tau_3}:=
\left\{
g: S^1\to \mathrm{SL}_{2}\C
\
\vert
\
g(-\nu)=\sigma g(\nu),
\ \
\tau_3 (g)(\nu)=g(\nu)\>
\right\},
$$
 where $\sigma$ is defined in \eqref{eq:sigma} and 
 $\tau_3(g)(\nu) = g(1/\bar{\nu})^{*-1}$.

 Moreover, the first formula in Proposition \ref{prop:Sym} can now be
 reinterpreted as
\begin{equation}\label{eq:SymH>1}
f_{\nu}:=\Phi 
\left(
\begin{array}{cc}
e^{-q/2} & 0\\
0 & e^{q/2}
\end{array}
\right)
\Phi^{*}.
\end{equation}
 This CMC immersion into $\mathbb H^3$ has mean curvature $H = \coth q$
 and the unit normal
\begin{equation}\label{eq:normalH>1}
n_{\nu}=
\Phi
\left(
\begin{array}{cc}
e^{-q/2} & 0\\
0 & -e^{q/2}
\end{array}
\right)
\Phi^{*}.
\end{equation}
 Note that the matrix $\diag (e^{-q/2}, e^{q/2})$ corresponds under
 the isomorphism \eqref{eq:basis} to the point
 $(\cosh(q/2),-\sinh(q/2),0,0) \in \mathbb H^3$.
 Conversely, the following result holds:
\begin{proposition}\label{prop:fundamentalH>1}
 Let $\mathscr H$ be a positive real number, $u$ a real valued function and
 $Q$ a holomorphic function on the simply connected domain 
 $\mathbb D \subset \C$. 
 Let $\nu$ be a complex parameter.
 Assume that the differential $1$-form 
 $\alpha = U\mathrm{d}z + V \mathrm{d}\bar z$ is
 integrable and let $\Phi$ denote a solution to $ \Phi^{-1} \mathrm{d} \Phi =
 \alpha$. Then \eqref{eq:SymH>1} defines for $\nu$ of absolute value
 $e^{-q/2}$ a CMC immersion into $\mathbb H^3$ with mean
 curvature $ H = \coth q$ and the unit normal is defined in 
 \eqref{eq:normalH>1}.
\end{proposition}

 The Lawson correspondence (see appendix \ref{sc:Lawsoncorrespond} for
 more details) between CMC surfaces in
 $\mathbb H^3$ with mean curvature $H>1$ and CMC surfaces in $\mathbb E^3$ 
 with mean curvature $\mathscr H$ has now in our setting the
 following simple explanation:
 Consider a CMC surface in $\mathbb H^3$ and let $U$ and $V$ denote the
 associated matrices, but now consider these  matrices as functions of
 $\nu \in S^1$. 
 Then it is straightforward to check that the $1$-form  $\alpha = U\mathrm{d}z +
 V \mathrm{d}\bar z$ has the decomposition as in \eqref{eq:decalpha} and
 satisfies all the conditions for being the Maurer-Cartan form of the
 extended frame of
 some CMC surface in $\mathbb E^3$ with Hopf differential $Q$ and mean
 curvature $\mathscr H$.

 The converse construction, starting from some CMC surface in $\mathbb E^3$ 
 and ending up with some CMC surface in $\mathbb H^3$ works out analogously.

 Thus for the case $H>1$, constructions of CMC surfaces in
 $\mathbb{H}^3$ are reduced to those for CMC surfaces in Euclidean
 3-space \cite{DPW}. In other words, we can construct CMC surfaces in
 $\mathbb{H}^3$ with $H>1$ via the generalized Weierstrass
 type representation (DPW method) for CMC surfaces in Euclidean 3-space.  
 See appendix \ref{sc:Lawsoncorrespond} for an other explanation for this
 fact. 

 \begin{remark}
{\rm 
 One can make this relation into a 1-1 relation by fixing the
 initial conditions of the extended frames at some base point.
}
\end{remark}

Let us  identify 
$\mathfrak{su}_2$ with Euclidean $3$-space via the 
correspondence
$$
x_{1}(i\Vec{e}_1)+x_{2}(i\Vec{e}_2)+
x_{3}(i\Vec{e}_3)
\longleftrightarrow
(x_1,x_2,x_3).
$$
The Euclidean inner product 
$\mathrm{d}x_{1}^{2}+\mathrm{d}x_{2}^{2}+\mathrm{d}x_{3}^{2}$ corresponds to 
the inner product
$$
\langle X,Y \rangle=-\frac{1}{2}\mathrm{tr}\>(XY),
\ \  X, Y \in \mathfrak{su}_2.
$$
Then one can see that 
$$
\varphi_{\nu}=\mathrm{Ad}(\Phi)(i\Vec{e}_{1}):\mathbb{D}\times C_r
\to \mathbb{S}^{2}\subset \mathfrak{su}_2
$$ is a loop of harmonic maps, where $C_r$ is a radius $r$ circle.

\section{CMC surfaces with $0 \leq H<1$}\label{sc:CMCH<1}
\subsection{}
 Now we start our study of CMC surfaces with mean curvature $H$ such
 that $0\leq H <1$. In this case, we may write $H=\mathrm{tanh}\>q, \;
 q \in \mathbb R_{\geq 0}$.

 We perform the gauge transformation:
\begin{equation}\label{eq:extendedH<1}
\Phi:=\tilde{\Phi}
\left(
\begin{array}{cc}
e^{(q+\pi{i})/4} & 0\\
0 & e^{-(q+\pi{i})/4}
\end{array}
\right).
\end{equation}
 We call $\Phi$ the \textit{extended frame} of a CMC immersion $f$ with
 $0 \leq H <1$. Moreover we put
$$
\mathscr{H}=ie^{-q}(H+1) \in i \mathbb R,
\ \
\nu=e^{-q/2}\lambda,
\ \
\mathscr{Q}=-iQ.
$$
 Note that $H = \tanh q$ implies, $\mathscr H = i e^{-q} (H+1) = -i e^q
 (H-1)$, and the Lax pair
$$
U=\Phi^{-1}\Phi_{z},\  \
V=\Phi^{-1}\Phi_{\bar z}
$$
 is given explicitly  by the matrices
\begin{equation}\label{LaxH<1}
U=\left(
\begin{array}{cc}
u_{z}/4 & -\frac{1}{2}\nu^{-1}\mathscr{H}e^{u/2}
\\
\nu^{-1}\mathscr{Q} e^{-u/2} & -u_{z}/4
\end{array}
\right),
\ \
V=
\left(
\begin{array}{cc}
-u_{\bar z}/4 & 
-\nu \overline{\mathscr{Q}}e^{-u/2}
\\
\frac{1}{2}\nu \mathscr{H}e^{u/2}
&
u_{\bar z}/4
\end{array}
\right).
\end{equation}
 Considering, as in section \ref{sc:CMCH>1}, everything as holomorphic 
 expressions in $\nu \in \C^{\times}$ it is straightforward to 
 check that the $1$-form $\alpha = U \mathrm{d} z + V \mathrm{d} \bar z$ is fixed by 
 the following automorphism of the loop algebra:
\begin{equation}\label{eq:invc4}
\tau_4:g(\nu)\longmapsto -\mathrm{Ad}
(\mathscr{R})
\{g(i/\bar{\nu})\}^{*},\ \ 
\mathscr{R}=\left(
\begin{array}{cc}
1/\sqrt{i} & 0\\
0 & \sqrt{i}
\end{array}
\right).
\end{equation}
 This automorphism is said to be of type $C_4$ (almost compact
 automorphism of the fourth kind), \cite{Kobayashi}. The Maurer-Cartan form
 $\alpha=\Phi^{-1}\mathrm{d}\Phi$ has the decomposition
\begin{equation}\label{eq:Maurernu}
\alpha=\nu^{-1}\alpha_{-1}+\alpha_{0}+\nu\alpha_{1}.
\end{equation}
 Moreover, $\tau_4 (\alpha) = \alpha$ translates into 
$$
\overline{\alpha_{0}}=-\alpha_{0},
\ \
\alpha_{-1}=
i\mathrm{Ad}(\mathscr{R})(\alpha_{1})^{*}.
$$
 The mapping $\Phi$ takes values in the twisted loop group
$$
\Lambda \mathrm{SL}_{2}\C_{\sigma,\tau_4}:=
\left\{
g:S^1\to \mathrm{SL}_{2}\C
\
\vert
\
g(-\nu)=\sigma g(\nu),
\ \
\tau_4 (g)(\nu)=g(\nu)\>
\right\}, 
$$
 where $\sigma$ is defined in \eqref{eq:sigma} and 
 \begin{equation}\label{eq:tau4}
 \tau_4 (g)(\nu) = \mathrm{Ad} (\mathscr{R}) \{g(i/\bar{\nu})\}^{*-1}.
 \end{equation}
 Moreover, the first formula in Proposition \ref{prop:Sym} can be
 reinterpreted as 
\begin{equation}\label{eq:SymH<1}
f_{\nu}:=\Phi \left(
\begin{array}{cc}
e^{-q/2} & 0\\
0 & e^{q/2}
\end{array}
\right)
\Phi^{*},
\end{equation}
 thus reproducing the given CMC immersion in $\mathbb H^3$ with $0\leq H
 = \tanh q<1$. Its unit normal can be written in the form
\begin{equation}\label{eq:normalH<1}
n_{\nu}=
\Phi
\left(
\begin{array}{cc}
e^{-q/2} & 0\\
0 & -e^{q/2}
\end{array}
\right)
\Phi^{*}.
\end{equation}
Conversely, the following result holds:
\begin{proposition}\label{prop:fundamentalH<1}
 Let $\mathscr H$ be a purely imaginary constant, $u$ a real function and
 $Q$ a holomorphic function on the simply connected domain $\mathbb D \subset \C$. 
 Let $\nu$ be a complex parameter.
 Assume that the differential $1$-form $\alpha = U\mathrm{d}z + V \mathrm{d}\bar z$ is
 integrable and let $\Phi$ denote a solution to $ \Phi^{-1} \mathrm{d} \Phi =
 \alpha$. Then \eqref{eq:SymH<1} defines for $\nu$ of absolute value
 $e^{-q/2}$ a CMC immersion into $\mathbb H^3$ with mean
 curvature $ H = \tanh q$ and normal as defined in \eqref{eq:normalH<1}.
\end{proposition}
\begin{remark}{\rm 
 As pointed out above,
 the Maurer-Cartan form
 $\alpha=\lambda^{-1}\alpha_{-1}+\alpha_{0}+\lambda\alpha_{1}$ is a type
 $C_4$ real form of the complex CMC surface equation \cite{Kobayashi}.
 But it does not correspond naturally to a CMC surface in $\mathbb E^3$, 
 since $\mathscr H$ is not real. 
}
\end{remark}

\subsection{}
 As indicated in the introduction
 the extended frame \eqref{eq:extendedH<1} 
 of a CMC surface $f$ with mean curvature $H=\tanh q$ 
 can also be considered as the extended frame of a minimal surface in 
 the hyperbolic $3$-space of sectional curvature $-1/\cosh^2 q$. 
 For simplicity, this paper primarily considers surfaces in $\mathbb H^3(-1)$.
 However, the formalism can easily be adjusted to fit surfaces in 
 $\mathbb H^3(c)$ with $c<0$.
 It suffices to ``scale'' a given surface, and thus $\mathbb H^3(-1)$ inside 
 $\mathrm{Her}_2 \mathbb C$, by the factor $1/\sqrt{|c|}$. The radial 
 deformation of the loop parameter $\lambda$ on unit circle
 to a radius $r$ circle changes the mean curvature for a CMC surface, 
 which is given by conjugation of a diagonal matrix to the extended frame.
 Combining the scaling and the radial deformation, 
 the Lawson correspondence for CMC surfaces 
 in $\mathbb H^3(c)$ is obtained.

 More precisely, let $\tilde{\Phi}=\tilde{\Phi}_{\lambda}$ be the
 $\Lambda\mathrm{SL}_{2}\mathbb{C}_{\sigma}$-valued map
 defined by (\ref{eq:tildephi}) which frames the associated family of
 a CMC surface $f:\mathbb{D}\to \mathbb{H}^3(-1)$ with mean curvature $H$. 
 Then for any real number $q$,
 \begin{equation}\label{defomedimmersion}
 f_{\lambda} = \left. \frac{1}{\cosh q - H \sinh q} 
 \tilde \Phi_{\lambda}
 \begin{pmatrix}
 e^{q/2} & 0 \\ 
 0 & e^{-q/2} 
 \end{pmatrix}
 \tilde \Phi_{\lambda}^*\right|_{\lambda=e^{q/2}}
 \end{equation}
 defines a CMC surface of mean curvature $H_{\lambda} = H \cosh q -\sinh q$
 in the hyperbolic space $\mathbb{H}^{3}(K_{\lambda})$ of sectional curvature 
 $K_{\lambda} = - (\cosh q - H \sinh q)^2$. 
 The surface $f_{\lambda}$ has the same metric and 
 the same Hopf differential as $f$. 
 By definition, $H_{\lambda}^2+K_{\lambda} = H^2-1$. 
 Thus $f_\lambda$ is a Lawson correspondent of 
 $f:\mathbb{D}\to \mathbb{H}^{3}(-1)$ in $\mathbb{H}^{3}(K_\lambda)$.

 Now we consider a CMC surface $f$ with mean curvature 
 $H=\tanh q$, then the Lawson correspondent $f_\lambda$  has the mean 
 curvature $H_{\lambda} =0$ and thus is a minimal surface.
 Note that the sectional curvature of 
 the ambient space is $K_\lambda=-1/\cosh^{2}q$ and the extended 
 frame of $f_{\lambda}$ is given as in \eqref{eq:extendedH<1}.
 In section \ref{subsc:primitiveminimal} we will show that 
 the minimality of the surface in $\mathbb H^3(c), c<0$ and primitivity of the 
 Gauss map of the surface are equivalent.
\section{$4$-symmetric structure of the unit tangent sphere bundle}\label{sc:4symUnittangent}

\subsection{}\label{sc:4symmetric}
 As we have seen in the preceding section, every CMC surface with
 $0 \leq H <1$ admits a loop group valued map $\Phi$ which is fixed under the
 type $C_4$ automorphism $\tau_4$.

 In this section we study the  automorphism $\tau$ of $\mathfrak{sl}_2 \C$ 
 which is obtained by first extending $\tau_4$ to the untwisted 
 loop algebra $\Lambda \mathfrak{sl}_2 \C$ using the formula 
 \eqref{eq:invc4} and then
 restricting it to $\mathfrak{sl}_2 \C \subset \Lambda
 \mathfrak{sl}_2 \C$.
 The automorphism $\tau$ is given on
 $\mathfrak{g}=\mathfrak{sl}_{2}\C$ by the formula
\begin{equation}\label{6.1}
\tau(X)=-\mathrm{Ad}(\mathscr{R})X^{*}, 
\end{equation} 
 where $\mathscr{R}$ is defined in \eqref{eq:invc4}.
 More explicitly, 
\begin{equation}\label{eq:tau}
\tau
\left(
\begin{array}{cc}
a & b\\
c & -a
\end{array}
\right)=
\left(
\begin{array}{cc}
-\bar{a} & i\bar{c}\\
-i\bar{b} & \bar{a}
\end{array}
\right).
\end{equation}

 It is easy to see that $\tau$ is of order $4$. The eigenspace
 decomposition of  the \textit{complexified Lie algebra}
 $\mathfrak{g}^{\C}$ with respect to $\tau$ is given by
$$
\mathfrak{g}^{\C}=\mathfrak{g}_{0}^{\C}\oplus
\mathfrak{g}_{1}^{\C}\oplus
\mathfrak{g}_{2}^{\C}\oplus
\mathfrak{g}_{3}^{\C},
$$
 where $\mathfrak{g}_{k}^{\C}$ is the eigenspace corresponding to the
 eigenvalue $i^k$. 
 Note that the complexified Lie algebra 
 $\mathfrak{g}^{\C}=(\mathfrak{sl}_{2}{\C})^{\C}$ is realized as
 $\mathfrak{g}\times\mathfrak{g}$. This construction can be described
 as follows: 
 
 Consider the map 
$$
 \iota: X \in \mathfrak{sl}_2 \C  \longmapsto (X,  \; \bar X) \in \mathfrak{sl}_2 \C \times
 \mathfrak{sl}_2 \C\;. 
$$
 This is an injective homomorphism of the real Lie algebra
 $\mathfrak{sl}_2 \C$ into $\mathfrak{sl}_2 \C \times
 \mathfrak{sl}_2 \C$. Since  $\iota(\mathfrak{sl}_2 \C) \cap 
 i \iota(\mathfrak{sl}_2 \C) = \{0\}$,
 the image $\iota(\mathfrak{sl}_2 \C)$ is a real form of 
 $\mathfrak{sl}_2 \C \times \mathfrak{sl}_2 \C$.
 The latter Lie algebra carries the natural complex structure 
 given by multiplying a complex number to each of the two factors.
 By transporting $\tau$ via $\iota$ we obtain
\begin{enumerate}  
\item[(i)] $\tau$ acts on $\iota(\mathfrak{sl}_2 \C)$ as $\tau(X, \bar
     X) = (\tau(X), \;  \overline{\tau(X)})$.
\end{enumerate}
 We now define the complex linear extension $\hat \tau$  to
 $\mathfrak{sl}_2 \C \times \mathfrak{sl}_2 \C$ as
 follows:
\begin{enumerate}
\item[(ii)] $\hat \tau(X,Y) = ( \tau(\bar Y),\; \overline{\tau(X)})$.
\end{enumerate}
 It is straightforward to show that, indeed, $\hat \tau$ is complex linear
 relative to $i$ acting on the first and the second factor equally by
 multiplication.
 
 Moreover, on $\iota (X)$ the new $\hat \tau$ acts like {\rm (i)}.
 Thus {\rm (ii)} is the complex linear extension of the original
 $\tau$ on $\mathfrak{sl}_2 \C$ to the product 
 $\mathfrak{sl}_2 \C \times \mathfrak{sl}_2 \C$.
 
 Since $\tau$ is an automorphism of order $4$, also $\hat \tau$
 is an automorphism of order $4$ and we have
\begin{lemma}\label{grading}
$$
[\mathfrak{g}_{k}^{\C},\mathfrak{g}_{l}^{\C}] \subset
\mathfrak{g}_{k+l}^{\C}\ \ (\mathrm{mod} \>4).
$$
\end{lemma}

 The automorphism has two \textit{real eigenvalues} $i^{0}=1$ and
 $i^{2}=-1$. The corresponding eigenspaces $\mathfrak g_0^{\C}$
 and $\mathfrak g_2^{\C}$ of
 $\mathfrak g^{\C} \cong \mathfrak g \times \mathfrak g$ are
 computed explicitly as
\begin{align*}
\mathfrak{g}_0^{\C} = 
\left\{
\left(
\begin{pmatrix}
x_1 & 0 \\ 0 & -x_1
\end{pmatrix},
\begin{pmatrix}
-x_1 & 0 \\ 0 & x_1
\end{pmatrix}
\right)
\right\},\;\;\;
\mathfrak{g}_2^{\C} = 
\left\{
\left(
\begin{pmatrix}
x_1 & 0 \\ 0 & -x_1
\end{pmatrix},
\begin{pmatrix}
x_1 & 0 \\ 0 & -x_1
\end{pmatrix}
\right)
\right\},
\end{align*}
 with $x_1 \in \C$.
 The corresponding real subspaces of $\mathfrak{sl}_2 \C$ are
 $$
 \mathfrak g \cap \mathfrak {g}_0^{\C} =
 \mathfrak h = \mathbb R  (i \Vec{e}_1)
 \; \mbox{and}\;\;  \mathfrak g \cap \mathfrak {g}_2^{\C} =
 \mathbb R (\Vec{e}_1).
 $$
 Moreover, the eigenspaces $\mathfrak g_1^{\C}$ and $\mathfrak g_3^{\C}$ 
 for eigenvalues $i$ and $-i$ are computed explicitly as
\begin{align*}
\mathfrak{g}_1^{\C} =
\left\{
\left(
\begin{pmatrix}
0 & x_2 \\ x_3 & 0
\end{pmatrix},
\begin{pmatrix}
0 & -x_3 \\ x_2 & 0
\end{pmatrix}
\right)
\right\}, \;\;\;
\mathfrak{g}_3^{\C} =
\left\{
\left(
\begin{pmatrix}
0 & x_2 \\ x_3 & 0
\end{pmatrix},
\begin{pmatrix}
0 & x_3 \\ -x_2 & 0
\end{pmatrix}
\right)
\right\},
\end{align*}
 with $x_2, x_3 \in \C$.
 The automorphism $\tau$ defines the semi-Riemannian $4$-symmetric
 space $\mathrm{SL}_2 \C/\mathrm{U}_1$. The space is isomorphic with
 the unit tangent sphere bundle $\mathrm{U}\mathbb H^3$ as shown in
 section \ref{sc:tangentspherebundle}. The complexified tangent space
 $\mathfrak{p}^{\C}=
 (T_{(\Vec{e}_0,\Vec{e}_1)}\mathrm{U}\mathbb{H}^3)^{\C}$ is given
 by
$$
\mathfrak{p}^{\C}=
\mathfrak{g}_{1}^{\C}\oplus\mathfrak{g}_{2}^{\C}
\oplus\mathfrak{g}_{3}^{\C}.
$$
 Comparing this with the fibration $\pi_{1}: \mathrm{U}\mathbb{H}^3\to
 \mathrm{Geo}(\mathbb{H}^3)$ discussed in section
 \ref{sc:tangentspherebundle}, we have
$$
\mathcal{H}_{1}^{\C}=\mathfrak{g}_{1}^{\C}\oplus\mathfrak{g}_{3}^{\C},
\ \
\mathcal{V}_{1}^{\C}=\mathfrak{g}_{2}^{\C}.
$$

\subsection{}\label{sc:contactmfd}
 Here we recall the notion of a contact manifold.
\begin{definition}{\rm
 A 1-form $\omega$ on a manifold $L$ of dimension $2n-1$ is said to be a
 \textit{contact form} if $(\mathrm{d}\omega)^{n-1} \wedge \omega \not=0$ on
 $M$. A hyperplane field $\mathcal{D}\subset TL$ on $L$ is called a
 \textit{contact structure} if for any point $p\in L$, there exists a
 contact form $\omega$ defined on a neighborhood $U_p$ of $p$ such that
 $\mathrm{Ker}\> \omega=\mathcal D$ on $U_p$. 
}
\end{definition}
 A $(2n-1)$-manifold $L$ with a contact structure  $\mathcal{D}$ is
 called a \textit{contact manifold}. If a contact manifold
 $(L,\mathcal{D})$ admits a globally defined contact form $\omega$ which
 annihilates $\mathcal{D}$,  \textit{i.e.},
 $\mathrm{Ker}\>\omega=\mathcal{D}$, then $(L,\mathcal{D})$ is said to
 be a contact manifold in the \textit{strict sense}.
\begin{definition}
{\rm
Let $M^n$ be an $n$-manifold and $F:M \to L$ a smooth map into
a contact $(2n-1)$-manifold. Then $F$ is said to be \textit{Legendre} if 
$\mathrm{d} F(TM)\subset \mathcal{D}$.
}
\end{definition}
 In particular, if $L$ admits a global contact form $\omega$, 
 then $F$ is Legendre if and only if  $F^{*}\omega=0$.

 Now let $(N^n,g)$ be a Riemannian $n$-manifold. Then its unit tangent
 sphere bundle $\mathrm{U}N$ admits a canonical contact structure. In case
 $N=\mathbb{H}^3$, one can check that the canonical contact structure of
 $\mathbb{H}^3$ is given by
$$
 \mathcal{D}=\mathcal{H}_{1}=(\mathfrak{g}_{1}^{\C}
 \oplus\mathfrak{g}_{3}^{\C}) \cap \mathfrak g.
$$
 By using this fact and results of section \ref{sc:tangentGauss} we
 obtain:

\begin{proposition}\label{prop:Legendre}
 Let $F=(f, n):\mathbb{D}\to\mathrm{U}\mathbb{H}^3$ be a smooth map with
 frame $\Phi:\mathbb{D}\to \mathrm{SL}_{2}\C$,
 \textit{i.e.}, $\Phi$ is a map satisfying 
$$
F=(\Phi \Phi^*, \Phi \Vec{e}_1\Phi^{*}).
$$
 Denote by $\alpha=\Phi^{-1}\mathrm{d}\Phi$
 the pull-back of the Maurer-Cartan form by $\Phi$ and decompose
 $(\alpha,\overline{\alpha})$ 
as
$$
(\alpha,\overline{\alpha})=
\alpha_{0}^\mathbb{C}
+
\alpha_{1}^\mathbb{C}
+
\alpha_{2}^\mathbb{C}
+
\alpha_{3}^\mathbb{C}
$$
according to the eigenspace decomposition with respect to 
$\tau$. Then $\alpha_{2}^{\mathbb{C}}$ is given by
$\alpha_{2}^{\mathbb{C}}=(\alpha_2,\overline{\alpha_2})$ with
$$
\alpha_{2}=\frac{1}{2}F^{*}\omega\>\Vec{e}_{1}.
$$
 Thus $F$ is Legendre if and only if $\alpha_2=0$.
 In particular, if $F=(f,n)$ is the Gauss map of a conformal
 immersion $f:\mathbb{D}\to\mathbb{H}^3$, then 
$$
\alpha_{2}=\frac{1}{2}F^{*}\omega\>\Vec{e}_{1} = \frac{1}{2}\langle \mathrm{d}f,n\rangle \>\Vec{e}_{1}=0.
$$
\end{proposition}
%
\subsection{} 
 Let $F: \mathbb{D}\to G/H = \mathrm{SL}_2 \C/\mathrm{U}_1$ be a
 Legendre map with frame $\Phi$. 
 Then we have the eigenspace decomposition
 of $\alpha=\Phi^{-1}\mathrm{d}\Phi$.
$$
(\alpha,\overline{\alpha})=
\alpha_{0}^{\mathbb{C}}+
\alpha_{1}^{\mathbb{C}}+
\alpha_{3}^{\mathbb{C}},
\ \ 
\alpha_{0}^{\mathbb{C}}=
(\alpha_0,\overline{\alpha_0}).
$$
 On the other hand, we have the decomposition
 $\alpha=\alpha_{\mathfrak{h}}+
 \alpha_{\mathfrak{p}}$. 
 We denote the first component of $\alpha_j^{\C}$ by $\alpha_j$.
 Comparing these decompositions, we get
\begin{equation}\label{eq:alphaLegendre}
\alpha=\alpha_{\mathfrak{h}}+
\alpha_{\mathfrak{p}},
\ \
\alpha_{\mathfrak{h}} = \alpha_0, \
\alpha_{\mathfrak{p}} =
\alpha_{\mathfrak{p}}^{\prime}+
\alpha_{\mathfrak{p}}^{\prime\prime}=
\alpha_{1}+\alpha_{3}. 
\end{equation}
 We express the type-decompositions of $\alpha_{1}$ and $\alpha_{3}$ 
 with respect to the conformal structure of $\mathbb{D}$ as
$$
\alpha_{1}=
\alpha_{1}^{\prime}+
\alpha_{1}^{\prime \prime},
\ \
\alpha_{3}=
\alpha_{3}^{\prime}+
\alpha_{3}^{\prime \prime}.
$$
Then from Lemma \ref{grading} and the integrability of $\alpha$, we derive 
$$
[\alpha^{\prime}_{\mathfrak{p}}
\wedge 
\alpha^{\prime\prime}_{\mathfrak{p}}
]_{\mathfrak{p}}=
[
(\alpha_{1}^{\prime}
+
\alpha_{3}^{\prime})
\wedge 
(
\alpha_{1}^{\prime\prime}
+
\alpha_{3}^{\prime \prime})
]
_{\mathfrak{p}}=
[\alpha_{1}^{\prime}
\wedge 
\alpha_{1}^{\prime\prime}]
+
[\alpha_{3}^{\prime}
\wedge 
\alpha_{3}^{\prime\prime}]
$$ 
 and 
$$
[\alpha_1 \wedge \alpha_1]+ 
[\alpha_3 \wedge \alpha_3]= 0.
$$
 Noting 
 $[\alpha_1 \wedge \alpha_1] = 2 [\alpha^{\prime \prime}_1 \wedge \alpha^{\prime}_1]$  and 
 $[\alpha_3 \wedge \alpha_3] =2 [\alpha^{\prime\prime}_3 \wedge \alpha^{\prime}_3]$, 
 we conclude that
$$
[\alpha^{\prime}_{\mathfrak{p}}
\wedge 
\alpha^{\prime\prime}_{\mathfrak{p}}]_{\mathfrak{p}}
=
[\alpha_{1}^{\prime}
\wedge 
\alpha_{1}^{\prime\prime}]
+
[\alpha_{3}^{\prime}
\wedge 
\alpha_{3}^{\prime\prime}]
=
\frac{1}{2} [\alpha_{1}
\wedge 
\alpha_{1}]
+
\frac{1}{2} [\alpha_{3}
\wedge 
\alpha_{3}]
=
0.
$$
 Now we arrive at the following 
 \textit{zero curvature representation} for Legendre harmonic maps.
\begin{proposition}\label{ZCR-Legendre}
 Let $F:\mathbb{D}\to \mathrm{SL}_2 \C/\mathrm{U}_1$ be a
 Legendre harmonic map with frame
 $\Phi:\mathbb{D}\to \mathrm{SL}_2 \C$.
 Then $\alpha_{\lambda}= \Phi^{-1}\mathrm{d}\Phi= 
 \alpha_{\mathfrak{h}}+ \lambda^{-1} \alpha_{\mathfrak{p}}^{\prime} + \lambda
 \alpha_{\mathfrak{p}}^{\prime \prime} $
satisfies
$$
\mathrm{d}\alpha_{\lambda}+\frac{1}{2}[\alpha_\lambda
 \wedge \alpha_\lambda]=0
$$
 for all $\lambda \in \C^{\times}$.
 Here $\alpha_{\mathfrak{h}}$ and $\alpha_{\mathfrak{p}}$ are
 defined in \eqref{eq:alphaLegendre} and $
 \alpha_{\mathfrak{p}}^{\prime}$ (resp. $
 \alpha_{\mathfrak{p}}^{\prime \prime}$) is the $(1, 0)$-part
 (resp. $(0,1)$-part)  of  $\alpha_{\mathfrak{p}}$.
\end{proposition}

\subsection{} 
 The square $\tau^{2}$ of $\tau$ is an involutive 
 automorphism of  $\mathfrak g$. The $1$-eigenspace and
 $(-1)$-eigenspace of $\tau^2$ on $\mathfrak g$
 are
$$
\mathfrak{g} \cap (\mathfrak{g}_0^{\C} \oplus
\mathfrak{g}_2^{\C}) = \mathfrak{d}
\ \mbox{and}\
\mathfrak{g} \cap (\mathfrak{g}_{1}^{\C}\oplus\mathfrak{g}_{3}^{\C})
 =\mathfrak{q},
$$
 where $\mathfrak{d}$ and $\mathfrak{q}$ are defined in \eqref{eq:d} and
 \eqref{eq:q}, respectively. Hence $(G,\tau^2)$ defines the
 semi-Riemannian symmetric space 
$$
 G/D=\mathrm{SL}_{2}\C/\C^{\times}
    =\mathrm{Gr}_{1,1}(\mathbb{E}^{1,3}).
$$
 The $4$-symmetric space $G/H = \mathrm{SL}_2 \C / \mathrm{U}_1$
 is a fiber bundle over $G/D$ with standard fiber
$$
D/H=\C^{\times}/\mathrm{U}_1=\mathbb{R}_{>0}.
$$
 Now let $F:\mathbb{D}\to G/H$ be a Legendre harmonic map with frame
 $\Phi: \mathbb D \to G$. Since $F$ is Legendre,
%
 Integrating $\Phi_{\lambda}^{-1}
 \mathrm{d}\Phi_{\lambda}=\alpha_{\lambda}$, 
 we get the associated
 family $\{F_{\lambda}\}$ of $F$.

 Decompose $\alpha$ according to the Lie algebra decomposition
 $\mathfrak{g}=\mathfrak{d}\oplus \mathfrak{q}$:
$$
\alpha_{\mathfrak{d}}=
\alpha_{0}+\alpha_{2} = \alpha_0,
\ \
\alpha_{\mathfrak{q}}=\alpha_{1}+\alpha_{3}.
$$
 Then the decomposition above can be rephrased as 
$$
\alpha_{\lambda}
=\alpha_{\mathfrak{d}}
+\lambda^{-1}
\alpha^{\prime}_{\mathfrak q}
+\lambda
\alpha^{\prime\prime}_{\mathfrak q}.
$$
 This formula implies that the projected map $\mathcal{G}=\pi_{1}\circ
 F$ is harmonic and 
$$
 \mathcal{G}_{\lambda}:=\pi_{1}\circ F_{\lambda}
$$
 gives the associated family $\{\mathcal{G}_\lambda\}$ of $\mathcal{G}$,
 where $\pi_1$ is the natural projection $\pi_1 : G/H \to G/D$.

\begin{remark}
 {\rm Proposition \ref{ZCR-Legendre} is valid for any horizontal harmonic maps 
 into semi-Riemannian $4$-symmetric spaces. More precisely, let $(G/H,\tau)$ 
 be a semi-Riemannian $4$-symmetric space with semi-Riemannian homogeneous 
 projection $\pi_1:G/H\to G/D$ onto the semi-Riemannian symmetric space 
 $(G/D,\tau^2)$. 
 Take a map $F:\mathbb{D}\to G/H$ which is horizontal with respect to $\pi_1$ 
 and let $\Phi:\mathbb{D}\to G$ be its frame.
 Decompose $\alpha=\Phi^{-1}\mathrm{d}\Phi$ as $\alpha=\sum_{j=0}^{3}\alpha_j$ 
 according to the eigenspace decomposition of $\tau$. Then one can see that 
 $F$ is horizontal if and only if $\alpha_2=0$. In addition, the harmonicity 
 of $F$ is equivalent to the flatness of the connections 
 $\mathrm{d}+\alpha_\lambda$, where $\alpha_\lambda$ is defined 
 as in Proposition \ref{ZCR-Legendre}.
}
\end{remark}
\subsection{}\label{subsc:primitiveminimal}
 In this section, we prove the following characterization of minimal surfaces.
\begin{proposition}\label{prop:prm}
  Let $\mathbb H^3(c)$ be the hyperbolic $3$-space of sectional curvature $c<0$.
 Then the unit tangent sphere bundle 
 $\mathrm{U} \mathbb H^3(c)$ is a $4$-symmetric space.
 Moreover, a surface $f: \mathbb D \to \mathbb H^3(c)$ is minimal 
 if and only if its Gauss map is a primitive map with 
 respect to the $4$-symmetric structure of $\mathrm{U}\mathbb H^3(c)$.
\end{proposition}
 {\bf Proof}. It is clear that $\mathrm{U} \mathbb H^3(c)$ is 
 a $4$-symmetric space with respect to $\tau$ defined in \eqref{eq:tau}. 
 Let $f:M\to \mathbb{H}^3(c)$ be a conformal 
 immersion with unit normal $n$.
 Take a simply connected coordinate domain $(\mathbb{D},z)\subset M$ and denote 
 by $\Psi$ the coordinate frame defined on $\mathbb{D}$. Let $\hat{\Phi}$ be 
 a lift of $\Psi$ to $\mathrm{SL}_{2}\mathbb{C}$ as in (\ref{eq:hatframe}).
 Then the Maurer-Cartan form $\hat{\alpha}=\hat{\Phi}^{-1}
 \mathrm{d}\hat{\Phi}$ is given by  (\ref{eq:hatMau1}) 
 with $H+1$ and $H-1$ replaced by $H+ \sqrt{|c|}$ and $H-\sqrt{|c|}$,  
 respectively.
 Now we decompose $(\hat \alpha, \overline{\hat{\alpha}})$
 according to the eigenspace decomposition of $\tau$. Then we have
$$
(\hat \alpha, \overline{\hat{\alpha}})=
\hat{\alpha}^{\mathbb{C}}_{0}+
\hat{\alpha}^{\mathbb{C}}_{1}+
\hat{\alpha}^{\mathbb{C}}_{2}+
\hat{\alpha}^{\mathbb{C}}_{3},
$$
where
\begin{align*}
\hat{\alpha}^{\mathbb{C}}_{0}
&=
\left(
\begin{pmatrix}x_1 & 0 \\ 0 & -x_1 \end{pmatrix}, 
\begin{pmatrix}-x_1 & 0 \\ 0 & x_1 \end{pmatrix}
\right), \;\;
\hat{\alpha}^{\mathbb{C}}_{2}
=(\Vec{0},\Vec{0}), \\
\hat{\alpha}^{\mathbb{C}}_{1}
&=
\left(
\begin{pmatrix}0 & x_2 \\ x_3 & 0 \end{pmatrix}, 
\begin{pmatrix}0  & -x_3 \\ x_2 & 0 \end{pmatrix}
\right), \;\;
\hat{\alpha}^{\mathbb{C}}_{3}
=
\left(
\begin{pmatrix}0 & x_4 \\ x_5 & 0 \end{pmatrix}, 
\begin{pmatrix}0  & x_5 \\ -x_4 & 0 \end{pmatrix}
\right)
\end{align*}
 with $x_1 = \frac{1}{4}(u_{z}\mathrm{d}z-u_{\bar z}\mathrm{d}\bar{z})$, 
 $x_2 =  \frac{\sqrt{|c|}}{2}e^{u/2}\mathrm{d}z$, 
 $x_3 =-Qe^{-u/2}\mathrm{d}z-\frac{H}{2}e^{u/2}\mathrm{d}\bar{z}$, 
 $x_4 = \frac{H}{2}e^{u/2}\mathrm{d}z+\bar{Q}e^{-u/2}\mathrm{d}\bar{z}$ 
 and $x_5 = \frac{\sqrt{|c|}}{2}e^{u/2}\mathrm{d}\bar{z}$.
 From these equations, we deduce that $f$ is minimal if and only if its 
 Gauss map $F$ is primitive. $\Box$

\begin{remark}
\mbox{}
 {\rm
\begin{enumerate}
\item The unit tangent sphere bundle $\mathrm{U}\mathbb{H}^3$ is 
 the twistor $CR$-manifold of $\mathbb{H}^3$ in the sense of \cite{LeBrun}. 
 There exist two standard $f$-structures $J_1$ and $J_2$ on $\mathrm{U}\mathbb{H}^3$, 
 \textit{i.e.}, endomorphism fields $J$ on  $\mathrm{U}\mathbb{H}^3$ 
 such that $J^3+J=0$.
 One can see that a map $F:M\to \mathrm{U}\mathbb{H}^3$
 from a Riemann surface to $\mathrm{U}\mathbb{H}^3$ 
 is $J_2$-holomorphic if and only if $F$ is a primitive map. 
 On the other hand, 
 for a conformal immersion $f:M\to \mathbb{H}^3$,
 its Gauss map is $J_1$-holomorphic if and only if 
 $f$ is totally umbilical (see \cite[Theorem 7.1]{Salamon}). 

\item  Setting $\tau (g) = {\rm Ad} \mathscr R (g^*)^{-1}$ we 
 obtain an automorphism of $G = {\rm SL}_2 \C$, the differential 
 of which coincides with $\tau$ as given in \eqref{6.1}. 
 By abuse of language we will
 use the same notation for the group level and for the Lie algebra level
 as well as for the corresponding complexified objects.

 \item In \cite{DMPW},  the original loop group approach \cite{DPW} 
 was extended to include primitive harmonic maps into compact $k$-symmetric 
 spaces. In our case the symmetric space under consideration is non-compact.
 This has far-reaching consequences. We have therefore included in sections 
 \ref{sc:Potential} and \ref{sc:DPW}
 a brief description of the corresponding technical details.
 In particular, the Iwasawa decomposition has not only one, but two open cells
 (the union of which is dense).
 Implications of this can already be seen in the examples presented in section 
 \ref{sc:Ex}.
\end{enumerate}
}
\end{remark}

\section{Potentials}\label{sc:Potential}
\subsection{}
 We recall loop groups and the Birkhoff decomposition.
 The twisted loop group is defined as
 \begin{equation}\label{eq:loopgroups}
 \Lambda{\rm SL}_2 \C_\sigma = \left\{ g:S^1  \to {\rm
 SL}_2 \C\; \left| \right. \; g \;\mbox{is continuous and}\;
 g(-\lambda) =  \sigma g(\lambda)  \; \right \} \; ,
 \end{equation}
 where $\sigma$ is defined in \eqref{eq:sigma}.
 More strictly, we assume that the coefficients of all $g \in
 \Lambda{\rm SL}_2 \C_\sigma$ are in the {\it Wiener algebra}
 $\mathcal A = \left\{ f(\lambda) = \sum_{n \in \mathbb Z} f_n
 \lambda^n\;  |\;  S^1 \to \C \;\; ; \;\; \sum_{n
 \in \mathbb Z}|f_n| < \infty \right\}$. 
  The Wiener algebra is a Banach algebra relative to the norm 
 $\| f\| = \sum |f_n|$, and $\mathcal A$
 consists of continuous functions. Thus 
 $\Lambda{\rm SL}_2 \C_\sigma$ is a Banach Lie group. 
 We denote the Lie algebra of
 $\Lambda{\rm SL}_2 \C_\sigma$ by $\Lambda\mathfrak{sl}_2
 \C_\sigma$, which consists of maps  $g : S^1 \to
 \mathfrak{sl}_2 \C$.

 We will need to consider two subgroups of $\Lambda{\rm SL}_2 \C_{\sigma}$, 
 the {\it twisted plus loop group} and the {\it minus loop group} 
 as follows:
 Let $\boldsymbol B$ a subgroup of ${\rm SL}_2 \C$.
 Let $ \Lambda_{B}^+ {\rm SL}_2 \C_\sigma$
 be the group of maps into $\Lambda {\rm SL}_2 \C_{\sigma}$
 which can be extended holomorphically to $D = \{\lambda \in \C
 \;|\;  |\lambda| < 1\}$ and which take values
 in $\boldsymbol B$ at $\lambda =0$. Similarly, let $ \Lambda_{B}^-
 {\rm SL}_2 \C_\sigma$ be the group of maps into $\Lambda {\rm
 SL}_2 \C_{\sigma}$ which can be extended holomorphically to
 $E =\{\lambda \in \C  \;|\; 1 <
 |\lambda| \} \cup \{\infty\}$ and take values
 in $\boldsymbol B$ at $\lambda =\infty$. If $\boldsymbol B = \{\rm Id\}$ 
 we write the subscript $*$  instead of $\boldsymbol B$, if $\boldsymbol B = 
 {\rm SL}_2 \C$ we  abbreviate  $\Lambda_{B}^+ {\rm SL}_2 \C_\sigma$ and 
 $\Lambda_{B}^- {\rm SL}_2 \C_\sigma$ by
 $\Lambda^+ {\rm SL}_2 \C_\sigma $ and 
 $\Lambda^- {\rm SL}_2 \C_\sigma$, respectively. 

 It is clear that the loop groups 
 $\Lambda {\rm SL}_2 \C_{\sigma, \tau_j}, \;j = 3, 4,$
 defined in sections \ref{sc:CMCH>1} and \ref{sc:CMCH<1} are also 
 subgroups of $\Lambda {\rm SL}_2 \C_\sigma$.  
 The Lie algebras of  $\Lambda {\rm SL}_2 \C_{\sigma, \tau_j}$
 are denoted by $\Lambda {\mathfrak{sl}}_2 \C_{\sigma, \tau_j}, \;j = 3, 4,$ 
 and those are subalgebras of $\Lambda {\mathfrak{sl}}_2 \C_{\sigma}$.
 \begin{theorem}[Birkhoff decomposition
  \cite{LoopGroup}]\label{thm:Birkhoff}
 The maps
\begin{equation*}
 \Lambda_{*}^- {\rm SL}_2 \C_\sigma \times \Lambda^+
 {\rm SL}_2 \C_\sigma \to \Lambda {\rm SL}_2 \C_\sigma
 \;\; \mbox{and}\;\;\Lambda_{*}^+ {\rm SL}_2 \C_\sigma \times \Lambda^-
 {\rm SL}_2 \C_\sigma \to \Lambda {\rm SL}_2 \C_\sigma \\
\end{equation*}
 are analytic diffeomorphisms onto the open dense subsets 
 $\Lambda_{*}^- {\rm SL}_2 \C_\sigma \cdot \Lambda^+
 {\rm SL}_2 \C_\sigma$ and $\Lambda_{*}^+ {\rm SL}_2 \C_\sigma 
 \cdot \Lambda^- {\rm SL}_2 \C_\sigma$ of 
 $\Lambda {\rm SL}_2 \C_\sigma$ respectively.
 The open dense subsets will be called the {\rm left big cell} and 
 the {\rm right big cell} respectively.
\end{theorem}
\begin{remark}
 {\rm In this paper the big cell always means the left big cell.}
\end{remark}

\subsection{}\label{subsc:holopot}
 The holomorphic potential for a CMC surface $f$ in $\mathbb H^3$
 is a $\Lambda \mathfrak{sl}_2 \C_{\sigma}$-valued holomorphic $1$-form
 determined from the extended frame of $f$, which is an analogue of 
 Weierstrass data for a minimal surface in $\mathbb R^3$. 
 The holomorphic potential reproduces the CMC surface $f$ 
 using the generalized Weierstrass
 type representation recalled in section \ref{sc:DPW}.
 
 Let $\mathbb D$ be a simply connected domain in $\C$, and 
 let $\Phi :\mathbb D\to \Lambda {\rm SL}_2 \C_{\sigma, \tau_j}, j = 3, 4,$ 
 the extended frame of a CMC surface $f$ in $\mathbb H^3$ with $H \neq 1$ 
 as defined in \eqref{eq:extendedH>1} or \eqref{eq:extendedH<1}, 
 respectively.
\begin{proposition}\label{prop:delbar}
 There exists a loop $g : \mathbb D \to \Lambda^+ {\rm SL}_2 \C_{\sigma}$, 
 holomorphic for $\lambda \in \C$, such that 
 $\Phi g : \mathbb D \to \Lambda {\rm SL}_2 \mathbb C_{\sigma}$ 
 is holomorphic in $\lambda \in \C$ and $z \in \mathbb D$.
\end{proposition}
 {\bf Proof}.
 The holomorphicity of $\Phi g$ is equivalent to 
 that 
 \begin{equation}\label{eq:delbar}
 g_{\bar z} + Vg =0\;,
 \end{equation}
 where $V$ is defined in \eqref{LaxH>1} or \eqref{LaxH<1} respectively.
 Since $V$ is real analytic in $z$, one can extend $V$ holomorphically to 
 $\mathbb D_{\epsilon} (p_0) \times \mathbb D_{\epsilon} (\bar p_0)$, 
 \textit{i.e.}, 
 there exists, for sufficiently small $\epsilon > 0$, a holomorphic matrix 
 function  
 $\tilde V(z, w) :\mathbb D_{\epsilon} (p_0) \times \mathbb
 D_{\epsilon} (\bar p_0) \to \Lambda^+ {\rm SL}_2 \C_{\sigma}$
 such that  $\tilde  V|_{(z, \bar z)} = V$, where $\mathbb D_{\epsilon}
 (p_0)$  (resp. $\mathbb D_{\epsilon} (\bar p_0)$) denotes the
 $\epsilon$-disk around $p_0 \in \mathbb D$ (resp. $\bar p_0 \in 
 \bar{\mathbb D}$).
 Let us consider the ordinary differential equation:
  \begin{equation}
  \tilde g_{w} + \tilde V \tilde g =0\;, \;\;\tilde g(z, \bar p_0) ={\rm
   Id}.
  \end{equation}
 This equation has, for every fixed $z$,  a unique solution 
 $\tilde g = \tilde g (z, w)$. Setting $g =  \tilde  g(z, \bar z)$, we obtain 
  \begin{equation}
  g_{\bar z} + V g =0\;, \;\; g(p_0, \bar p_0) ={\rm Id}.
  \end{equation}
 Therefore, on every $U_{\alpha}$ of an open cover $(U_{\alpha})$ 
 of $\mathbb D$, there exists 
 a real analytic solution $g_{\alpha}$. On $U_{\alpha} \cap U_{\beta}$, 
 we define 
  \begin{equation}
  h_{\alpha \beta} = g_{\alpha}^{-1} g_{\beta} : U_{\alpha} \cap
   U_{\beta} \to \Lambda^+ {\rm SL}_2 \C_{\sigma}.
  \end{equation}
 It is easy to verify that $h_{\alpha \beta}$ is holomorphic in $z$. 
 Moreover, the $h_{\alpha \beta}$ satisfy the co-cycle condition 
 $h_{\alpha  \beta} h_{\beta \gamma} = h_{\alpha \gamma}$. Thus 
 this data defines a holomorphic $\Lambda^+ {\rm SL}_2 \C_{\sigma}$-principal 
 fiber bundle $P \to \mathbb D$. Since $\mathbb D$ is simply connected,
 by a generalization of Grauert's theorem \cite{Bungart}, the
 holomorphic bundle $P$ is trivial. Thus $h_{\alpha \beta}$ splits, 
 $h_{\alpha \beta} = h_{\alpha} h_{\beta}^{-1}$ on $U_{\alpha} \cap U_{\beta}$,  where $h_{\gamma}: U_{\gamma }\to \Lambda^+ {\rm SL}_2 \C_{\sigma}$ is 
 holomorphic. 
 Therefore $g=  g_{\alpha} h_{\alpha}$  is well-defined on $\mathbb D$
 and \eqref{eq:delbar} holds. This completes the proof.  $\Box$
\begin{corollary}[Existence of a holomorphic potential]
 Let $\Phi :\mathbb D\to \Lambda {\rm SL}_2 \C_{\sigma, \tau_j} 
 (j = 3\;\mbox{or}\;4)$
 be the extended frame for some CMC surface $f$ in $\mathbb H^3$ with $H \neq 1$ 
 as defined in \eqref{eq:extendedH>1} or \eqref{eq:extendedH<1}.
 Then there exist $g: \mathbb D \to \Lambda^+ {\rm SL}_2 \C_{\sigma}$
 and a $\Lambda \mathfrak{sl}_2 \C$-valued $1$-form $\eta$
 such that
 \begin{equation}\label{eq:holopot}
  \eta = (\Phi g)^{-1} \mathrm{d} (\Phi g) 
       = \sum_{j= -1}^{\infty} \lambda^j \eta_j\;,
  \end{equation} 
 where the $\eta_j$ are $\lambda$-independent diagonal (resp. off-diagonal)
 holomorphic $1$-forms if $j$ is even (resp. odd). In particular, 
 the $1$-form $\lambda \eta$ is holomorphic in $z \in \mathbb D$ and 
 $\lambda \in \C$ and the upper right entry of $\eta_{-1}$ does not vanish 
 on $\mathbb D$.
\end{corollary}
 {\bf Proof}. Let $g$ be the loop defined in Proposition \ref{prop:delbar} 
 and set $\eta = (\Phi g)^{-1} \mathrm{d} (\Phi g)$. Then $\eta$
 defines a $\Lambda  \mathfrak{sl}_2 \C_{\sigma}$-valued holomorphic $1$-form
 and can be computed as 
 \begin{equation}\label{eq:holopot2}
 \eta = (g^{-1} U g + g^{-1} g_{z}) \; \mathrm{d}z\;,
 \end{equation}
 where $U$ is defined in \eqref{LaxH>1} and \eqref{LaxH<1}, respectively. 
 Thus $\eta$
 has the form  $\eta = \sum_{j= -1}^{\infty} \lambda^j \eta_j$.
 Moreover, since $\eta$ satisfies the twisting condition, $\eta_j$ is
 diagonal  (resp. off-diagonal) if $j$ is even  (resp. odd). $\Box$
 \begin{remark}
\mbox{}

{\rm 
 \begin{enumerate}
 \item The holomorphic $1$-form $\eta$ in \eqref{eq:holopot} 
       will be called  the {\it holomorphic potential} 
       of the immersion $f$.

 \item Holomorphic potentials are not unique, since the right 
       multiplication of $C$ by some holomorphic loop $C_+ : \mathbb D \to
       \Lambda^+       
       {\rm SL}_2  \C_{\sigma}$ gives  another holomorphic
       potential, \textit{i.e.}, $\tilde \eta = \tilde C^{-1} \mathrm{d} \tilde C$
       is a new holomorphic potential, where $\tilde  C = C C_+$, $C_+$ holomorphic.

 \item The generalization of Grauert's theorem has so far only been 
       proven for non-compact Riemann surfaces. Therefore, at this point 
       we can only infer that  holomorphic potentials 
       exist for all CMC immersions  defined on a non-compact 
       Riemann surface. We would expect that the analogous result holds 
       for a compact Riemann surface. But in this case the potential will 
       probably only be meromorphic.
 \end{enumerate}
}
 \end{remark}

 \subsection{}
 In section \ref{subsc:holopot}, the holomorphic potential was 
 derived from a CMC surface $f$.
 In this section, we give 
 the normalized potential for $f$ which is the 
 $\Lambda \mathfrak{sl}_2 \C$-valued meromorphic $1$-form 
 determined from the extended frame of $f$. Unlike the holomorphic 
 potential, the normalized potential is not holomorphic, 
 however, the Fourier expansion of it has only one coefficient.
 \begin{proposition}[Existence of a normalized potential]\label{prop:norpot}
  Let $\Phi :\mathbb D\to \Lambda {\rm SL}_2 \C_{\sigma, \tau_j}, 
  j = 3, 4$
  be the extended frame for some CMC surface $f$ in $\mathbb H^3$ with $H \neq 1$ 
  as defined in \eqref{eq:extendedH>1} or \eqref{eq:extendedH<1}.
  Then there exist $\Phi_+: \mathbb D \to \Lambda^+ {\rm SL}_2 \C_{\sigma}$
  and a $\Lambda \mathfrak{sl}_2 \C$-valued $1$-form $\xi$
  such that
  \begin{equation}\label{eq:normalizedpot}
  \xi = (\Phi \Phi_+)^{-1} \mathrm{d} (\Phi \Phi_+) = \lambda^{-1} \xi_{-1},
  \end{equation}
   where $\xi_{-1}$ is an off-diagonal $1$-form which is  meromorphic on
  $\mathbb D$. 
 \end{proposition}
 {\bf Proof}.
 Let us consider the Birkhoff decomposition of $\Phi$:
 \begin{equation}
 \Phi = \Phi_{-} \Phi_{+}, \;\; \Phi_{-} \in \Lambda_{*}^{-} {\rm SL}_2
  \C_{\sigma} \;\mbox{and}\;\Phi_{+} \in\Lambda^{+} {\rm SL}_2 \C_{\sigma}.
 \end{equation}
 It can be shown as in \cite{DPW} that this decomposition holds in $\mathbb D
 \setminus S$, where $S$ is a discrete subset of $\mathbb D$. 
 Moreover, this decomposition can be extended meromorphically across $S$.
 Differentiating $\Phi$ with respect to $\bar z$, one obtains 
 \begin{equation}\label{eq:Phiz}
 \Phi_{-}^{-1} \Phi_{-, \bar z} = \Phi_{+} V \Phi_{+}^{-1} - \Phi_{+, \bar
  z} \Phi_{+}^{-1}\;.
 \end{equation}
 Clearly, $\Phi^{-1}_{-} \Phi_{-, z} \in
 \Lambda^{-}\mathfrak{sl}_2 \mathbb C_{\sigma}$, and in view of 
 $\Phi_{-} \to {\rm Id}$ as $\lambda \to \infty$, the coefficient matrix 
 at $\lambda^{0}$ vanishes.
 Hence the left side of \eqref{eq:Phiz} contains only powers 
 $\lambda^{k}$ with $k <0$.
 Moreover, since $V$ only contains $\lambda^{0}$ and $\lambda^{1}$, 
 the right side of \eqref{eq:Phiz} contains only powers 
 $\lambda^{k}$, $k \geq 0$. 
 Therefore both sides of \eqref{eq:Phiz} vanishes.
 Thus $\Phi_{-, \bar  z} =0$, \textit{i.e.}, $\Phi_{-}$ is holomorphic 
 in $z$, where it is non-singular. Set $\xi = \Phi_{-}^{-1} \mathrm{d}
  \Phi_{-}$. Differentiating $\Phi$ with respect to $z$, one obtains
 similar to \eqref{eq:Phiz}
 \begin{equation}
 \xi = (\Phi_{+} U \Phi_{+}^{-1} - \Phi_{+, z} \Phi_{+}^{-1})\mathrm{d}z.
 \end{equation}
 Since $U$ only contains $\lambda^{-1}$ and $\lambda^{0}$, 
 the right hand side does not have $\lambda^{j}$ with $j \leq  -2$. 
 Thus $\xi$ has the form $\xi = \lambda^{-1} \xi_{-1}$,
 and the twisting condition implies that $\xi_{-1}$ is off-diagonal.
 $\Box$
\begin{remark}
{\rm 
  The meromorphic $1$-form $\xi$ in \eqref{eq:normalizedpot}
  will be called  the {\it normalized potential} of the immersion $f$.
}
\end{remark}
\section{Generalized Weierstrass type representation}\label{sc:DPW}
\subsection{}
 In Proposition \ref{prop:delbar} and \ref{prop:norpot} we have
 considered objects which are
 holomorphic in $\lambda \in \C^{\times}$. For the construction of
  CMC surfaces we need to obtain frames.
 Therefore we consider the double loop groups $\mathcal H = \Lambda
 {\rm SL}_2 \C_\sigma \times \Lambda {\rm SL}_2 \C_\sigma$.
 Then the subgroups $\mathcal H_+$ and $\mathcal
 H_-$ of $\mathcal H$ are defined as follows:  
\begin{align*}
\mathcal H_+ = \Lambda^+ {\rm SL}_2  \C_\sigma
 \times \Lambda^- {\rm SL}_2 \C_\sigma,\;\;
 \mathcal H_- = 
 \left\{ (g_1,\;\;g_2) \in \mathcal H \;
 \left| 
 \;\mbox{$g_1= g_2$}
 \right. 
 \right\}.
\end{align*}
 We quote Theorem 2.6 in \cite{DW:CMC-loop}.
\begin{theorem}[Generalized Iwasawa decomposition]\label{thm:generalIwasawa}
 The map $\mathcal H_- \times \mathcal H_+ \rightarrow \mathcal H_-
 \mathcal H_+$ is an analytic diffeomorphism. The image is open and
 dense in $\mathcal H$. Moreover,
 $$
 \mathcal H = \bigcup_{n=0}^{\infty} \mathcal H_- w_n \mathcal H_+\;\;,
 $$
 where $ w_n =   \left({\rm Id},\;\;\left(\begin{smallmatrix}\lambda^{n} & 
 0\\ 0 & \lambda^{-n}\end{smallmatrix}\right)\right) $ if $n=2k$ and 
 $  \left( {\rm Id},\;\; \left( \begin{smallmatrix} 0 & 
 \lambda^{n} \\ -\lambda^{-n} & 0 \end{smallmatrix}\right)\right)$ if
 $n=2k+1$.
\end{theorem}

 \subsection{}
 We recall that $\Lambda {\rm SL}_2 \C_{\sigma, \tau_j},\; j = 3, 4$, 
 defined in sections \ref{sc:CMCH>1} and \ref{sc:CMCH<1}, 
 are real forms of $\Lambda  {\rm SL}_2 \C_{\sigma}$, \cite{Kobayashi}.
 These real forms naturally induce Iwasawa decompositions of 
 $\Lambda {\rm SL}_2 \C_{\sigma}$:
 \begin{theorem}[Iwasawa decomposition for $\tau_3$, \cite{LoopGroup}]
 \label{thm:Iwasawa}
 The map 
 \begin{equation}
 \Lambda {\rm SL}_2 \C_{\sigma, \tau_3} \times \Lambda^+ {\rm
  SL}_2 \C_{\sigma}  \to  \Lambda {\rm SL}_2 \C_{\sigma}
 \end{equation}
 is an analytic diffeomorphism. 
 \end{theorem}
  The automorphism $\tau_3$ in the theorem above is induced from some 
  automorphism of ${\rm SL}_2 \mathbb C$. And this is used essentially in 
 the proof of this theorem. For the automorphism $\tau_4$ the situation 
 is completely different, since $\tau_4$ is not induced by some automorphism 
 of the (finite dimensional) Lie group ${\rm SL}_2 \mathbb C$.
\begin{theorem}[Iwasawa decomposition for $\tau_4$]\label{thm:Iwasawa2}
 Let $\omega_0 =  \left(\begin{smallmatrix} 0 & \lambda^{-1} \\
  - \lambda & 0\end{smallmatrix}\right)$. Then the map
 \begin{equation}
 \Lambda {\rm SL}_2 \C_{\sigma, \tau_4} \times
 \left\{{\rm Id}, \omega_0\right\} \times \Lambda^+ {\rm
  SL}_2 \C_{\sigma}  \to  \Lambda {\rm SL}_2 \C_{\sigma},\;\;
 (g, \delta, w_+) \mapsto g \delta w_+,
 \end{equation}
 is an analytic diffeomorphism onto an open dense subset
 of $\Lambda {\rm SL}_2 \C_\sigma$. The open dense subset 
 will be called  the {\rm Iwasawa core} with two open cells.
\end{theorem}
{\bf Proof}.
 Consider the map $g \in \Lambda {\rm SL}_2 \C_{\sigma} \mapsto
 (g, \tau_4 (g))  \in \Lambda {\rm SL}_2 \C_{\sigma} \times \Lambda {\rm
 SL}_2 \C_{\sigma}$, where $\tau_4$ is defined in \eqref{eq:tau4}.
 Applying the generalized Iwasawa decomposition of 
 Theorem \ref{thm:generalIwasawa}, we have
 \begin{equation*}
 (g, \tau_4(g)) = (\Phi, \Phi) ({\rm Id}, W) (V_+,  V_-),
 \end{equation*}
 Since $\Phi = g V_+^{-1}$ and $\Phi = \tau_4(g) V_{-}^{-1}W^{-1}$, 
 we obtain
 \begin{equation}\label{eq:splitggsharp}
 g^{-1}\tau_4(g) =V_+^{-1} W V_-\;\;\mbox{on $S^1$.}
 \end{equation}
 Since 
 $g^{-1}\tau_4(g)$ is in general not positive definite, the middle term 
 $W$ is not the identity element for the Birkhoff  decomposition 
 \eqref{eq:splitggsharp} in general.
 Consider the injective real analytic group homomorphism into 
 $$
 \Lambda {\rm SL}_2 \C_{\sigma, \tau_4}  \to \mathrm{Gr}(H), 
 \;\;g \mapsto A_{g} \cdot H_+,
 $$
 where $H = L^2 (S^1,\C)$ with polarization $H = H_+ \oplus H_-$, 
 \cite[Chap. 2]{DPW}. 
 Let $s$ be the non-trivial holomorphic section of the dual of the 
 determinant bundle over $\mathrm{Gr}(H)$. 
 We now set $\ell (g) =
 s (A_{g^{-1} \tau_4(g)} \cdot H_{+})$. It is known that 
 $\ell(g) \neq 0$ if and only if $g^{-1} \tau_4(g)$ is 
 in the big cell, \textit{i.e.}, $W = {\rm Id}$, 
 \cite[Corollary 2.5]{DPW}. 
 The set  given by 
 $\ell(g) =0$  in $\Lambda {\rm SL}_2 \C_\sigma$
 can not contain any open subset and is closed, 
 since $\ell (g)$ is  a real analytic function and $\ell ({\rm Id}) \neq 0$.
 Thus the set given by
 $\ell(g) \neq 0$ in $\Lambda {\rm SL}_2 \C_\sigma$
 is an open dense subset.
 Let's assume $\ell(g)\neq 0$ and $\ell({\rm Id}) \neq 0$, 
 \textit{i.e.}, $W ={\rm Id}$. 
 Then taking $\tau_4$ on both sides of \eqref{eq:splitggsharp}, we obtain
 $\tau_4(V_{+}^{-1} V_{-}) = (V_{+}^{-1} V_{-})^{-1}$
 and hence
 \begin{equation*}\label{eq:V+V-}
 \tau_4(V_{+}) = k^{-1} V_{-}\;\mbox{and}\;\tau_4(V_{-}) = k^{-1}
  V_{+},
 \end{equation*}
 where $k$ is a $\lambda$-independent diagonal matrix with entries
 $k_0, k_0^{-1} \in \mathbb R$.  Then from these symmetries, $\Phi$ 
 acquires the symmetry  $\tau_4(\Phi) = \Phi k$.
 If $k_0 >0$, we let $\tilde k$ be a real matrix such that 
 $\tilde k^2 =k$ holds. Then with 
 $\tilde \Phi = \Phi \tilde k$ and $\tilde V_+ = \tilde k^{-1} V_+$, 
 we obtain
 \begin{equation*}
 g = \tilde \Phi \tilde V_+, \;\;\tilde \Phi \in 
 \Lambda {\rm SL}_2 \C_{\sigma, \tau_4}, \;\; \tilde V_+ 
 \in \Lambda^+ {\rm  SL}_2 \C_{\sigma}.
 \end{equation*}
 If $k_0 < 0$, we set 
  $\omega_0 =  \left(\begin{smallmatrix} 0 & \lambda^{-1}  
 \\  - \lambda & 0\end{smallmatrix}\right)$, and choose
 a real matrix  $\check k$  such that 
 $\check k^2 = - k$ holds.
 Then with $\check \Phi = \Phi \check k \omega_0$ and 
 $\check V_+ = - \check k^{-1} V_+$, we have
 \begin{equation*}
 g = \check \Phi\omega_0 \check V_+, \;\;\check \Phi \in 
 \Lambda{\rm SL}_2 \C_{\sigma, \tau_4}, \;\; \check V_+ 
 \in \Lambda^+ {\rm  SL}_2 \C_{\sigma}.
 \end{equation*}
 To show that the sets $ \Lambda {\rm SL}_2 \C_{\sigma, \tau_4} \cdot 
 \Lambda^+ {\rm SL}_2 \C_\sigma $ and
 $\Lambda {\rm SL}_2 \C_{\sigma, \tau_4}
  \omega_0\Lambda^+ {\rm SL}_2 \C_\sigma$ are open in 
 $\Lambda{\rm SL}_2 \C_\sigma$, it suffices to show 
 that 
 \begin{equation}\label{eq:openset1}
 {\rm Lie}(\Lambda {\rm SL}_2 \C_\sigma)  = 
 {\rm Lie}(\Lambda {\rm SL}_2 \C_{\sigma, \tau_4})  \oplus 
 {\rm Lie}(\Lambda^+ {\rm SL}_2 \C_\sigma)
 \end{equation}
 and 
 \begin{equation}\label{eq:openset2}
 {\rm Lie}(\Lambda {\rm SL}_2 \C_\sigma)  = 
 \omega_0^{-1} {\rm Lie}(\Lambda {\rm SL}_2 \C_{\sigma, \tau_4})\omega_0 \oplus 
 {\rm Lie}(\Lambda^+ {\rm SL}_2 \C_\sigma),
 \end{equation}
 where ${\rm Lie} (G)$ denotes the Lie algebra of 
 $G$ \cite[Chap. II, Lemma 2.4]{Hel:Diff}. 
 But these two equations can be proven by a straightforward computation.
 This completes the proof. $\Box$
\begin{remark}
{\rm 
 The Iwasawa decomposition for $\tau_4$ in Theorem \ref{thm:Iwasawa2} can 
 be rephrased as follows: Let $\hat \tau$ be an extension of $\tau$ 
 in \eqref{6.1} to $G^{\C} = {\rm SL}_2 \C \times  {\rm SL}_2 \C$ 
 in section \ref{sc:4symmetric} and 
\begin{align*}
 \Lambda G^{\C}&= 
 \Lambda {\rm SL}_2 \C \times \Lambda {\rm SL}_2 \C, \\
 \Lambda G^{\C}_{\hat \tau} &= 
 \{(g, h) \in  \Lambda G^{\C} \;|\; 
 \hat{\tau}(g(\lambda),h(\lambda))=(g(-i \lambda),h(-i \lambda))\},  \\
 \Lambda G_{\hat \tau} &= 
 \left\{(g, h) \in  \Lambda G^{\C}_{\hat \tau} \;|\; 
 h(\lambda) = \overline{ g(1/\bar \lambda)} \right\}, \\
 \Lambda^+ G^{\C}_{\hat \tau}  &= 
 \{(g, h) \in  \Lambda G^{\C}_{\hat \tau} \;|\;  
 g \in  \Lambda^+ {\rm SL}_2 \C,\;\; h \in  \Lambda^- {\rm SL}_2 \C \}.
\end{align*}
 Since $\hat \tau^2$ is the involution $\hat \sigma = (\sigma, \sigma)$ 
 with $\sigma$ defined in \eqref{eq:sigma}, 
 we have 
 $$
 \Lambda G^{\C}_{\hat \tau} \subset \Lambda G^{\C}_{\hat \sigma} = 
 \Lambda {\rm SL}_2 \C_{\sigma} \times \Lambda {\rm SL}_2 \C_{\sigma}.
 $$
 Then the Iwasawa decomposition Theorem \ref{thm:Iwasawa2}
 can be rephrased as follows:
 The multiplication map 
 \begin{equation}
 \Lambda G_{\hat \tau} \times
 \left\{({\rm Id}, {\rm Id}), (\omega_0, \omega_0)\right\} \times 
 \Lambda^+ G^{\C}_{\hat \tau} 
 \to   \Lambda G^{\C}_{\hat \tau}, \;\; 
 (g, \delta, w_+) \mapsto g \delta w_+
 \end{equation}
 is an analytic diffeomorphism onto an open dense subset of 
 $\Lambda G^{\C}_{\hat \tau}$. 
}
\end{remark}

 \subsection{}
 In section \ref{sc:Potential}, we discussed potentials associated with CMC
 surfaces  in $\mathbb H^3$ with $H \neq 1$. In this section,
 we give  conversely a  construction of CMC surfaces from potentials, the
  {\it Generalized Weierstrass type representation}:

{\bf The case $H>1$:}

\noindent{\bf Step 1:} Take a holomorphic potential $\eta$ as defined 
 in \eqref{eq:holopot}. 

\noindent{\bf Step 2:} Solve the ordinary differential equation 
 $\mathrm{d} C = C \eta$ with $C(z_*, \lambda) = {\rm Id}$, where $z_*$ 
 is some base point.
                      
\noindent{\bf Step 3:} 
 Perform the Iwasawa decomposition for $C = \Phi V_{+}$ by 
 Theorem \ref{thm:Iwasawa}, where 
 $\Phi \in \Lambda {\rm SL}_2 \C_{\sigma, \tau_3}$ 
 and $V_{+} \in \Lambda^+ {\rm SL}_2 \C_{\sigma}$.

\begin{theorem}\label{thm:ExtFrame}
 Up to a diagonal gauge $D \in {\mathrm U}_1$
 and a change of coordinates, 
 the matrix  $\Phi$ obtained in
 Step 3 is  the extended frame of some CMC
 surface with $H>1$.
\end{theorem}

\noindent{\bf Step 4:} Inserting $\Phi$ obtained in Step 3
 into the Sym formula in  \eqref{eq:SymH>1} we obtain a CMC surface 
 in $\mathbb H^3$ with $H = \coth q >1$.

 {\bf Proof of Theorem \ref{thm:ExtFrame}}.
  Let us compute the Maurer-Cartan form of $\Phi = C V_{+}^{-1}$, where
 $\Phi \in \Lambda {\rm SL}_2 \C_{\sigma, \tau_3}$. 
 A direct computation shows that 
 \begin{equation}
 \alpha = \Phi^{-1} \mathrm{d} \Phi = V_{+} C^{-1} \mathrm{d} 
 C V_{+}^{-1} - \mathrm{d} V_{+} V_{+}^{-1}.
 \end{equation}
 Thus $\alpha$ has the form $\alpha = \sum_{j = -1} ^{\infty}
 \lambda^{j} \alpha_{j}$. Since $\Phi$ is an element in
 $\Lambda {\rm SL}_2 \C_{\sigma, \tau_3}$, $\alpha$ satisfies 
 that $\tau_3 (\alpha) = \alpha$, where $\tau_3$ 
 is defined in \eqref{eq:invc3}. 
 Therefore 
 \begin{equation}
 \alpha = \lambda^{-1} \alpha_{-1} + \alpha_0 + \lambda \alpha_1\;,
 \end{equation}
 where $\overline{\alpha_0} = - \alpha_0$ and $\alpha_{-1} =
 {\alpha^*}_1$.
 It is easy to check that $\alpha$ has the form \eqref{LaxH>1}
 up to a diagonal gauge 
 $D =\diag (e^{i \theta}, e^{-i \theta})$ and a change of coordinates
 \cite[Section A.8]{DH:Mero}. 
 Since $D \in {\mathrm U}_{1}$, $C = \Phi D \cdot D^{-1} V_+$ is
 also an Iwasawa decomposition. Thus, by Proposition
 \ref{prop:fundamentalH>1}, $\Phi D$ is
 the extended frame of some CMC surface  with $H>1$. 
 This completes the proof. $\Box$

{\bf The case $0 \leq H<1$:}

\noindent{\bf Step 1:} Take a holomorphic potential $\eta$ as defined 
 in \eqref{eq:holopot}. 

\noindent{\bf Step 2:} Solve the ordinary differential equation 
 $\mathrm{d} C = C \eta$ with $C(z_*, \lambda) = {\rm Id}$, where $z_*$ 
 is some base point.
                      
\noindent{\bf Step 3:} 
 Perform the Iwasawa decomposition for $C$, \textit{i.e.}, 
 $C = \Phi V_{+}$ or 
 $C = \Phi \omega_0 V_{+}$
 for all $z$ in $\mathbb D$ such that $C(z, \lambda)$
 is in $\mathbb D_1 \subset \mathbb D$ 
 or $\mathbb D_2 \subset \mathbb D$ respectively, where 
 $\Phi \in \Lambda {\rm SL}_2 \C_{\sigma, \tau_4}$, 
 $\omega_0 = \left(\begin{smallmatrix} 0 & \lambda^{-1} 
 \\ -\lambda & 0\end{smallmatrix}\right)$, 
 $V_{+} \in \Lambda^+ {\rm SL}_2 \C_{\sigma}$, and
 $\mathbb D_1$ and $\mathbb D_2$ are open subsets in $\mathbb D$ 
 such that the real valued functions $k_0$ defined in 
 Theorem \ref{thm:Iwasawa2} are positive and negative on 
 $z \in \mathbb D_1$ and $z \in \mathbb D_2$ respectively.

\begin{theorem}\label{thm:ExtFrame2}
 Up to a diagonal gauge $D \in {\mathrm U}_1$ and a change of coordinates, 
 the matrices $\Phi$ and $\Phi \omega_0$ obtained in
 Step 3 are the extended frames of some CMC
 surface with $0 \leq H <1$. 
\end{theorem}

\noindent{\bf Step 4:} Inserting $\Phi$ or $\Phi\omega_0$ 
 obtained in Step 3 into the Sym formula in \eqref{eq:SymH<1}, 
 we obtain a CMC surface in $\mathbb H^3$ with $0\leq H <1$.

{\bf Proof of Theorem \ref{thm:ExtFrame2}}.
 Let us compute the Maurer-Cartan form of $\Phi$ if 
  $C = \Phi V_+ $ and for $\Phi\omega_0$ if $C = \Phi \omega_0 V_+$ with 
 $\Phi  \in \Lambda {\rm SL}_2 \C_{\sigma, \tau_4}$ and $V_+ 
 \in \Lambda {\rm SL}_2 \C_{\sigma, \tau_4}$. 
 Since $\Phi$ is an element in $\Lambda {\rm SL}_2 \C_{\sigma, \tau_4}$, 
 $\alpha$ satisfies that $\tau_4(\alpha) = \alpha$, where 
 $\tau_4$ is defined in \eqref{eq:invc4}.
 Similar to the proof of Theorem \ref{thm:ExtFrame},  
 the Maurer-Cartan form of $\Phi$ and $\Phi \omega_0$ respectively
 has the form
 \begin{equation}
 \alpha = \lambda^{-1} \alpha_{-1} + \alpha_0 + \lambda \alpha_1\;,
 \end{equation}
 where $\overline{\alpha_0} = - \alpha_0$ and $\alpha_{-1} =
 i {\rm Ad}(\mathscr{R})(\alpha_1)^*$.
 It is easy to check that 
 up to a diagonal gauge 
 $D =\diag (e^{i \theta}, e^{-i \theta})$ 
 and a change of coordinates \cite[Section A.8]{DH:Mero},
 $\alpha$ has the form \eqref{LaxH<1}.
 Since $D \in {\mathrm U}_{1}$ and $D \omega_0 D = \omega_0$,
 $C = \Phi D \cdot D^{-1} V_+$ and 
 $C = \Phi \omega_0 D \cdot D^{-1} V_+ = 
 \Phi D^{-1} \omega_0 \cdot D^{-1} V_+$ are 
 also Iwasawa decompositions. Thus, by Proposition
 \ref{prop:fundamentalH<1}, $\Phi D$ and $\Phi D^{-1} \omega_0 $ are
 the extended frame of some CMC surface  with $0\leq  H<1$. 
$\Box$

\subsection{}
 Let $\eta$ be a potential, defined on $\mathbb D$, for a CMC immersion
 into $\mathbb H^3$ with mean curvature $0\leq H <1$. Let $C$ be defined 
 by $\mathrm{d} C = C\eta$ and $C(z_*, \lambda) ={\rm Id}$, 
 where $z_* \in \mathbb D$
 is some base point. In the proof of Theorem \ref{thm:Iwasawa2}, 
 we have seen that 
 there exists a real analytic function $k_0 : \mathbb D \to \mathbb R$
 such that the ``singular set'' $S_0 =\{ z \in \mathbb D \;\vert\; 
 k_0 (z, \bar z) =0\}$ divides $\mathbb D$ into two open subsets 
 $\mathbb D_1$ and $\mathbb D_2$, where $C(z)$ is in the open Iwasawa cell 
 containing ${\rm Id}$ for $z \in \mathbb D_1$ and in the other open Iwasawa 
 cell for $z \in \mathbb D_2$. Let $C  = \Phi V_+$ be an Iwasawa decomposition 
 on the open Iwasawa cell containing ${\rm Id}$ for $z \in \mathbb D_1$. 
 Examples like the one given in 
 section \ref{subsec:sphere} show that 
 the frame $\Phi$ associated with $C$ will be generically singular along $S_0$.
 Hence it does not seem to make sense to extend the immersion $f$ associated 
 with $\eta$ in $\mathbb D_1$ across $S_0$ to $\mathbb D_2$. 
 However, in view of 
 \cite[Theorem 3.2]{DK:cyl},  one can extend $\Phi l $ meromorphically to 
 $\mathbb D \times \bar{\mathbb D}$,
 if $l$ is an appropriately chosen diagonal matrix independent of $\lambda$.
 Actually, with $\tilde k$ as in the proof of Theorem \ref{thm:Iwasawa2} 
 the matrices, 
 $\hat \Phi = \Phi \tilde k^{-1}$ and $\tilde k^2$ have meromorphic 
 extension to $(z, w) \in \mathbb D \times \bar{\mathbb D}$. As a consequence, 
 consider the CMC immersion of mean curvature $H = \tanh q$ given by the 
 Sym formula $f = \Phi D_0 \Phi^* $ for $z\in \mathbb D_1$, where 
 $D_0 = \diag (e^{-q/2}, e^{q/2})$.
 Then we rephrase the Sym formula as  
\begin{eqnarray*}
 f &=& \Phi D_0 \Phi^*=(\Phi \tilde k^{-1}) \tilde k^{2}D_0 
       (\Phi \tilde k^{-1})^*.
\end{eqnarray*}
 Then by the argument above, $\Phi \tilde k^{-1}$ and  $\tilde k^2$ 
 can be extended meromorphically  to 
 $(z, w) \in  \mathbb D \times \bar{\mathbb D}$.
 Let's put $\hat \Phi = \Phi \tilde k^{-1}$  and $k=\tilde k^2$. 
 Then $\hat f = \hat \Phi k D_0 \hat \Phi^*$
 has a meromorphic extension to $(z, w) \in \mathbb D \times \bar {\mathbb D}$.
 We note that, on $w = \bar z$, 
 $\tau_4(\hat \Phi) = \hat \Phi k$ and $k$ has real diagonal entries,
 where $\tau_4$ is defined in \eqref{eq:tau4}.
  We can now consider the 
 analytic continuation from the one open Iwasawa cell for
 $z \in \mathbb D_1$, to  the other open Iwasawa cell for $z \in \mathbb D_2$, 
 \textit{i.e.}, 
 $C = \tilde \Phi \omega_0 \tilde V_+$ is the Iwasawa 
 decomposition for $z \in \mathbb D_2$.
 Then, it is easy to check that
 the entries of $k$ are negative on $z \in \mathbb D_2$. 
 Thus we can rephrase $k$ as $k = - \check k^2$, 
 where $\check k$ is the diagonal matrix with positive entries 
 which are independent of $\lambda$ defined as in Theorem \ref{thm:Iwasawa2}. 
 Moreover
 $C =(\hat \Phi \check k \omega_0) \omega_0 (- \check k^{-1} \hat V_+)$
 is the Iwasawa decomposition for $z \in \mathbb D_2$,  
 where $\hat V_+$ is the meromorphic extension
  of $\tilde k V_+$. Using an obvious abbreviation
 $\tilde \Phi =\hat \Phi \check k \omega_0$, 
 we can rephrase again the Sym formula as  
\begin{eqnarray*}
 \hat f &=& \tilde \Phi  \omega_0 
 \check k^{-2}  k D_0 \omega_0^*  \tilde \Phi^*= -(\tilde \Phi \omega_0) D_0 
     (\tilde \Phi\omega_0)^*.
\end{eqnarray*}
 Therefore, it is natural to use for $z$ in $\mathbb D_2$ 
 the negative of the Sym formula \eqref{eq:SymH<1} and to 
 use this formula for $\tilde\Phi \omega_0$.
 Thus in the second open Iwasawa cell actually $\tilde\Phi \omega_0$ 
 is the ``frame'' to use.

\subsection{}
 We have seen in section \ref{sc:CMCH<1} that coefficient matrices of
 the form \eqref{LaxH<1} define CMC surfaces in $\mathbb H^3$ with mean
 curvature $0 \leq H < 1$.

 In section \ref{sc:DPW} we have shown that all such matrices can be
 obtained via a generalized Weierstrass type representation procedure from
 ``potentials''. Thus we have a loop group procedure to construct all CMC
 surfaces in $\mathbb H^3$ with mean curvature $0 \leq H < 1$.

 On the other hand, a simple minded approach would follow \cite{DPW} more
 closely, since we know (Ishihara, Proposition \ref{prop:Ishihara}) that
 the surfaces under consideration are characterized by the fact that
 their ``Gauss map'' \eqref{eq:mapsinUH3} is a harmonic map into
 $\mathrm{U}\mathbb{H}^3 = \mathrm{SL}_2 \C/\mathrm{U}_1$. Since
 $\mathrm{U}\mathbb{H}^3$ is not a symmetric space, only a $4$-symmetric
 space, \cite{DPW} or \cite{BD1} cannot be applied directly. 

 In section
 \ref{sc:4symUnittangent} we have shown that
 one can indeed introduce a loop parameter into the Maurer-Cartan form
 $\alpha$ as in \cite[section 3.2]{BP},
 since the Gauss map is Legendre:
$$
 \alpha = \alpha_{\mathfrak h} + \nu^{-1} \alpha^{\prime}_{\mathfrak p}
 + \nu \alpha^{\prime \prime}_{\mathfrak p}, \;\; \nu \in \C^{\times},
$$
 where $\alpha_{\mathfrak h}$, $\alpha_{\mathfrak p}$
 are defined in \eqref{eq:alphaLegendre}. 
 Integration of $\alpha$ yields an extended frame $\Phi = \Phi(z,\bar z,\nu)$.

 There is no loop group method known for general harmonic maps into
 $k$-symmetric spaces if $k>2$. However, since we were able to introduce
 a loop parameter one could decompose $\Phi$ anyway, say \'a la Birkhoff,
 see Theorem \ref{thm:Birkhoff}:
$$
\Phi = \Phi_{1-} \Phi_{2+} = \Phi_{1+} \Phi_{2-}, \;\;\Phi_{1\pm} \in
\Lambda^{\pm}_{*} \mathrm{SL}_2 \C_{\sigma}, \;\Phi_{2\pm} \in
\Lambda^{\pm} \mathrm{SL}_2 \C_{\sigma}.
$$
 This would yield a potential $\xi = (\xi_-, \xi_+) = (\Phi_{1-}^{-1}
 \mathrm{d} \Phi_{1-}, \Phi_{1+}^{-1}\mathrm{d} \Phi_{1+})$.

 A priori, the only information we have about $\xi_\pm$ is that they are 
 elements in $\Lambda^{\pm} \mathfrak{sl}_2 \C_{\sigma}$ with
 simple loop behavior. 
 
 If one would want to construct, conversely, extended frames from
 potentials of this type, one would run irreparably into the problem of
 violating the admissibility condition \eqref{admissible}.
 The observation that $\alpha$ is fixed under $\tau_4$, however, permits to
 fix this problem. Since $\tau_4 (\Phi_{1-}) = \Phi_{1+}$, we have 
$$
\tau_4(\xi_-)=\xi_+\;.
$$
 Similar to \cite{Kobayashi}, one can now show that for such potentials
 one can obtain an extended frame which is fixed under $\tau_4$.
 
 It is not difficult to see that following \cite{Kobayashi}
 corresponds exactly to the approach (using an Iwasawa decomposition)
 explained in section \ref{sc:CMCH<1} and section \ref{sc:DPW}.

\subsection{} \label{subsec:finitetype}
 Among the classes of surfaces that can be constructed by integrable 
 system methods the surfaces of ``finite type'' have received special
 attention. This is due to the fact that these surfaces can be constructed
 fairly explicitly. Moreover, the algebro-geometric methods used for these 
 surfaces are classical and beautiful. These methods have been highly 
 successful for the construction of tori. On the other hand, since until
 now only surfaces without umbilical points can be constructed by 
 algebro-geometric methods, it is necessary to develop more general
 methods which also work when umbilical points are present.

 In the context of this paper the property of being of finite type 
 can be expressed on all levels of our construction scheme. 
 We sketch  these constructions, but will not present any proofs.

 The starting point for the construction of a CMC surface in $\mathbb H^3$
 with $0 \leq H <1$ in this paper are ``potentials''. In the spirit
 of \cite{BP} and \cite{DoSter}, we call a potential $\eta$ to be of
 ``finite type'' if there exists a Laurent polynomial 
 ${\xi}^o(\lambda)$
 satisfying $\tau_4 \xi^o = \xi^o$ such that 
 $\eta = \lambda^{d-1} \xi^o(\lambda)
 \mathrm{d}z$, where $d$ is odd and the maximal degree of $\lambda$ 
 occurring in $\xi^o$. In our construction scheme one obtains from 
 $\eta$ a CMC surface in $\mathbb H^3$ with $0 \leq H<1$. 
 Surfaces constructed from some  $\xi^o$ as above are usually 
 called to be of {\it Symes finite type} in the sense of \cite{DoSter}. 

 Our construction scheme produces from $\eta$ an ODE solution 
 $\mathrm{d} C = C \eta$
 satisfying $C(z_*, \lambda) = {\rm Id}$ and from $C$ via Iwasawa decomposition 
 $C = \Phi V_+$, $\Phi(z_*, \bar z_*, \lambda) ={\rm Id}$, 
 the frame $\Phi$ of the CMC  surface constructed from $\eta$. 
 As in \cite{BP}  (see also \cite{DoSter}) one can show 
 the following theorem:
\begin{theorem}
 The Maurer-Cartan form $\alpha = \Phi^{-1} \mathrm{d} \Phi$ can be represented 
 in the form 
\begin{equation}
 \alpha = \alpha_{\xi} = 
 (\lambda^{-1}\xi_{-d}+\frac{1}{2} \xi_{-d+1})\mathrm{d}z +
 \tau_4 (\lambda^{-1}\xi_{-d}+\frac{1}{2} \xi_{-d+1})\mathrm{d}\bar z, 
\end{equation}
 where $\xi = \sum_{n \in \mathbb Z} \xi_n$ satisfies $\tau_4 \xi = \xi$ 
 and the Lax equation 
 \begin{equation}\label{eq:phixi}
 d \xi =[\xi, \alpha_{\xi}], \;\;\xi(z_*) = \xi^o.
 \end{equation}
 Moreover the $\xi$ is polynomial in 
 $\lambda$ and satisfies
\begin{equation}
 \xi = \Phi^{-1} \xi^o \Phi.
\end{equation}
 Such a  $\xi$ is called the 
 {\rm polynomial Killing field}.
\end{theorem}
\begin{remark}
{\rm 
 In view of $\Phi = C V_+^{-1}$ one obtains that the degree of $\xi$ 
 is equal to the degree of $\xi^o$.
 Again following \cite{BP} (see also \cite{DoSter}) one calls a frame $\Phi$ 
 to be of {\it finite type} if there exists some $\xi^o$ such that the 
 solution to \eqref{eq:phixi} satisfies $\Phi^{-1} \mathrm{d} \Phi = \alpha_{\xi}$. 
 Moreover, a CMC surface in $\mathbb H^3$ with $0 \leq H <1$ is 
 called of {\it finite type} if its frame is of finite type. With 
 this notation the theorem above states that finite type surfaces 
 are equivalent to Symes finite type surfaces. 
}
\end{remark}
 There are basically two ways to treat finite type surfaces 
 algebro-geometrically. One way is to consider the spectral curve 
 $Y =\{ (\lambda, \mu) \in \mathbb C^{\times} \times \mathbb C \;\vert \;
 \mathrm{det}(\mu \cdot {\rm Id}- \xi(z, \lambda)) =0\}$ and the eigenline 
 bundle $\mathscr L$ of $\xi$ over $Y$. 
 This approach was carried out by I.~McIntosh for 
 finite type harmonic maps into $\mathbb C P^n$ \cite{Mc:CPn} 
 and it seems that his method should also apply to our case. 
 The second approach is due to A.~Bobenko, it produces algebro-geometric
 solutions to the complex sine-Gordon equation from line bundles over 
 hyperelliptic curves via the corresponding theta functions. For the 
 surfaces considered in this paper this was carried out in \cite{BaBo}. 
 
 In the case of CMC surfaces in $\mathbb R^3$, the relation between 
 the loop group method as used in this paper and the approach of Bobenko
 was established in \cite{DoH}. It would be interesting to see whether also for 
 CMC surfaces of finite type in $\mathbb H^3$ a relation with Bobenko's 
 approach can be realized in a similar fashion.

\section{\bf Examples of CMC immersions with $0 \leq H <1$ in  $\mathbb H^3$}
 \label{sc:Ex}
 In \cite{KKRS}, examples  of {\sc CMC}-immersions with $H>1$ have
 been constructed from  holomorphic potentials. In this section, we
 present some examples of {\sc
 CMC}-immersions with $0 \leq H <1$ starting from holomorphic potentials.
 From now on we assume always $q \in \mathbb R_{\geq 0}$ and use the symbol
 $\lambda$ instead of $\nu$. We recall that the resulting immersion by the
 Sym formula $f_{\lambda = e^{-q/2}}$ in \eqref{eq:SymH<1} has
 the mean curvature $0 \leq H = \tanh q<1$. 
 
 All the examples given below have the feature that there exist some curves 
 $\mathcal C$ along which are immersions which tend to infinity. 
 Each of the connected components of the immersion from $M \setminus \mathcal C$
 is well defined  and takes values in the forward $\mathbb H^3$ or 
 the backward  $\mathbb H^3$ (in the light cone picture of hyperbolic space). 
 We hope to study this behavior in detail in a separate publication.

\begin{figure}[h]
\begin{center}
 \includegraphics[width=0.4\textwidth]{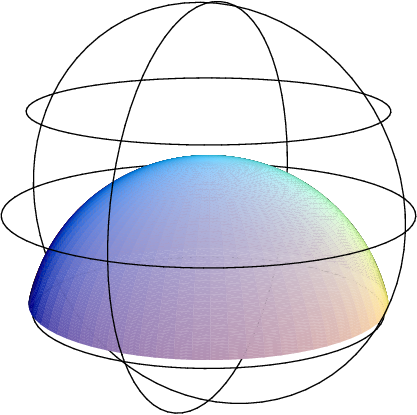}\hspace{2cm}
 \includegraphics[width=0.4\textwidth]{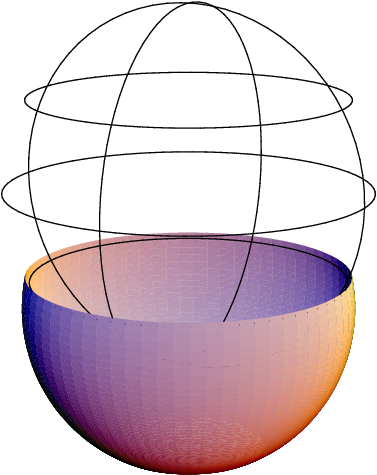}
 \caption{
 An equidistance surface with $H = \tanh (1)$, $|z|<1$ (left) and 
 an equidistance surface with $H = \tanh (1)$, $|z|>1$  (right).
 Surfaces are shown in the Poincar{\'e} ball model and the outside 
 of the Poincar{\'e} ball model respectively.
 }\label{fig:one}
\end{center}
\end{figure}

\subsection{Umbilical surfaces (Figure \ref{fig:one})} \label{subsec:sphere}
 Assume an CMC immersion $f$ with $0 \leq H<1$ to be totally umbilical.
 It is easy to see that the Hopf differential $Q \mathrm{d}z^2$ of an
 umbilic CMC immersion with $0 \leq H < 1$ vanishes, \textit{i.e.}, $Q =0$. 
 Then the corresponding normalized potential $\xi$ as in
 \eqref{eq:normalizedpot} has the form 
 $$
 \xi = \lambda^{-1} \begin{pmatrix}0 & h(z) \\ 0 & 0 \end{pmatrix} \mathrm{d}z,
 $$
 where $h(z)$ is a non-vanishing meromorphic
 function. Using the coordinate change
  $w = \int h \mathrm{d}z$, $\xi$ has the form defined in
 \begin{equation}\label{eq:umbiliceta}
 \eta =
 \begin{pmatrix}
 0 & \lambda^{-1} \\
 0 & 0
 \end{pmatrix} \mathrm{d}z. 
 \end{equation}
 Let us compute totally umbilic surfaces explicitly from holomorphic
 potentials $\eta$  defined in \eqref{eq:umbiliceta}. 
 A solution $C$ to the holomorphic differential equation 
 $\mathrm{d} C = C \eta$ with the above $\eta$ 
 and the Iwasawa decomposition for $C = \Phi V_+$
 can be computed explicitly as follows:
 \begin{align*}
 C =
 \begin{pmatrix} 1 & \lambda^{-1} z \\ 0 &1 \end{pmatrix}
   = \Phi V_{+} = 
  \begin{pmatrix}
 \frac{1}{\sqrt{1-|z|^2}}  & \frac{\lambda^{-1} z}{\sqrt{1-|z|^2}} \\
  \frac{\lambda \bar z }{\sqrt{1-|z|^2}} & \frac{1}{\sqrt{1-|z|^2}}
 \end{pmatrix}
 \begin{pmatrix}
  \frac{1}{\sqrt{1-|z|^2}}   & 0 \\
 \frac{-\lambda \bar z}{\sqrt{1-|z|^2}} & \sqrt{1-|z|^2}
 \end{pmatrix}\;. 
\end{align*}
 The Sym formula \eqref{eq:SymH<1} for the above $\Phi$ with $\lambda =
 e^{-q/2}$ gives  
\begin{equation*}
f_{\lambda} = \frac{1}{1-|z|^2}\begin{pmatrix}
e^{-q/2} + e^{3q/2 }|z|^2 & (e^{-q} + e^{q}) z\\ (e^{-q}+e^{q}) \bar z&
e^{q/2}+ e^{-3q/2 } |z|^2\end{pmatrix}
\;.
\end{equation*}
 Thus the resulting immersion is a totally umbilic surface, which is 
 called an {\it equidistance surface} for $0 < H <1$ and a {\it
 totally geodesic surface} for 
 $H=0$  \cite[Theorem 29, p.~77]{Spivak}.

\begin{figure}[h]
\begin{center}
 \includegraphics[width=0.4\textwidth]{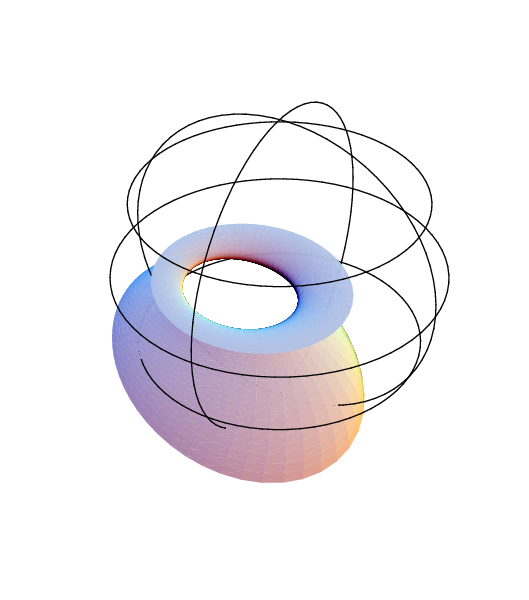}\hspace{2cm}
 \includegraphics[width=0.4\textwidth]{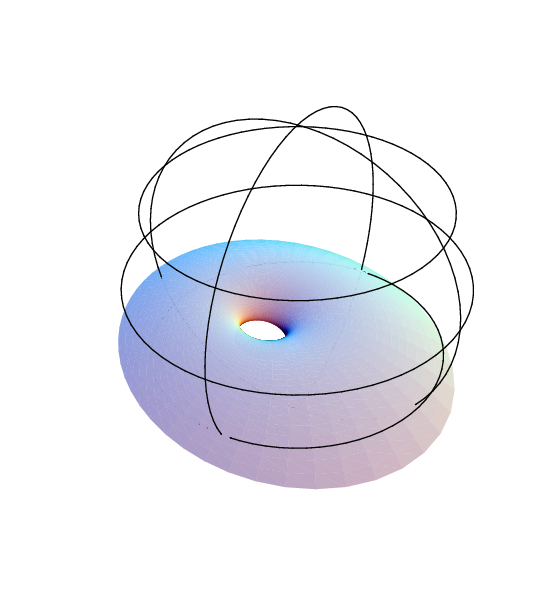}
 \caption{
 A portion of a surface of revolution with $H = \tanh (0.3)$ 
 and $a =0.3$ (left) and 
 a portion of a minimal surface of revolution with $a =0.2$ (right).
 Surfaces are shown in the Poincar{\'e} ball model and the outside 
 of the Poincar{\'e} ball model.
 }\label{fig:two}
\end{center}
\end{figure}

 \subsection{Surfaces of revolution (Figure \ref{fig:two})} \label{subsec:Equivariant}
 We set
 \begin{equation}\label{eq:Delpotential}
 \eta = A \mathrm{d}z =
  \begin{pmatrix}
   0 & \lambda^{-1} a + \lambda b \\
  \lambda^{-1} b - \lambda a& 0
 \end{pmatrix} \mathrm{d}z\;, 
\end{equation}
 where $a, b \in \mathbb R$ and $b^2 -a ^2 + a b(e^{q} - e^{-q}) =1/4$.
 A solution to the holomorphic differential equation $\mathrm{d} C = C \eta$ with 
 the above potential $\eta$ is 
\begin{equation}\label{eq:Delsol}
 C = \exp \left(A z \right)\;.
\end{equation}
 Let $\gamma : z \mapsto z + p$ be an automorphism of $\C$ with $p  \in \C^{\times}$. 
 Noting that $\tau_4 A = - A$, it is easy to check that 
 $\gamma^* C = C (p, \lambda) C$ and $C(p, \lambda) \in
  \Lambda {\rm SL}_2 \C_{\sigma, \tau_4}$ for all 
 $p \in i \mathbb R$ and $|\lambda| = e^{-q/2}$. 
  
 Let $p$ be purely imaginary and let $C = \Phi  V_+$ denote the
 Iwasawa decomposition of $C$. 
 Then the Iwasawa decomposition for $\gamma^* C =  \hat \Phi \hat V_+$
 can be  computed as $\hat \Phi= C (p, \lambda) \Phi$ and  $\hat V_+ = V_+$. 
 Thus $C_{q}(p) = C(p, \lambda=e^{-q/2})$ acts on the resulting
 immersion $f_{\lambda}$  at $\lambda = e^{-q/2}$, which is denoted by
 $f_q$, as a $1$-parameter group of isometries with group
 parameter $p \in i  \mathbb R$:
 \begin{equation}\label{eq:oneparameter}
  \gamma^* f_{q} = C_{q}(p) f_{q} C_{q}(p)^*.
 \end{equation}
 Moreover, $C_q(p)$ has the eigenvalues
 $$
 \exp (\pm p \sqrt{b^2 -a^2 + a b(\lambda^{-2} - \lambda^2)})|_{\lambda
 =  e^{-q/2}} = \exp (\pm  p/2).
 $$ 
 Since  $C_q(p) = \pm {\rm  Id}$ for $p = 2k \pi i$ with 
 $k \in \mathbb Z$, the resulting immersion closes up for $p= 2k \pi i$, 
 \textit{i.e.}, $\gamma^* f_{q} = f_{q}$. 
 The resulting immersion defines a surface
 of revolution with respect to the $1$-parameter group
 \eqref{eq:oneparameter}. As its profile curve one can choose:
  \begin{equation*}
  g (x) = \Phi \begin{pmatrix} e^{-q/2} & 0\\ 0& e^{q/2}\end{pmatrix}
          {\Phi}^* |_{z = x, \lambda = e^{-q/2}},
  \end{equation*}
  where $x \in \mathbb R$. For the axis of the surface of revolution we obtain
 \begin{equation*}
  \ell =\{ X \in {\rm Her}_2 \C\;|\; C_q(p)X -
   XC_q(p)^{*-1} = 0, \; p \in i\mathbb R\} \cap \mathbb H^3.
 \end{equation*}
 Even though the axis $\ell$ could be an empty set, it is natural to call
 the resulting immersion a surface of revolution \cite{Mori}.

\subsection{Radially symmetric surfaces}
\label{subsec:radiall}
 We set 
\begin{equation*}
 \eta = \begin{pmatrix} 0 &  \lambda^{-1}\\
        z^k \lambda^{-1} & 0\end{pmatrix} \mathrm{d}z\;,
\end{equation*}
 where $k \in \mathbb N$. Set $\epsilon : (z, \lambda) \mapsto (a z,
 a^{\frac{k+2}{2}}\lambda)$. Then $\epsilon^* \eta  = T \eta T^{-1}$, 
 where $T= \diag (a^{\frac{-k}{4}}, a^{\frac{k}{4}})$. Moreover,
 we assume that $a \in S^1$. Let $C$ denote the solution to $\mathrm{d} C = C
 \eta$, $C(0, \lambda) ={\rm Id}$. Let $C = \Phi V_+$ be the Iwasawa
 decomposition of $C$ and assume that $\Phi$ is in the first open Iwasawa
 cell. Then the Iwasawa decomposition for $\epsilon^* C$ is 
 $(\epsilon^* \Phi) (\epsilon^* V_+)$, thus $\epsilon^* \Phi$ is 
 again in the first open Iwasawa cell. Moreover $\epsilon^* \Phi$ and 
 $\epsilon^* V_+$ can be computed
 as $\epsilon^* \Phi = T \Phi T^{-1}$ and
 $\epsilon^* V_+ = T V_+ T^{-1}$.  Since $V_+(\lambda =0)$ is diagonal,
 the last relation shows that $V_{+}(\lambda =0)$ only depends on $|z|$.
 In view of $C = \Phi V_+$, comparison of the potential considered in
 this section with \eqref{LaxH<1} shows that the metric only depends on
 $|z|$. The resulting surface is a radially symmetric surface. The minimal 
 surface of this type has been investigated in \cite{Novokshenov}.

\subsection{Cylinders}\label{subsec:cylinders}
 We set 
\begin{equation*}
 \eta = S \mathrm{d}z = \begin{pmatrix} 0 &  (a \lambda^{-1} + b \lambda) h(z)\\
	       (- b \lambda^{-1} + a \lambda) \overline{h(\bar z)}&
	       0\end{pmatrix} \mathrm{d}z\;, 
\end{equation*}
 where $h(z)$ is a periodic holomorphic function on $\C$ with 
 period $p \in \mathbb R^{\times}$ and $a, b \in \mathbb R^{\times}$. Moreover we
 assume $ b = a e^{-q}$. It is easy to check that $S$ has the following two
 properties: $\tau_4 S = S$ and
 $S ( z, \lambda = e^{-q/2})$ is an upper triangular matrix.
 
 Let $C(z, \lambda)$ denote the solution to $\mathrm{d} C = C \eta$, $C(z=0,
 \lambda) ={\rm Id}$. By the Picard-Lindel\"{o}f iteration, 
 $C$ can be computed as
 \begin{equation*}
 C = {\rm Id} + \int_0^z C \eta \mathrm{d}t+ 
     \int_0^z \left(\int_0^t \eta \mathrm{d}t_1
     \right) \eta  \mathrm{d}t+ \cdots\;.
 \end{equation*}
 Set $\gamma : z \mapsto z +p$ with $p \in \mathbb R^*$. 
 Then $\eta$ is invariant under $\gamma$, 
 and we obtain $\gamma^* C = C(p) C$.
 From the properties of $S(z, \lambda)$ stated above, 
 we obtain $C(p) \in \Lambda {\rm SL}_2 \C_{\sigma, \tau_4}$
 and 
 \begin{equation*}
 C(p)|_{\lambda = e^{-q/2}} = \begin{pmatrix} 1 & \int_0^p h(t)
			       \mathrm{d}t \\ 0 & 1
		       \end{pmatrix}\;.
 \end{equation*}
 Let $C = \Phi V_+$ be the Iwasawa decomposition of $C$. Then 
 $\gamma^* C$ has the Iwasawa decomposition $\gamma^* C = \hat \Phi \hat
 V_+$, where $\hat \Phi = C(p) F$ and $\hat V_+ = V_+$, since 
 $C(p) \in \Lambda {\rm SL}_2 \C_{\sigma, \tau_4}$. Therefore the 
 immersion $f_q$ obtained by the Sym formula
 closes up, \textit{i.e.}, $\gamma^* f_q = f_q$, if and  only if
 $\int_0^p h(t)  \mathrm{d}t =0$. It is not difficult to find 
 periodic holomorphic functions $h(z)$ with period $p$ satisfying this
 condition.

\subsection{Totally symmetric surfaces}
 \label{subsec:trinoids}
 We set 
\begin{equation}\label{eq:trinoidspot}
 \eta = 
 - \begin{pmatrix}
    0 & \lambda^{-1} \\ \lambda X^2 &0  
   \end{pmatrix} 
   \frac{\mathrm{d}z}{z}  + 
   \begin{pmatrix}
    0 & \lambda^{-1} \\  \lambda X^2 &0  
   \end{pmatrix} 
   \frac{\mathrm{d}z}{z-1}  +  
   \begin{pmatrix}  
    0 & 0 \\ \lambda (X^2-1/4) &0 
   \end{pmatrix} 
   \mathrm{d}z,
\end{equation}
 where $a, b \in \mathbb R$, $X = \sqrt{ b^2 -a ^2 + a b (\lambda^{-2} -
 \lambda^2)}$ with  $b^2 -a ^2 + a b(e^{q} - e^{-q}) =1/4$. We assume 
 that 
 \begin{equation}\label{eq:trinoidcondition}  
  b^2\not\equiv a^2 \;(\mbox{mod}\;\boldsymbol S) \;\;\mbox{and}\;\;
  -1+ 4 \sin^2 (\pi X) \neq 0 \;\;\mbox{on}\;\;\lambda \in S^1,
 \end{equation}
 where $\boldsymbol S = \{y \in \mathbb Z \;|\; \mbox{
 $x^2 = y$ for some $x \in \mathbb Z$}\}$.
 It is known that one can construct CMC immersions from $\C P^1
 \setminus \{0, 1, \infty\}$ into $\mathbb E^3$ from potentials $\eta$
 which have the form \eqref{eq:trinoidspot}. However, in this case $a$
 and $b$ need to be chosen differently, see \cite{DW-n-noids} and
 \cite{KKRS}. It turns out that 
 the conditions \eqref{eq:trinoidcondition} are equivalent to that
 monodromy matrices can be simultaneously  conjugated into 
 $\Lambda {\rm SL}_2 \C_{\sigma, \tau_4} 
 \subset \Lambda {\rm SL}_2 \C_{\sigma}$.
 In our case, following an analogous procedure, 
 it can be shown that the resulting
 CMC immersion with $H = \tanh q$ in $\mathbb H^3$  given by the holomorphic
 potential in \eqref{eq:trinoidspot} is well-defined on $\C P^1
 \setminus \{0, 1, \infty\}$. The detailed computation has been discussed in 
 \cite{Kobayashi:tri}.

\appendix
\section{Gauss-Codazzi equations}\label{sc:Gauss-Codazzi}
 Let $f:\mathbb{D}\to \mathbb{H}^3\subset \mathbb{E}^{1,3}$ be a
 conformal immersion of a simply connected region $\mathbb{D}$ in
 $\C$ with unit normal $n$. Define a frame field
 $s=(f,f_z,f_{\bar z},n)$ along $f$. Then  we have
$$
\mathcal{U}=s^{-1}s_{z}=
\left(
\begin{array}{cccc}
0 & 0 & \frac{1}{2}e^u & 0 \\
1 & u_z & 0 & -H \\
0 & 0 & 0 & -2Qe^{-u} \\
0 & Q & \frac{H}{2}e^u & 0 
\end{array}
\right),
$$
$$
\mathcal{V}=s^{-1}s_{\bar z}=
\left(
\begin{array}{cccc}
0 & \frac{1}{2}e^u & 0 & 0 \\
0 & 0 & 0 & -2{\bar Q}e^{-u} \\
1 & 0 & u_{\bar z} & -H \\
0 & \frac{H}{2}e^u & {\bar Q} & 0 
\end{array}
\right).
$$
 The Gauss-Codazzi equations
$$
\mathcal{V}_{z}-\mathcal{U}_{\bar z}+
[\mathcal{U},\mathcal{V}]=0
$$
 of $f$ are equivalent with
$$
 u_{z{\bar z}}+\frac{1}{2}(H^2-1)e^u -2|Q|^2 e^{-u}=0, \;\;\;
 Q_{\bar z}=\frac{1}{2}H_{z}e^u.
$$
\section{Surfaces in $3$-dimensional space forms}\label{sc:Lawsoncorrespond}
 Let us denote by $\mathscr{M}^{3}(c)$ the $3$-dimensional simply
 connected Riemannian space form.
 
 Let $f:M\to \mathscr{M}^3(c)$ be a conformal immersion of a  Riemann
 surface $M$ into the space form $\mathscr{M}^{3}(c)$. Denote by $n$ the
 unit normal vector field to $f$. Take a simply connected local complex
 coordinate region $(\mathbb{D},z)$. Express the induced metric
 $\mathrm{I}$ as
\begin{equation}
\mathrm{I}=e^{u}
\mathrm{d}z\mathrm{d}\bar{z}.
\end{equation}
 Then the Gauss-Codazzi equations of $f$ are given by
\begin{equation*}
 u_{z{\bar z}}+\frac{1}{2}(H^2+c)e^u-2|Q|^2 e^{-u}=0, \;\;\;
 Q_{\bar z}=\frac{1}{2}H_{z}e^u.
\end{equation*}
 Here $H$ is the mean curvature of $f$. The Hopf differential of $f$ is
 defined by $Q\>\mathrm{d}z^2$ with $Q=\langle f_{zz},n\rangle$.

 The Gauss-Codazzi equations imply the following fact:
\begin{proposition}
 Let $f:\mathbb{D}\to \mathscr{M}^{3}(c)$ be a simply connected surface
 of constant mean curvature $H$. Take a pair $(\tilde{c}, \tilde{H})$ of
 real numbers such that $H^2+c={\tilde H}^2+{\tilde c}$. Then there
 exists a conformal immersion $\tilde{f}: \mathbb{D}\to
 \mathscr{M}^{3}(\tilde{c})$ with constant mean curvature $\tilde{H}$
 whose induced metric is the original metric of $(\mathbb{D},f)$.
\end{proposition}
 The correspondence $f\mapsto \tilde{f}$ is frequently called the
 \textit{Lawson correspondence}.  In particular, surfaces of constant
 mean curvature $H=\pm 1$ in $\mathbb H^3$ correspond to minimal
 surfaces in $\mathbb{E}^3$. A.~Fujioka gave a generalization of the Lawson
 correspondence in \cite{F1}. Moreover, he proved that the Lawson 
 correspondence and the dressing action of loop groups are equivariant
 \cite{F2}.

 Via the Lawson correspondences, one can construct constant mean
 curvature surfaces in $\mathbb S^3$ or $\mathbb H^3$ by the 
 generalized Weierstrass type representation (DPW method) 
 of CMC surfaces in $\mathbb{E}^3$.  For instance, in \cite{KKRS},
 some CMC surfaces in $\mathbb S^3$ or $\mathbb H^3$ are
 constructed.

 However, surfaces of constant mean curvature $H$ in $\mathbb H^3$ such
 that $0 \leq H <1$ have no corresponding surfaces in $\mathbb{E}^3$ or
 $\mathbb S^3$.  This has motivated us to establish a DPW method
 for CMC surfaces in  $\mathbb H^3$ without any restrictions on the
 range of the mean curvature.

\section{Semi-Riemannian submersion}\label{sc:submersion}
 Let $(P,g_P)$, $(N,g_N)$ be semi-Riemannian manifolds. 
 According to O'Neill \cite{Oneill:Semi-Riemannian}, a submersion $\pi :
 P \to N$  is said to be a \textit{semi-Riemannian submersion} if 
\begin{enumerate}
\item[(i)] The fibers $\pi^{-1} \{x\}, x \in N$ are semi-Riemannian
 submanifolds of $P$,
\item[(ii)] $\mathrm{d}\pi$ preserves scalar products of vectors normal to
 fibers.
\end{enumerate}
 The kernel $\mathrm{Ker}(\pi_{*p})$ is called the \textit{vertical
 subspace} at $p$ and denoted by $\mathcal {V}_p$. The orthogonal 
 complement $\mathcal H_p = \mathcal V_p^{\perp}$ is called the 
 \textit{horizontal subspace} at $p$. By definition, $\mathrm{d} \pi :  
 T_p P \to T_{\pi(p)} N$ is a linear isometry. 

 Let $\pi: P \to N$ a semi-Riemannian submersion and $\psi : M \to P$ 
 a smooth map from a Riemannian manifold $M$. Then $\psi$ is said to be
 \textit{horizontal} with respect to $\pi$ ($\pi$-horizontal in short)
 if $\mathrm{d} \psi (TM) \subset \mathcal H$ . Then the tension field
 $\tau (\psi)$ of a $\pi$-horizontal map satisfies the following
 formula:
 \begin{equation*} 
 \mathrm{d} \pi (\tau (\psi)) = \tau (\pi \circ \psi).
 \end{equation*}
Hence we have the following result (\textit{cf}. \cite[Proposition 2.36,
p.~203]{Ur}).
\begin{proposition}\label{submersionLemma}
 Let $M$ be a Riemannian manifold and $\psi:M\to P$ 
 a harmonic map which is $\pi$-horizontal.
 Then $\pi\circ \psi :M\to N$ is also harmonic.
\end{proposition}
 This Proposition will be used in appendix \ref{AppendixC}.

\section{Gauss maps}\label{AppendixC}
 In the surface geometry of hyperbolic $3$-space, several notions of
 ``Gauss map'' have been introduced. From our point of view (integrable
 system approach), the Gauss map $F=(f,n)$ is the \textit{best}
 one. In this appendix, we collect some other ``Gauss maps'' and
 discuss how they are not compatible with the generalized Weierstrass
 type representation.
 
\subsection{Obata's Gauss map}\label{App:Obata}
 Let us denote by $\mathcal{S}$ the space of all oriented totally
 geodesic surfaces in $\mathbb{H}^3$. This space is identified with the
 Grassmann manifold
 $\mathrm{Gr}_{1,2}(\mathbb{E}^{1,3})$ of
 $3$-dimensional oriented  \textit{timelike} subspaces (\textit{i.e.,}
 linear subspaces with signature $(-,+,+)$) in $\mathbb E^{1,3}$.  In fact,
 take such a timelike
 subspace $W$, then its intersection $\gamma:=W\cap \mathbb{H}^3$ is a
 totally geodesic surface in $\mathbb{H}^3$. Under the identification
 $\gamma$ with $W$, the space $\mathcal{S}$ is identified with
 $\mathrm{Gr}_{1,2}(\mathbb{E}^{1,3})$. Next, via taking the orthogonal
 complements, $\mathrm{Gr}_{1,2}(\mathbb{E}^{1,3})$ is identified with
 the space $\mathrm{Gr}_{1}^{+}(\mathbb{E}^{1,3})$ of all oriented
 \textit{spacelike} lines. 
$$
\perp: \mathrm{Gr}_{1,2}(\mathbb{E}^{1,3})
\to \mathrm{Gr}_{1}^{+}(\mathbb{E}^{1,3});
\ \
W\longmapsto W^{\perp}.
$$
 Obviously, the Grassmann manifold
 $\mathrm{Gr}_{1}^{+}(\mathbb{E}^{1,3})$ is the \textit{de Sitter
 $3$-space}
$$
\mathbb{S}^{1, 2}=\{(x_0,x_1,x_2,x_3)\in \mathbb{E}^{1,3}
\
\vert
\ -x_{0}^{2}+x_{1}^{2}+x_{2}^{2}+x_{3}^{2}=1\}. 
$$
 Hence, the space $\mathcal{S}$ is identified with de Sitter $3$-space.

 Let $f:M\to \mathbb{H}^3\subset \mathbb{E}^{1,3}$ be a conformal
 immersion. To $f$, one can associate a map $\gamma:M\to \mathcal{S}$
 as follows: Define a $3$-dimensional linear subspace $W(p)\subset
 \mathbb{E}^{1,3}$ by
$$
W(p):=\mathrm{d}f_{p}(T_{p}M)\oplus \mathbb{R}f(p).
$$
 Take the intersection $\gamma(p):=W(p)\cap \mathbb{H}^{3}$ of $W(p)$ 
 with the hyperbolic $3$-space. Then $\gamma(p)$ is a totally geodesic
 subspace of $\mathbb{H}^3$. The resulting map
 $\gamma:M\to\mathcal{S}=\mathrm{Gr}_{1,2}(\mathbb{E}^{1,3})$ is called
 \textit{Obata's Gauss map} \cite{Obata}. Under the identification
 $\gamma\leftrightarrow \gamma^\perp$, Obata's Gauss map can be
 regarded as a map into the de Sitter 3-space $\mathbb{S}^{1, 2}$.

 Now let $F=(f,n):M\to \mathrm{U}\mathbb{H}^3$ be the Gauss map of
 $f$. Then we can easily check that (see section \ref{sc:fibering})
$$
\pi_{3}\circ F=\gamma^{\perp}, \ \
\pi_{3}:\mathrm{U}\mathbb{H}^3
\to \mathbb{S}^{1, 2}
$$
 Thus Obata's Gauss map is identified with $n:M\to \mathbb{S}^{1, 2}$.
 The harmonicity of Obata's Gauss map is characterized as follows \cite{Ishihara}, \cite{Obata}:
\begin{proposition}
 Obata's Gauss map $n:M\to \mathbb{S}^{1, 2}$ is harmonic if and only if
 $f$ is minimal. The map $n:M\to \mathbb{S}^{1, 2}$ is a singular 
 spacelike surface with zero mean curvature. 
\end{proposition}
 Note that $n:M\to \mathbb{S}^{1, 2}$ is referred to as the \textit{polar
 variety} of $(M,f)$. Compare with the $\mathbb{S}^3$-case \cite{Lawson}.

\subsection{Normal Gauss map}
\subsubsection{}
 The hyperbolic $3$-space $\mathbb{H}^3$ is represented by
 $\mathbb{H}^3=\mathrm{SL}_{2}\C/\mathrm{SU}_2$ as a Riemannian
 symmetric space in section \ref{sc:unit tangent sphere bundle}. On the
 other hand, $\mathrm{SL}_{2}\C$ admits an Iwasawa decomposition:
$$
 \mathrm{SL}_{2}\C=\mathscr{S}\cdot \mathrm{SU}_2,
$$
 where $\mathscr{S}$ is a $3$-dimensional solvable Lie subgroup of
 $\mathrm{SL}_{2}\C$. Explicitly,
$$
\mathscr{S}=\left\{
\left(
\begin{array}{cc}
\sqrt{u_3} & (u_{1}+iu_{2})/\sqrt{u_3}\\
0 & 1/\sqrt{u_3}
\end{array}
\right)
\
\biggr
\vert
\
u_1,u_2\in \mathbb{R},
u_3>0
\right\}.
$$
 Hyperbolic $3$-space $\mathbb{H}^3$ is identified with the solvable Lie 
 group $\mathscr{S}$ equipped with a special left invariant metric.
 In fact, $\mathscr{S}$ and this $\mathbb H^3$ can be identified with
 the upper half space model
$$
(\{(u_1,u_2,u_3)
\
\vert
\ u_3>0\},\
(\mathrm{d}u_1^2+\mathrm{d}u_2^2+\mathrm{d}u_3^2)/u_{3}^2\>).
$$
 The mapping 
$$
\phi:\mathscr{S}\to \mathbb{H}^{3}\subset \mathrm{Her}_2 \mathbb C;
\ \
\phi(s)=ss^{*}
$$
 is an isometry between $\mathscr{S}$ and $\mathbb{H}^3$.
 Hence $\mathbb{H}^3$ can be represented as
$$
\mathbb{H}^{3}=\{ss^{*}\
\vert
\ s\in \mathscr{S}\}
\subset
\mathrm{Her}_2 \C.
$$
 Denote by $\mathfrak{s}$ the Lie algebra of $\mathscr{S}$, then we have
 the following diagram (briefly explained below):
\begin{equation}
\begin{CD}
T_{ss^*}\mathbb{H}^3 @>
(\mathrm{d}\phi_{s})^{-1}
>> T_{s}\mathscr{S}\\
@V \mathrm{dA}_{s}^{-1} VV  
@VV \mathrm{d}L_{s}^{-1}
V
\\
\mathfrak{m}=T_{\Vec{e}_0}\mathbb{H}^3
@>
(\mathrm{d}\phi_{\Vec{e}_0})^{-1}
>>
T_{\Vec{e}_0}\mathscr{S}=\mathfrak{s}
\end{CD}
\end{equation}
 Take a tangent vector $X_{ss^*}$ of $\mathbb{H}^3$ at $ss^*$,
 then the corresponding tangent vector of $\mathscr{S}$ at $s$ is
 $\tilde{X}_{s}=(\mathrm{d}\phi_{s})^{-1}X_{ss^*}$. One can check that
\begin{equation}\label{C.1}
\mathrm{d}\phi_{\Vec{e}_0}(\mathrm{d}L_{s}^{-1}\tilde{X}_{s})=
\mathrm{dA}_{s}^{-1}(X_{ss^*}),
\end{equation}
 where $\mathrm{A}$ denotes the action of $\mathrm{SL}_{2}\C$ on
 $\mathrm{Her}_2 \C$. By using this \textit{solvable Lie group model}
 $\mathscr{S}$ of $\mathbb{H}^3$, we will introduce below the notion of
 a ``normal Gauss map'' for surfaces in  $\mathbb{H}^3$.

\subsubsection{}
 In this section, we recall some fundamental facts about the tangent bundle
 $TG$ of $G$.
 
 Let $G$ be a real Lie group with multiplication map
$$
\mu:G\times G\to G,\ \ \
\mu(a, b)=ab.
$$
 Then the differential ${\mathrm d}\mu$ of $\mu$ defines a
 multiplication on the tangent bundle $TG$:
\begin{equation}\label{TG}
{\mathrm d}\mu:TG \times TG \to TG;
\
{\mathrm d}\mu
((a;A_a),(b;B_b))=
(ab;\mathrm{d}L_{a}B_{b}+\mathrm{d}R_{b}A_{a}),
\end{equation}
 where $L_{a}$ denotes left multiplication by $a$ in $G$ and 
 $R_{b}$ denotes right multiplication by $b$ in $G$.
 With respect to ${\mathrm d}\mu$, the tangent bundle $TG$ becomes a Lie
 group which is called the \textit{tangent group} of $G$.

 Set $\mathfrak g = T_{e}G$, where $e$ denotes the identity element of
 $G$. Then the tangent bundle $TG$ of $G$ is identified via the left
 Maurer-Cartan form $\vartheta=\vartheta^{L}$ of $G$ with $G \times
 \mathfrak{g}$:
$$
(a;A_{a})\longmapsto (a,\vartheta^{L}(A_{a}))=(a,\mathrm{d}L_{a}^{-1}A_{a}).
$$
 As usual we identify  $\mathfrak{g} = T_{e}G$ with the Lie algebra of
 all smooth ``left invariant'' vector fields on $G$. Under this
 identification, the group structure of $TG$ is transferred to $G \times
 \mathfrak g$ as:
$$
(a,A)(b,B)=(ab, B+\mathrm{Ad}(b^{-1})A).
$$
 Hereafter, we denote this semi-direct product group by $G\ltimes
 \mathfrak{g}$.

 The Lie group $G$ is imbedded in $TG$ as set of all zero sections
 $\{(g,0)\ | \ g \in G\}$. The Lie algebra $\mathfrak{g}$ is identified
 with the normal subgroup $\{(e,A) \ | \ A \in \mathfrak{g}\}$.
\begin{remark}
{\rm 
 Analogously, by using the right Maurer-Cartan form $\vartheta^{R}$, 
 $TG$ can be  identified with $G \rtimes \mathfrak{g}$:
$$
(a;A_a)\longmapsto
(a,\vartheta^{R}(A_a))=
(a,\mathrm{d}R_{a}^{-1}(A_a))\in G \times \mathfrak{g},
$$
$$
(a,A)\longmapsto
(a,\mathrm{d}R_{a}A)\in TG. 
$$
 Under this identification, the multiplication law of $TG$ is transferred
 to $G \times \mathfrak g$ as
$$
 (a,A)(b,B)=(ab,A+\mathrm{Ad}(a)B).
$$
 Here we consider $\mathfrak{g}$ as the Lie algebra of all smooth
 ``right invariant'' vector fields on $G$. We denote this semi-direct
 product group by $G\rtimes \mathfrak{g}$. The identification $TG=G
 \rtimes \mathfrak{g}$ is used by V.~Balan and the first named author 
 \cite{BD1}.  
}
\end{remark}

\subsubsection{}
 Let $G=(G,\langle\cdot,\cdot\rangle)$ be a $3$-dimensional real Lie group
 equipped with a left invariant Riemannian metric. Let $f:M\to G$ be a
 conformally immersed surface with unit normal $n$.
 Then the unit normal vector field $n$ is regarded as a map $F=(f;n):M
 \to TG$. Under the identification $TG=G\ltimes \mathfrak{g}$, $F$
 induces a map $\tilde F=(f;\kappa):M\to G\ltimes
 \mathfrak{g}$. The map
 $\kappa=\mathrm{d}L_{f}^{-1}\>n:M\to \mathbb{S}^2\subset \mathfrak{g}$
 is called the \textit{normal Gauss map} of $f$. Here $\mathbb{S}^2$
 is the unit sphere in $\mathfrak{g}$ centered at the origin.

\subsection{}
 Now we apply these fundamental observations to
 $G=\mathscr{S} \cong \mathbb{H}^3$.
 Denote by $\mathfrak{s}$ the Lie algebra of the solvable Lie group
 $\mathscr{S}$. Then the tangent group of $\mathbb{H}^3$ is
 $T\mathbb{H}^3=\mathscr{S}\ltimes \mathfrak{s}$. The unit tangent
 sphere bundle $\mathrm{U}\mathbb{H}^3$ is identified with
$$
\mathscr{S}\times \mathbb{S}^{2}\subset \mathscr{S}\ltimes \mathfrak{s}.
$$
 Let $f:\mathbb{D}\to \mathbb{H}^3=\mathscr{S}$ be a surface with unit
 normal $n$. Then the Gauss map $F$ of $f$ is given by
 $F=(f,f^{-1}n)\in
 \mathscr{S}\times \mathbb{S}^2$. Now we compute the normal Gauss map
 $\kappa$ in terms of the Sym formula:
$$
f=\tilde \Phi\tilde \Phi^{*},\ \ n=\tilde \Phi \> \Vec{e}_{1}\tilde \Phi^{*},
$$
 where $\tilde \Phi$ is defined in \eqref{eq:tildephi}.
 According to the Iwasawa decomposition
 $\mathrm{SL}_{2}\C=\mathscr{S}\cdot \mathrm{SU}_2$,
 $\tilde \Phi$ is decomposed as $\tilde \Phi=su$, $s\in \mathscr{S}$ and
 $u\in \mathrm{SU}_2$. Hence $f=\tilde \Phi\tilde \Phi^{*}=ss^*$,
 since  $u$ is unitary. By using (\ref{C.1}), we have 
\begin{align*}
\mathrm{d}\phi_{\Vec{e}_0}(\kappa)
&=
\mathrm{d}\phi_{\Vec{e}_0}(\mathrm{d}L_{s}^{-1}\mathrm{d}\phi_{s}^{-1}n)=
\mathrm{dA}_{s}^{-1}(n)
\\
&= u\Vec{e}_{1}u^{*}=\mathrm{Ad}(u)\Vec{e}_{1}.
\end{align*}
 The Lie algebra $\mathfrak{sl}_{2}\C$ is decomposed as (see
 section \ref{sc:Minkowski}):
$$
\mathfrak{sl}_{2}\C=\mathfrak{su}_2\oplus 
\mathfrak{m},
\ \
\mathfrak{m}=i\mathfrak{su}_2=\mathbb{R}\Vec{e}_1
\oplus \mathbb{R}\Vec{e}_{2}
\oplus \mathbb{R}\Vec{e}_{3}.
$$
 The special unitary group $\mathrm{SU}_2$ acts on $\mathfrak{m}$ via
 the $\mathrm{Ad}$-action. The unit $2$-sphere $\mathbb{S}^2\subset
 \mathfrak{m}$ is the $\mathrm{Ad}(\mathrm{SU}_2)$-orbit of $\Vec{e}_1$. The
 isotropy
 subgroup of $\mathrm{SU}_2$ at $\Vec{e}_1$ is $H=\mathrm{U}_1$. Hence
 we have $\mathbb{S}^{2}=\mathrm{SU}_2/\mathrm{U}_1$. The natural
 projection $\pi_{\mathbb{S}^2}:\mathrm{SU}_2\to \mathbb{S}^2$
 is given explicitly by
 $\pi_{\mathbb{S}^2}(g)=g\Vec{e}_{1}g^*$.
 Hence the normal Gauss map is the (smooth) map
 $\kappa_{f}:={\rm Ad}(u)\Vec{e}_{1}:M \to \mathbb{S}^2$.

 G{\'o}es and Sim{\~o}es \cite{GS} obtained an integral representation
 formula for minimal surfaces in $\mathbb{H}^3$ and $\mathbb{H}^4$. The
 data for the formula due to G{\'o}es and Sim{\~o}es are smooth maps into
 $\mathbb{S}^2$ which are  solutions to certain second order elliptic
 partial differential equation. Independently, Kokubu \cite{Kok}
 obtained  same representation formula. Moreover, he showed that the data
 of the  representaion formula are normal Gauss maps which are  
\textit{harmonic} with respect to a certain singular Riemannian metric  
on $\mathbb{S}^2$. The singular metric (now referred as the
 \textit{Kokubu metric}) is not homogeneous. Thus one can not apply
 the DPW method to harmonic maps into $\mathbb{S}^2$ equipped with the Kokubu
 metric. The representaion formula due to G{\'o}es--Sim{\~o}es--Kokubu
 is generalized to the model space $\mathrm{Sol}$ of ``solvgeometry'' by
 Lee and the second named author \cite{IL}.

\begin{corollary}
 Let $\tilde \Phi$ be the extended
 frame defined in \eqref{eq:tildephi} for a minimal immersion 
 $f : \mathbb D \to \mathbb H^3$.
 Split $\tilde \Phi=su$, $s \in \mathscr S, u \in \mathrm{SU}_2$ according to
 the Iwasawa decomposition of $\mathrm{SL}_{2}\C = \mathscr S \cdot
 \mathrm{SU}_2$. Then
$$
 \kappa_{f}:\mathbb{D}\times C_{r}\to \mathbb{S}^{2}\subset \mathfrak{m}
$$ 
 is a loop of harmonic maps into the unit $2$-sphere equipped with the
 Kokubu metric. 
\end{corollary}

\subsection{Hyperbolic Gauss map}
 It is classically known  that, 
 by using the Poincar{\'e} ball
 model of $\mathbb{H}^3$, the space $\mathrm{Geo}(\mathbb{H}^3)$ of
 oriented geodesics in $\mathbb{H}^3$ is identified with the space 
$$
\{ (p,q) \in \mathbb{S}^{2}\times \mathbb{S}^{2}
\
\vert
\ p\not=q\}
$$
 of distinct point pairs in $\mathbb{S}^2=\partial\mathbb{H}^3$, where
 we identify $\mathbb S^2$ with the \textit{ideal boundary} 
 $\partial\mathbb{H}^3$
 of $\mathbb{H}^3$ \cite{Epstein}. An immersion $\ell:M\to
 \mathrm{Geo}(\mathbb{H}^3)$ of a $2$-manifold $M$ is said to be an
 \textit{oriented geodesic congruence} in $\mathbb{H}^3$.

 Let $f:M\to \mathbb{H}^3$ be a conformal immersion with unit normal
 $n$. Then there exists an oriented geodesic congruence
 $\mathcal{G}=(g_{L},g_{R}):M\to \mathrm{Geo}(\mathbb{H}^3)$ which
 satisfies the following condition:

 For any $p\in M$, the oriented geodesic $\gamma$ starting from $g_{L}(p)$
 and ending at $g_{R}(p)$ coincides
 with the oriented normal geodesic at $f(p)$ \cite{Epstein}.

 The map $\mathcal{G}$ is said to be the \textit{oriented normal
 geodesic congruence} of $(M,f,n)$. Each component map of $\mathcal{G}$
 is called a \textit{hyperbolic Gauss map}. 

 Let us recall the
 fibration $\pi_{1}:\mathrm{U}\mathbb{H}^{3}\to
 \mathrm{Geo}(\mathbb{H}^3)$.  By using the quadratic model
 of $\mathbb{H}^{3}\subset \mathbb{E}^{1,3}$ \cite{Spivak}, the space
 $\mathrm{Geo}(\mathbb{H}^3)$ is identified with the Grassmann manifold
 $\mathrm{Gr}_{1,1}(\mathbb{E}^{1,3})$. Under this model, the
 oriented normal geodesic congruence $\mathcal{G}$ is represented in the
 form \cite{Br}:
$$
\mathcal{G}=(f+n,f-n) = \pi_1 \circ F,
$$
 where  $F=(f,n)$ denotes the Gauss map. 
 By using the Sym formula, $\mathcal{G}$ is rewritten as  
$$
\mathcal{G}=(\tilde \Phi(\Vec{e}_0+\Vec{e}_1) \tilde \Phi^*,
\tilde \Phi(\Vec{e}_0-\Vec{e}_1)\tilde \Phi^{*}),
$$
 where $\tilde \Phi$ is defined in \eqref{eq:tildephi}.
 Proposition \ref{submersionLemma} implies that for a CMC surface
 $f:M\to \mathbb{H}^3$, its oriented normal geodesic congruence is a 
 harmonic map into $\mathrm{Geo}(\mathbb{H}^3)$. Moreover,  
 $\mathcal{G}$ is \textit{Lagrangian}, \textit{i.e.},
 $\mathcal{G}^{*}\Omega=0$ since $F$ is Legendrian. Here $\Omega$ is the 
 canonical symplectic form of $\mathrm{Geo} (\mathbb H^3)$. In fact, we have
\begin{equation*}
\mathcal{G}^* \Omega = (\pi_1 \circ F)^* \Omega = F^*(\pi_1^* \Omega)
 = F^* (\mathrm{d} \omega) = \mathrm{d} (F^* \omega) =0,
\end{equation*}
 since $F$ is Legendrian. Thus we have shown the following.
\begin{corollary}
 Let $\tilde \Phi$ be the extended frame defined in \eqref{eq:tildephi} 
 for a Legendre harmonic map $F:\mathbb{D} \to \mathrm{U}\mathbb{H}^3$. 
 Then 
$$
\mathcal{G}:\mathbb{D}\times C_{r}
\to \mathrm{Geo}(\mathbb{H}^3);
\ \
(\tilde \Phi (\Vec{e}_0+\Vec{e}_1)\tilde \Phi^*, 
 \tilde \Phi(\Vec{e}_0-\Vec{e}_1)\tilde \Phi^{*})
$$
 is a loop of Lagrangian harmonic maps.
\end{corollary}

\subsection{Generalized Gauss map}
 Let us denote by $\mathrm{Gr}_{2}^{+}(\mathbb{E}^{1,3})$ the Grassmann
 manifold of all oriented \textit{spacelike} planes in
 $\mathbb{E}^{1,3}$. Then one can see that the operation of taking 
 perpendicular subspaces:
$$
\perp:\mathrm{Gr}_{1,1}(\mathbb{E}^{1,3})\to
\mathrm{Gr}_{2}^{+}(\mathbb{E}^{1,3})
$$
 is an isometry.

 For a conformal immersion $f:M\to \mathbb{H}^3$, we can associate a map
 $\hat{F}:M\to \mathrm{Gr}_{2}^{+}(\mathbb{E}^{1,3})$ by 
$$
\hat{F}(z,\bar{z})=\mathrm{d}f_{(z,\bar{z})}(T_{(z,\bar{z})}M)\subset 
T_{f(z,\bar{z})}\mathbb{H}^3.
$$
 The map $\hat{F}$ is called the \textit{generalized Gauss map} of
 $f$.

 As we saw before, the space $\mathrm{Geo}(\mathbb{H}^3)$ of all
 oriented geodesics is identified with
 $\mathrm{Gr}_{1,1}(\mathbb{E}^{1,3})$. Thus the generalized Gauss map
 can be considered as a map into  $\mathrm{Geo}(\mathbb{H}^3)$. The
 resulting map into $\mathrm{Geo}(\mathbb{H}^3)$ is 
 the pair of \textit{hyperbolic Gauss maps}.

\section{Fronts}
\subsection{}
 Let us recall that a map $F: M \to \mathrm{U}\mathbb H^3$ is called 
 a Legendre map if $F^*\omega =0$, where $\omega$ has been defined 
 in section \ref{sc:unit tangent sphere bundle}.
 A smooth map $f:M \to \mathbb{H}^3$ is said to be a (\textit{wave})
 \textit{front} if there exists a Legendre immersion  $F:M \to
 \mathrm{U}\mathbb{H}^3$ such that $\pi_{2}\circ F=f$, where $\pi_2$ is 
 defined in section \ref{sc:fibering}. In other words,
 $f$ is a front if and only if there exists a map $n:M\to \mathbb{E}^{1,3}$ 
 such that 
\begin{enumerate}
\item[(i)] $\langle f,n\rangle=0$, $\langle n,n \rangle=1$,
\item[(ii)] $F:=(f,n)$ is an immersion into $\mathrm{U}\mathbb{H}^3$, 
\item[(iii)] $\langle \mathrm{d}f,n\rangle=0$. 
\end{enumerate}
 It is easy to see that an immersion $f:M\to \mathbb{H}^3$ is a front if and
 only if $M$ is orientable. Even if a front $f$ is not an immersion, the
 maps $n$ and $F$ are referred to as a \textit{unit normal vector field} and
 the \textit{Gauss map} of $f$, respectively.
\begin{definition}
{\rm
 A smooth map $f:M\to \mathbb{H}^3$ is said to be a CMC \textit{front} if
\begin{enumerate}
\item[(i)] $f$ is a front with Gauss map $F=(f,n)$ 
\item[(ii)] the Gauss map is a harmonic map. 
\end{enumerate}
}
\end{definition} 
 By the generalized Weierstrass type representation \cite{KKRS} and the
 results of this paper, we can construct CMC fronts from  prescribed potentials.

\subsection{}
 Let $f:M\to \mathbb{H}^3\subset \mathbb{E}^{1,3}$ be a front. Then
 the \textit{parallel front} $f_{r}$ in $\mathbb H^3$ at the distance
 $r$ of $f$ is defined by
$$
f_{r}=\cosh (r)\>f+\sinh (r)\>n.
$$
 On the other hand, the parallel front $\check{f}_r$ of
 $f$ at the distance $r$ in the de Sitter $3$-space $\mathbb{S}^{1, 2}$ is
 defined by
$$
\check{f}_{r}=\sinh (r)\>f+\cosh (r)\>n.
$$
 If $f$ is an immersion, then $f_r$ and $\check{f}_r$ are called
 \textit{parallel surfaces} of $f$.

 Now let $f$ be a CMC surface in $\mathbb H^3 $ with mean curvature
 $0 \leq H = \tanh q <1$ and Hopf differential $Q\mathrm{d}z^2$. Then the
 parallel surface $\check{f}_{q}$
 is a spacelike CMC surface in $\mathbb{S}^{1, 2}$ with Hopf differential
 $Q \mathrm{d}z^2$. Now let us denote by
 $\check{M}$ the Riemann surface which is obtained by reversing the
 orientation of $M$. Then $\check{f}_{q}:\check{M}\to \mathbb{S}^{1, 2}$ is
 a conformal spacelike CMC immersion with Hopf differential
 $\bar{Q}\mathrm{d} \bar z^2$.
\begin{corollary}
 Let $\Phi$ be the extended frame of some CMC immersion $f$ into
 $\mathbb H^3$ with mean curvature $0 \leq H = \tanh q <1$. Then 
$$
\check{f}_{q}=\Phi 
\left(
\begin{array}{cc}
e^{q/2} & 0\\
0 & -e^{-q/2}
\end{array}
\right)
\Phi^{*} 
$$ 
 is a loop of spacelike CMC immersion with mean curvature $H=\tanh q$
 in $\mathbb S^{1, 2}$.
\end{corollary}


\bibliographystyle{plain}
\def\cprime{$'$}

\end{document}